\documentclass[reqno,10pt]{amsart}
\usepackage{xcolor}

\usepackage{amssymb}
\usepackage{comment}
\usepackage[all,cmtip]{xy}
\newtheorem{MainTheorem}{Theorem}
\newtheorem*{UnnumberedCorollary}{Corollary}
\newtheorem{Proposition}{Proposition}[section]
\newtheorem{Definition}[Proposition]{Definition}
\newtheorem{Lemma}[Proposition]{Lemma}
\newtheorem{Theorem}[Proposition]{Theorem}
\newtheorem{Corollary}[Proposition]{Corollary}
\newtheorem{Remark}[Proposition]{Remark}

\DeclareMathOperator{\Val}{Val}

\DeclareMathOperator{\Gr}{Gr}
\DeclareMathOperator{\Fl}{Fl}
\DeclareMathOperator{\Res}{Res}
\DeclareMathOperator{\Jac}{Jac}
\renewcommand{\Re}{\mathrm {Re}}
\renewcommand{\Im}{\mathrm {Im}}
\DeclareMathOperator{\OO}{\mathrm {O}}
\DeclareMathOperator{\SO}{\mathrm {SO}}
\DeclareMathOperator{\vol}{vol}

\DeclareMathOperator{\Dens}{Dens}
\DeclareMathOperator{\Ad}{Ad}
\DeclareMathOperator{\Span}{Span}

\DeclareMathOperator{\arcosh}{arcosh}

\DeclareMathOperator{\Ker}{Ker}
\DeclareMathOperator{\Hom}{Hom}
\DeclareMathOperator{\Stab}{Stab}
\DeclareMathOperator{\Sym}{Sym}
\DeclareMathOperator{\Kl}{Kl}
\DeclareMathOperator{\Sc}{Sc}
\DeclareMathOperator{\Cr}{Cr}
\DeclareMathOperator{\WF}{WF}

\DeclareMathOperator{\GL}{GL}
\DeclareMathOperator{\sign}{sign}
\newcommand{\R}{\mathbb{R}}
\newcommand{\KS}{\mathrm{KS}}
\newcommand{\C}{\mathbb{C}}

\newcommand{\largewedge}{\mbox{\Large $\wedge$}}

\title{Valuation theory of indefinite orthogonal groups}

\author{Andreas Bernig}
\author{Dmitry Faifman}
 
\email{bernig@math.uni-frankfurt.de}
\email{dfaifman@math.utoronto.ca}
\date{}

\address{Institut f\"ur Mathematik, Goethe-Universit\"at Frankfurt,
Robert-Mayer-Str. 10, 60629 Frankfurt, Germany}
\address{Department of Mathematics, University of Toronto, Bahen Centre, 40 St. George St., Toronto ON M5S 2E4, Canada}

\thanks{A. B. was supported by DFG grants BE 2484/5-1 and BE 2484/5-2. D.F. was partially supported by a post-doctoral fellowship grant of the Rothschild Foundation and by an NSERC Discovery grant.\\ AMS 2010 {\it Mathematics subject
classification}: 
52B45; 
53C65
}

\begin{document}

\begin{abstract}
Let $\mathrm{SO}^+(p,q)$ denote the identity connected component of the real orthogonal group with signature $(p,q)$. We give a complete description of the spaces of continuous and generalized translation- and $\mathrm{SO}^+(p,q)$-invariant valuations, generalizing Hadwiger's classification of Euclidean isometry-invariant valuations. As a result of independent interest, we identify within the space of translation-invariant valuations the class of Klain-Schneider continuous valuations, which strictly contains all continuous translation-invariant valuations. The operations of pull-back and push-forward by a linear map extend naturally to this class.
\end{abstract}

\maketitle 
\tableofcontents

\section{Introduction and statement of main results}

\subsection{Background}

\subsubsection{Valuations}

Valuation theory is in many ways a generalization of measure theory. Roughly
speaking, one relaxes countable additivity to finite additivity, and at the same time replaces the sigma-algebra of measurable sets by a smaller family of geometrically nice sets, such as convex bodies or manifolds with corners. In many cases, one adds a different kind of analytic requirement, e.g. continuity w.r.t. the Hausdorff metric. This leads to a theory which unifies many seemingly different notions such as volume, surface area, Euler characteristic, function etc.

In his 3rd problem, Hilbert asked whether a definition of volume for polytopes using finite additivity only can be given. As Dehn proved shortly afterwards, such a definition is impossible, as there are many other functionals which are finitely but not countably additive, and vanish on polytopes of positive codimension. The ensemble of those discontinuous valuations on polytopes is known today as the Dehn invariants. 

The theory of valuations had a great impact on convex geometry. Almost all natural invariants in convex geometry can be interpreted as valuations. For example, the volume, the surface area, the mean width, mixed volumes, but also the affine surface area, projection and intersection bodies and the Steiner point are valuations \cite{schneider_book14}. Valuations are also fundamental for integral geometry, for example, Crofton's formula and Weyl's tube formula can be considered as statements about valuations.

A more recent line of research, initiated by Z\"ahle and Fu, is to study valuations on certain non-convex sets, such as sets of positive reach, manifolds (possibly with boundary or corners), subanalytic sets, or the so called WDC-sets \cite{fu94, pokorny_rataj, zaehle86}. In this theory, convexity arguments are replaced by tools from geometric measure theory. This opened the way for the investigation, initiated by Alesker, of valuations on manifolds, see \cite{alesker_val_man1, alesker_val_man2, alesker_survey07, alesker_val_man4, alesker_intgeo, alesker_bernig_convolution, alesker_bernig, alesker_val_man3, bernig_broecker07}. 

For applications in complex geometry, it is natural to work in the complex projective space. Recently, Alesker's theory of valuations on manifolds was used in the study of different concepts from algebraic geometry: Chern classes, Bezout's theorem and tube formulas \cite{bernig_fu_solanes}. The fundamental insight is that the Euler characteristic can be interpreted as a valuation on manifolds. This follows easily from Chern's proof of the Chern-Gauss-Bonnet theorem. Using the Euler characteristic, one can build many other interesting valuations on complex projective space. The recently found kinematic formulas on complex projective spaces can be seen as vast generalizations of Bezout's theorem. 
From a different perspective on algebraic geometry, valuation theory of polytopes on lattices appears naturally in the theory of toric varieties and in Ehrhart's theory of lattice point counting \cite{barvinok,gruber_book}. 

\subsubsection{Specific background}

Let us collect some definitions and results which are more specifically relevant for the present paper. 

Let $V$ be a finite-dimensional vector space and $\mathcal{K}(V)$ the set of compact convex bodies in $V$. A valuation (sometimes also called convex valuation) is a map $\mu:\mathcal{K}(V) \to A$ such that 
\begin{displaymath}
 \mu(K \cup L)+\mu(K \cap L)=\mu(K)+\mu(L)
\end{displaymath}
whenever $K,L,K \cup L \in \mathcal{K}(V)$. Here $A$ is an abelian semi-group. 

In this article, we will restrict to the cases $A=\R,\C$, but we invite the reader to look at the references \cite{abardia12, bernig_fu_solanes, hug_schneider_localtensor, ludwig_2005, schuster_wannerer, wannerer_area_measures, wannerer_unitary_module} for recent developments for other abelian semi-groups $A$. 

Examples of valuations on a Euclidean vector space are the intrinsic volumes \cite{klain_rota}, or the mixed volumes $K \mapsto V(K[i],L_{i+1},\ldots,L_n)$ for fixed $L_{i+1},\ldots,L_n \in \mathcal{K}(V)$. Both examples are translation-invariant in the obvious sense and continuous with respect to the topology induced by the Hausdorff metric. 

One of the most influential theorems in integral geometry is Hadwiger's theorem, stating that the vector space $\Val^{\mathrm{SO}(n)}$ of all rotation- and translation-invariant continuous valuations is spanned by the intrinsic volumes. Many theorems in integral geometry, like kinematic formulas, Kubota's formula, Steiner's formula etc. are easy consequences of Hadwiger's theorem, compare \cite{klain_rota} for a nice introduction to this topic. 

The theory of continuous and translation-invariant valuations is fundamental for an understanding of all valuations (on affine spaces and even on manifolds). A breakthrough was achieved by Alesker who confirmed McMullen's conjecture that mixed volumes span a dense subspace in the space of continuous, translation-invariant valuations. Based on this theorem, several algebraic structures on a certain dense subspace of smooth valuations (product, convolution, Alesker-Fourier transform) were constructed, and these structures were used in an algebraic treatment of integral geometric questions. We refer to \cite{bernig_aig10, fu_barcelona} for surveys on this topic. 

The Alesker product of smooth valuations satisfies a version of Poincar\'e duality which can be used to introduce the large class of generalized valuations \cite{alesker_val_man4}. Generalized valuations on manifolds are important for our understanding of kinematic formulas, see \cite{alesker_bernig_convolution, alesker_bernig, bernig_fu_solanes}. Recently we showed that generalized translation-invariant valuations form a partial algebra which contains McMullen's polytope algebra \cite{bernig_faifman}. Such valuations will be essential in the present paper. 

Alesker \cite{alesker_survey07} showed that, given a compact group $G\subset \mathrm{GL}(V)$, the space $\Val^G$ of $G$-invariant and translation-invariant continuous valuations is finite-dimensional if and only if $G$ acts transitively on the unit sphere. The connected groups acting effectively and transitively on some unit sphere were classified by Montgomery-Samelson and Borel \cite{borel49, montgomery_samelson43}. Besides the euclidean rotation group, there are complex and quaternionic versions of rotation groups, $\mathrm{U}(n), \mathrm{SU}(n), \mathrm{Sp}(n), \mathrm{Sp}(n) \cdot \mathrm{U}(1), \mathrm{Sp}(n) \cdot \mathrm{Sp}(1)$, as well as three exceptional cases $\mathrm{G}_2, \mathrm{Spin}(7), \mathrm{Spin}(9)$. 

As with Hadwiger's theorem, finite-dimensionality of $\Val^G$ implies the existence of integral geometric formulas, like kinematic formulas. The program, initiated by Alesker, to obtain Hadwiger-type theorems and kinematic formulas for these groups, has seen a lot of progress in recent years. Algebraic operations on valuations play an essential role in the explicit computation of such formulas. We refer to \cite{alesker03_un, alesker_su2_04, alesker08, bernig_g2, bernig_fu_hig, bernig_fu_solanes, bernig_solanes, bernig_voide, fu06} for more information. 

Skipping the compactness assumption for the group $G$, one usually has to weaken the continuity assumption to obtain interesting characterization theorems. As an example, Ludwig and Reitzner \cite{ludwig_reitzner99, ludwig_reitzner10} showed that the space of translation-invariant, $\mathrm{SL}(n)$-invariant and semi-continuous valuations is spanned by the affine surface area, the Euler characteristic and the volume. 

In \cite{alesker_faifman}, Alesker and the second named author studied valuations invariant under the indefinite orthogonal group $\mathrm{SO}^+(n-1,1)$, also known as the Lorentz group. They showed that there are relatively few invariant continuous valuations:  apart from the Euler characteristic and the volume, they only appear in degree of homogeneity $(n-1)$, and the dimension of the corresponding space is $3$ (provided $n \geq 3$),  or $2$ if only even valuations are considered. Replacing the space of continuous valuations by the slightly larger space of generalized valuations, their classification becomes more similar to Hadwiger's characterization: for each degree of homogeneity between $1$ and $(n-1)$, the space of even invariant generalized valuations is $2$-dimensional. These valuations are constructed by some Crofton formulas with generalized Crofton measures. The case of odd generalized translation-invariant valuations was not treated.   

\subsection{Results of the present paper}
In this paper, we give a non-compact version of Hadwiger's theorem. More precisely, we characterize continuous and generalized translation-invariant valuations invariant under some indefinite orthogonal group $\mathrm{SO}^+(p,q)$. This group is the connected component of the identity of the group $\mathrm{O}(p,q)$ of linear maps leaving invariant some non-degenerate quadratic form $Q$ on $V$ of signature $(p,q)$. We will often write $\mathbb R^{p,q}$ instead of $V$ to emphasize that we work in the category of linear spaces equipped with a quadratic form. For computations we typically work with the standard $(p,q)$ form on $\mathbb R^n$.

Our first main theorem describes the dimension of the spaces of $k$-homogeneous, $\mathrm{SO}^+(p,q)$-invariant translation-invariant continuous valuations. It generalizes the classical Hadwiger theorem (which is the case $\min(p,q)=0$) and the theorem by Alesker and the second named author on the Lorentz group \cite{alesker_faifman} (which is the case $\min(p,q)=1$). 

\begin{MainTheorem} \label{mainthm_continuous_case} 
For $k \in \{0,n\}$, 
\begin{displaymath}
 \dim \Val_k(V)^{\mathrm{SO}^+(p,q)}=1;
\end{displaymath}
for $1 \leq k \leq n-2$,
\begin{displaymath}
\dim \Val_k(V)^{\mathrm{SO}^+(p,q)}=
\begin{cases} 1 & \min(p,q)=0,\\
0 & \min(p,q) \geq 1;
\end{cases} 
\end{displaymath}
and for $k=n-1$, 
\begin{displaymath}
\dim \Val_k(V)^{\mathrm{SO}^+(p,q)}=
\begin{cases} 1 & \min(p,q)=0,\\
4 &  p=q=1,\\
3 & \min(p,q)=1,n\geq 3,\\
2 & \min(p,q) \geq 2.
\end{cases} 
\end{displaymath}
In all cases with $\min(p,q)\geq1$, $\dim \Val_{n-1}^+(V)^{\mathrm{SO}^+(p,q)}=\dim \Val_{n-1}(V)^{\mathrm O(p,q)}=2$.
In particular,  
\begin{displaymath}
\dim \Val_{n-1}^{-}(V)^{\mathrm{SO}^+(p,q)}=
\begin{cases} 0 & \min(p,q)=0,\\
2 &  p=q=1,\\
1 & \min(p,q)=1, n \geq 3,\\
0 & \min(p,q) \geq 2.
\end{cases} 
\end{displaymath}
\end{MainTheorem}

Our second main theorem describes the spaces of $k$-homogeneous, $\mathrm{SO}^+(p,q)$-invariant generalized translation-invariant valuations,  which is denoted by \linebreak $\Val_k^{-\infty}(V)^{\mathrm{SO}^+(p,q)}$. Here again, the result in the case $\min(p,q)=0$ follows more or less directly from Hadwiger's theorem. The case $\min(p,q)=1$ was treated in \cite{alesker_faifman}, but only for even valuations. Here we give an independent treatment and complete the previous results by working out the dimensions in the odd case too. In the case $\min(p,q) \geq 2$ no previous results seem to be known.   

\begin{MainTheorem} \label{thm_dimensions}
For $k \in \{0,n\}$, 
\begin{displaymath}
 \dim \Val_k^{-\infty}(V)^{\mathrm{SO}^+(p,q)}=1;
\end{displaymath}
and for $1 \leq k \leq n-1$,
 \begin{displaymath}
\dim \Val_k^{-\infty}(V)^{\mathrm{SO}^+(p,q)}=
\begin{cases} 1 & \min(p,q)=0,\\
4 &  p=q=1,\\
3 & \min(p,q)=1, n \geq 3,\\
2 & \min(p,q) \geq 2.
\end{cases} 
\end{displaymath}
In all cases with $\min(p,q)\geq1$ and $1\leq k\leq n-1$, $\dim \Val_{k}^{+,-\infty}(V)^{\mathrm{SO}^+(p,q)}=\dim \Val_{k}^{-\infty}(V)^{\mathrm O(p,q)}=2$.
In particular for $1 \leq k \leq n-1$, 
 \begin{displaymath}
 \dim \Val_k^{-,-\infty}(V)^{\mathrm{SO}^+(p,q)}=
 \begin{cases} 0 & \min(p,q)=0,\\
 2 &  p=q=1,\\
 1 & \min(p,q)=1, n \geq 3,\\
 0 & \min(p,q) \geq 2.
 \end{cases} 
 \end{displaymath}
\end{MainTheorem}

In order to describe our third main theorem, we have to introduce some notation. 

Let $V$ be an $n$-dimensional vector space. For $0 \leq k \leq n$ let $\Gr_k(V)$ be the Grassmann manifold of $k$-planes in $V$. If $\phi \in \Val_k^+(V)$ is an even continuous translation-invariant valuation, then the restriction of $\phi$ to a $k$-plane $E$ is a multiple of the Lebesgue measure (here we use for simplicity an auxiliary Euclidean structure). Putting $\Kl_\phi(E)$ for the proportionality coefficient, we obtain a function $\Kl_\phi \in C(\Gr_k(V))$, which is called the Klain function of $\phi$. By a fundamental result due to Klain \cite{klain00}, the map $\Kl:\Val^+_k(V) \to C(\Gr_k(V))$ is injective. If $\phi$ is not continuous, but a generalized even translation-invariant valuation, then $\Kl_\phi$ may still be defined as a generalized function on $\Gr_k(V)$ \cite{alesker_faifman}. 

Assume for simplicity that $V=\R^n=\R^{p,q}$ is equipped with the standard Euclidean form, as well as the standard $(p,q)$-form $Q=x_1^2+\dots+x_p^2-x_{p+1}^2-\dots-x_n^2$. In order to describe the Klain functions of the $\mathrm{SO}^+(p,q)$-invariant and translation-invariant even valuations, we need the orbit structure of $\Gr_k(\R^{p,q})$ under the larger group $\mathrm{O}(p,q)$. These orbits are given by the sets $X^k_{a,b}$ consisting of all $k$-planes such that the restriction of the quadratic form has signature $(a,b)$. Here $a,b$ range over all integers with $\max(0,k-q) \leq a \leq p, \max(0,k-p) \leq b \leq q, a+b \leq k$. 

The open orbits are those where $a+b=k$. If $E$ belongs to such an orbit, we fix a Euclidean-orthonormal basis $v_1,\ldots,v_k$ of $E$ and set 

\begin{displaymath}
 \kappa_a(E):=\begin{cases} \left|\det (Q(v_i,v_j))_{i,j=1}^k \right|^\frac12 &  \text{if }E \in X^k_{a,k-a},\\0 & \text{otherwise.}\end{cases}
\end{displaymath}

Our third main theorem describes the even $\mathrm{SO}^+(p,q)$--invariant generalized translation-invariant valuations in terms of their Klain functions. 

\begin{MainTheorem} \label{main_thm_klain}
The Klain function of a valuation $\phi \in \Val_k^{+,-\infty}(\R^{p,q})^{\mathrm{SO}^+(p,q)}$ is a linear combination of $\{\kappa_a\}$. The function $\sum_{a=\max(0,k-q)}^{\min(k,p)} c_a \kappa_a$ is in the image of the Klain map if and only if 
\begin{displaymath} 
 c_{a+1}+c_{a-1}=0, \quad \max(0,k-q)<a<\min(k,p).
\end{displaymath}
\end{MainTheorem}

We also introduce the class of Klain-Schneider continuous (or KS-continuous) valuations, denoted $\Val^{\KS}(V)$, and study some of its properties. This class plays an important role in the present study, but is interesting in its own right. The KS-continuous class is comprised of those generalized valuations that have continuous Schneider sections, with the topology of uniform convergence on the Schneider sections. Non-formally, a $k$-homogeneous generalized valuations is KS-continuous if it can be naturally evaluated on a $(k+1)$-dimensional convex body, and the latter value depends continuously on the body. The main properties are summarized in the following theorem.

\begin{MainTheorem} \label{mainthm_klaincont}
The pull-back and push-forward by a linear map $T:U \to V$, originally acting between the corresponding spaces of continuous valuations, extend by continuity to maps between the corresponding spaces of KS-continuous valuations. Moreover, the class of even KS-continuous valuations is invariant under the Alesker-Fourier transform.
\end{MainTheorem}

Thus this class enjoys many nice properties, characteristic of both smooth valuations (such as invariance under the Alesker-Fourier transform in the even case) and continuous valuations (such as functoriality in the linear category), while being quite large at the same time. 

We prove that all the invariant valuations are in fact KS-continuous.

\begin{UnnumberedCorollary}
All elements in  $\Val_k^{-\infty}(\R^{p,q})^{\mathrm{SO}^+(p,q)}$ are KS-continuous. For $p' \leq p, q' \leq q$, $0\leq k\leq p'+q'-1$, let $i: \mathbb R^{p', q'} \hookrightarrow \mathbb R^{p,q}$ be a fixed isometric inclusion. The restriction map
\begin{displaymath}
i^*: \Val_k^{+,-\infty}(\R^{p,q})^{\mathrm{SO}^+(p,q)} \to \Val_k^{+,-\infty}(\R^{p',q'})^{\mathrm{SO}^+(p',q')}
\end{displaymath}
is surjective, and an isomorphism if $p',q'>0$.
Similarly, \begin{displaymath}
	i^*: \Val_k^{-,-\infty}(\R^{p,1})^{\mathrm{SO}^+(p,1)} \to \Val_k^{-,-\infty}(\R^{p',1})^{\mathrm{SO}^+(p',1)}
	\end{displaymath} 
is an isomorphism.
\end{UnnumberedCorollary}

Finally, let us comment on how this paper relates to previous work. Hadwiger considered the compact orthogonal group. In the recent work Alesker and the second named author, the Lorentz signature $\min(p,q)=1$ was considered. This case has two notable features facilitating the study of invariant generalized valuations, compared to the general signature. One is the large maximal compact subgroup $\mathrm{SO}(n-1)$, which, in combination with the description of the kernel of the cosine transform, together with its relation to the description of valuations through Crofton measures (due to Alesker-Bernstein \cite{alesker_bernstein04}), reduces the analysis to the study of invariant Crofton measures. This however only applies to even valuations, and indeed the case of odd valuations mostly remained untreated. In the present work, we use the description of valuations through invariant currents, which allows us to apply the same analysis simultaneously for all signatures, as well as for both even and odd parities.
The simplifying feature of $\mathrm O(n-1,1)$ is that the number of open orbits on each Grassmannian coincides with the dimension of the space of invariant valuations when $\min(p,q)\leq1$, but is greater when $\min(p,q)\geq 2$. This is another obstacle that only appears for general signature.

Let us also note that the classification of continuous valuations for the general signature which is carried out in this work, constitutes in fact a reduction to the Lorentz group case, which was done in \cite{alesker_faifman}. 

In a follow-up paper by the second named author \cite{faifman_crofton}, the Crofton formulas associated with $\OO(p,q)$-invariant valuations are studied.

\subsection*{Plan of the paper}

Section \ref{sec_basic_definitions} contains some basic definitions related to valuations, like continuous and generalized translation-invariant valuations, Crofton measure and the Klain embedding. 

In Section \ref{sec_klain_continuous}, we introduce the KS-continuous translation-invariant valuations as the completion of continuous translation-invariant valuations with respect to some convenient topology. We show that the space $\Val^{+,\KS}$ of even KS-continuous translation-invariant valuations is closed under the Alesker-Fourier transform, while  $\Val^{\KS}$ is closed under pull-back and push-forward by linear maps. 

The importance of this class for the main subject of the present paper comes from the easy observation that all even $\mathrm{SO}^+(p,q)$-invariant generalized translation-invariant valuations are KS-continuous.   

In Section \ref{sec_orbits} we describe the orbits on the Grassmann manifold under the action of the group $\mathrm{O}(p,q)$.

The technical heart of the paper is Section \ref{sec_dimension_generalized}, where we compute the dimension of the space of $\mathrm{SO}^+(p,q)$-invariant generalized translation-invariant valuations. 

In Section \ref{sec_klain_functions}, we prove Theorem \ref{main_thm_klain} by induction over $p+q$. For the induction base $\R^{2,2}$, we will already use some results from Section \ref{sec_r22}.  

In Section \ref{sec_classification_continuous} we classify  all continuous $\mathrm{SO}^+(p,q)$-invariant translation-invariant valuations. We use the description of the Klain functions from the previous section, and suitable pull-backs and push-forwards to reduce to the Euclidean or Lorentz case.  We also prove that odd invariant valuations are KS-continuous. This is less straightforward than in the even case.

The special case $\mathrm{SO}^+(2,2)$ will be treated in detail in Section \ref{sec_r22}. We will show that each invariant generalized valuation admits an invariant generalized Crofton measure. Using recent results from hermitian integral geometry, we will compute their Klain functions. This step will complete the proof of Theorem \ref{main_thm_klain}.  

In the appendix, we state and prove a result on generalized invariant sections on a manifold.  

\subsection*{Acknowledgments}
We wish to thank Semyon Alesker for many helpful remarks throughout the development of this project, and also Vitali Milman, Franz Schuster, Mykhailo Saienko, Gil Solanes and Thomas Wannerer for helpful comments on the first draft of this paper. A.B. thanks Gautier Berck for useful explanations. D.F. thanks Yael Karshon for numerous valuable suggestions, and owes a debt to the wonderful working atmospheres of IHP and IHES, where a large part of this work was carried out. We thank the referee for the careful checking of a first draft of this paper and many useful comments.

\section{Basic definitions}
\label{sec_basic_definitions}

\subsection{Linear algebra}

Let $\Gr_k(V)$ denote the Grassmann manifold of linear $k$-planes. We let $\Dens(V)$ denote the $\R$-span of a Lebesgue measure on the finite-dimensional vector space $V$. Given a manifold $M$, we let $|\omega|_M$ (or simply $|\omega|$) denote the linear bundle of densities over $M$, whose fiber over $x\in M$ is $\Dens(T_xM)$.

We will use some natural isomorphisms related to densities and refer to \cite[Subsection 2.1]{alesker_fourier} for more information. 

First, we have
\begin{equation} \label{eq_density_dual}
\Dens(V^*) \cong \Dens(V)^*.
\end{equation}

For $E\in\Gr_k(V)$ one has 
\begin{equation} \label{eq_exact_sequence}
\Dens(V) \cong \Dens(E) \otimes \Dens(V/E).
\end{equation}

We let $E^\perp \in \Gr_{n-k}(V^*)$ denote the annihilator of $E$. There is a natural isomorphism 
\begin{equation} \label{eq_dual_anni}
(E^\perp)^* \cong V/E.
\end{equation}

We will use the following fact from linear algebra
\begin{Lemma}\label{lem:long_quotient}
Given a hyperplane $H\subset V$, there is a natural isomorphism 
\begin{displaymath}
\largewedge^k V/\largewedge^k H \cong \largewedge^{k-1} H\otimes V/H.
\end{displaymath}
\end{Lemma}	

\proof
Note that $\dim V/H=1$, and define  $f:\largewedge^{k-1} H\otimes V/H\to\largewedge^k V/\largewedge^k H$ by setting $f(h_1\wedge\dots\wedge h_{k-1} \otimes (v+H))=h_1\wedge\dots\wedge h_{k-1}\wedge v +\largewedge^k H$. $f$ is clearly well-defined and linear. It is easily seen to be injective, and comparing the dimensions we conclude $f$ is an isomorphism.
\endproof

Let now $Q$ be a non-degenerate quadratic form of signature $(p,q)$ on $V$, with associated bilinear form $Q(\bullet,\bullet)$. Given $E \in\Gr_k(V)$, we denote by $E^Q:=\{v \in V: Q(v,e)=0, \forall e \in E\}$ the $Q$-orthogonal complement. Note that $E \cap E^Q = \{0\}$ if and only if $Q|_E$ is non-degenerate.

The quadratic form $Q$ allows us to identify $V$ and $V^*$. Under this identification, $E^Q$ corresponds to $E^\perp \in \Gr_{n-k}(V^*)$. 

A basis $v_1,\ldots,v_n$ of $V$ will be called $Q$-orthonormal if $Q(v_i,v_j)=0, i \neq j$, $Q(v_i) \in \{\pm 1\}$ for all $i$ and $Q(v_i)\geq Q(v_j)$ for all $i\leq j$. It is elementary that such a basis exists, which then identifies $(V,Q)$ with $\R^{p,q}:=(\R^n, x_1^2+\dots+x_p^2-\dots-x_n^2)$. 

The form $Q$ induces a natural choice of density $\vol_Q$ on $V$ given by 
\begin{displaymath}
\vol_Q(v_1 \wedge \ldots \wedge v_n)=|\det Q(v_i,v_j)|^\frac{1}{2}.
\end{displaymath}
which may be used to $\OO(Q)$-equivariantly identify 
\begin{equation} \label{eq_densities_trivial}
\Dens(V) \cong \R.
\end{equation}

From \eqref{eq_density_dual}, \eqref{eq_exact_sequence}, \eqref{eq_dual_anni} and \eqref{eq_densities_trivial} it follows that $\OO(Q)$-equivariantly one has
\begin{equation} \label{eq_density_q_dual}
\Dens(E) \cong \Dens(E^Q). 
\end{equation}

\begin{Lemma} \label{lemma_signatures}
Let $E \in \Gr_k(V)$. Let $(a,b)$ be the signature of $Q|_E$. Then
\begin{displaymath}
a+b \leq k, \quad \max(0,k-q) \leq a \leq p, \quad \max(0,k-p) \leq b \leq q.
\end{displaymath}
Conversely, any $(a,b)$ satisfying these inequalities is the signature of $Q|_E$ for some $E \in \Gr_k(V)$.   
\end{Lemma}

\proof
We may decompose $E=E_+ \oplus E_- \oplus E_0$ as a $Q$-orthogonal sum such that $Q|_{E_+}>0, Q|_{E_-}<0, Q|_{E_0}=0$. From the definition of the signature it follows immediately that $a=\dim E_+ \leq p, b=\dim E_-\leq q$. Since $Q|_{E_- \oplus E_0}$ is negative definite, we have $k-a=\dim E_- \oplus E_0 \leq q$ and similarly $k-b \leq p$.

Let now $(a,b)$ be given subject to the displayed inequalities and set $r:=k-a-b \geq 0$. By assumption on $(a,b)$ we have $r+a \leq p$ and $p+r+b \leq n$. Choose a $Q$-orthonormal basis $v_1,\ldots,v_p,v_{p+1},\ldots,v_n$ of $V$ and set 
\begin{displaymath}
E:=\mathrm{span}\{v_1+v_{p+1},\ldots,v_r+v_{p+r},v_{r+1},\ldots,v_{r+a},v_{p+r+1},\ldots,v_{p+r+b}\}.
\end{displaymath} 
 The signature of $E$ is clearly $(a,b)$, its dimension is $r+a+b=k$.
\endproof

If the signature of $Q|_E$ is $(0,0)$, i.e. if $E \subset E^Q$, then $E$ is called $Q$-isotropic. The collection of isotropic lines is the light cone.

Let us consider the split case $p=q$ more carefully. 
	Set
	\begin{displaymath}
	\mathcal{J}:= \{j \in \mathrm{GL}(V): j^*Q=-Q, j^2=\mathrm{Id}\}.
	\end{displaymath}
	If $Q(z, w)=\sum_{i=1}^p z_i^2-\sum_{i=1}^p w_i^2$ is the standard form of signature $(p,p)$ on $\mathbb R^{p,p}$, then the map $j(z,w)=(w,z)$ belongs to $ \mathcal{J}$. 
	
\begin{Definition} \label{def_split}
Let $X$ be a complex $\mathrm{GL}(n)$-module and $j \in \mathcal{J}$. Since $\mathrm O(Q)=\mathrm O(-Q)$, $j$ is an automorphism of $X^{\mathrm O(Q)}$ with eigenvalues $\pm 1$. We may thus decompose into $j$-even and $j$-odd elements,
\begin{equation} \label{eq_split}
X^{\mathrm O(Q)}=X^{\mathrm O(Q),j} \oplus X^{\mathrm O(Q),-j},
\end{equation}
where 
$X^{\mathrm O(Q),\pm j}:=\{v \in X^{\mathrm O(Q)}:jv=\pm v\}$.
\end{Definition}

\begin{Lemma}
The decomposition \eqref{eq_split} is independent of the choice of $j \in \mathcal{J}$. 
\end{Lemma}
	
\proof
Let $j \in \mathcal{J}$ be fixed. We claim that any other $j' \in \mathcal{J}$ is of the form $j'=gj$ with $g \in \mathrm{O}(Q)$ such that $gj=jg^{-1}$.
	
	To prove the claim, first note that $(j')^*Q=-Q=j^*Q$, hence $g:=j' j^{-1} \in \mathrm{O}(Q)$. Set $h:=jgj \in \mathrm{GL}(V)$. Then $\mathrm{Id}=(j')^2=gjgj=gh$, hence $h=g^{-1} \in \mathrm{O}(Q)$.
	
	From the claim and using $\mathrm O(Q)$-invariance, it follows that 
	\begin{align*}
		X^{\mathrm O(Q),\pm j'} & =\{v \in X^{\mathrm O(Q)}:j'v=\pm v\} =\{v \in X^{\mathrm O(Q)}:gjv=\pm v\}\\
		& =\{v \in X^{\mathrm O(Q)}:jv=\pm v\}=X^{\mathrm O(Q),\pm j}
	\end{align*}
	\endproof 

The isometry group $\OO(Q)=\mathrm{O}(p,q)$ has four connected components when $p,q\geq 1$, corresponding to the orientations of maximal positive- and negative definite subspaces. The connected component of the identity is denoted $\SO^+(p,q)$. A maximal compact subgroup of $\mathrm{O}(Q)$ is isomorphic to $\mathrm{O}(p) \times \mathrm{O}(q)$ - for a proof, see e.g. \cite{wolf_book2011}.

Witt's theorem \cite[Thm. 5.3]{scharlau85} asserts that if $E,F\in\Gr_k(V)$ and $T:E\to F$ is a linear map such that $T^*Q|_F=Q|_E$ then one can find $\tilde T\in\OO(Q)$ with $\tilde T|_E=T$.

\subsection{The interplay between Euclidean and $(p,q)$- geometries}
\label{subsec_interplay}

Let $P$ be an auxiliary Euclidean scalar product on $V$. One has a unique simultaneous diagonalization of the forms: $V=\bigoplus_{j=1}^{N(P,Q)} V_j$ s.t. $Q|_{V_j}=\lambda_j P|_{V_j}$ with all $\lambda_j$ distinct and all $V_j$ pairwise $P$- and $Q$- orthogonal. Note that $N(P,Q) \geq 2$ if $Q$ is indefinite, which we assume from now on. 
Denote by $S:V\to V$ the map given by $Q(u,v)=P(Su,v)$.

\begin{Lemma} \label{lem:maximal_compact}
	The following are equivalent:
	\begin{enumerate}
		\item $N(P,Q)=2$.
		\item $\mathrm{O}(P)\cap \mathrm{O}(Q) \subset \mathrm{O}(Q)$ is a maximal compact subgroup.
	\end{enumerate}
\end{Lemma}

\proof
\mbox{}
\begin{enumerate}
	\item[(i) $\implies$ (ii)] Assume the decomposition is $V=V_p\oplus V_q$ with $\lambda_p>0$ and $\lambda_q<0$. Then $\mathrm{O}(V_p) \times \mathrm{O}(V_q) \subset \mathrm{O}(P) \cap \mathrm{O}(Q)$. Since $\mathrm{O}(V_p) \times \mathrm{O}(V_q) \cong \mathrm{O}(p) \times \mathrm{O}(q)$ is a maximal compact subgroup of $\mathrm{O}(Q)$ and $\mathrm{O}(P) \cap \mathrm{O}(Q)$ is a compact subgroup of $\mathrm{O}(Q)$, we must have equality. 
	\item[(ii) $\implies$ (i)] The group $\mathrm{O}(P)\cap \mathrm{O}(Q)$ leaves each $V_j$ invariant. It is then easy to see that $\mathrm{O}(P)\cap \mathrm{O}(Q)=\mathrm{O}(V_1)\times \dots \times \OO(V_N)$.  As $\OO(P) \cap \OO(Q) \cong \OO(p)\times\OO(q)$ by assumption, it follows that $N=2$.
\end{enumerate}
\endproof

We will also need the following description:
\begin{Lemma} \label{lem:compatible_definition}
	The following are equivalent:
	\begin{enumerate}
		\item $N(P,Q)=2$  and  $\lambda_p=1$, $\lambda_q=-1$.
		\item $S$ satisfies $S^2=\mathrm{Id}$.
	\end{enumerate}
\end{Lemma}

\proof
\mbox{}
\begin{enumerate}
	\item[(i) $\implies$ (ii)] It is immediate that $S|_{V_p}=\mathrm {Id}$ and $S|_{V_q}=-\mathrm{Id}$.
	\item[(ii) $\implies$ (i)] Since $P(Su,v)=Q(u,v)=Q(v,u)=P(Sv,u)=P(u,Sv)$, $S$ is self-adjoint with respect to $P$. Choose an eigenbasis $e_j$ of $S$ which is $P$-orthonormal. Since $S^2=\mathrm{Id}$, the eigenvalues of $S$ are $\pm 1$. Therefore, $Q(e_i,e_j)=P(e_i, Se_j)=\pm \delta_{ij}$. Thus $N(P,Q)=2$ with $V_p$ the positive eigenspace of $S$ and $V_q$ the negative eigenspace. 
\end{enumerate}

\endproof

\begin{Definition}  \label{def_compatible_forms}
	A Euclidean form $P$ is called compatible with $Q$ if $N(p,q)=2$ and $\lambda_p=1$, $\lambda_q=-1$. We then have the decomposition $V=V_p\oplus V_q$, where $\dim V_p=p$, $\dim V_q=q$, $V_p=\{x: P(x)=Q(x)\}$, $V_q=\{x: P(x)=-Q(x)\}$.
\end{Definition}

Clearly such a compatible Euclidean form $P$ exists for any $Q$. Let us fix a $Q$-compatible Euclidean form $P$, with the associated involution $S:V\to V$ s.t. $Q(u,v)=P(u, Sv)$.

\begin{Lemma}\label{Lem:TwoFormsOrthogonalComplements}
	\begin{enumerate}
		\item $Q(S u,S v)=Q(u,v)$, $P(S u,S v)=P(u,v)$.
		\item $P(v,S u)=P(S v, u)=Q(v,u)$, $Q(v,S u)=Q(S v, u)=P(v,u)$.
		\item For any subspace $E \subset V$, $(SE)^P=S(E^P)=E^Q$ and $(SE)^Q=S(E^Q)=E^P$,  $SE=(E^Q)^P=(E^P)^Q$.
		\item $E= SE \iff E^P=E^Q$.
		\item If $E$ is $Q$-isotropic then $E$ is $P$-orthogonal to $SE$.
		
	\end{enumerate}
\end{Lemma}

\proof
The first two items follow directly from the definition of $S$. The third item is then an easy consequence.

Let us show (iv). If $E=SE$, then $E^P=E^Q$ follows immediately from (iii). If $E^P=E^Q$, then $SE=(E^Q)^P=(E^P)^P=E$. For (v), we use (iii) to obtain $E \subset E^Q =(SE)^P$.
\endproof

\begin{Lemma} \label{lemma_isotropic_complement}
Let $Q$ be a non-degenerate quadratic form on $V$, where $\dim V=2r$. Given $E \in \Gr_r(V)$ isotropic, there exists $F \in \Gr_r(V)$ isotropic such that $V=E \oplus F$.
\end{Lemma}

\proof
Since $E$ is isotropic, we have $E \subset E^Q$. Since both spaces are of dimension $r$, we actually have $E=E^Q$. Fix a compatible Euclidean form $P$ and set $F:=E^P$. By Lemma \ref{Lem:TwoFormsOrthogonalComplements}iii), $F=E^P=(E^Q)^P=(E^P)^Q=F^Q$, hence $F$ is isotropic. 
\endproof

\subsection{Smooth and generalized valuations}

The space of continuous and translation-invariant valuations on a vector space $V$ is denoted by $\Val(V)$. Examples of elements in $\Val$ are mixed volumes, i.e. valuations of the form $K \mapsto V(K[i],L_{i+1},\ldots,L_n)$ for fixed $L_{i+1},\ldots,L_n \in \mathcal{K}(V)$. By Alesker's solution of McMullen's conjecture \cite{alesker_mcullenconj01}, mixed volumes span a dense subset in $\Val(V)$, the latter space equipped with a certain natural locally convex topology.

A valuation $\mu \in \Val(V)$ is of degree $k$ if $\mu(tK)=t^k\mu(K)$ for all $t \geq 0$ and all $K$. It is even if $\mu(-K)=\mu(K)$ and odd if $\mu(-K)=-\mu(K)$. The corresponding subspaces of $\Val(V)$ are denoted by $\Val_k^+(V), \Val_k^-(V)$. 

McMullen \cite{mcmullen77} proved the decomposition
\begin{displaymath} 
  \Val(V)=\bigoplus_{\substack{k=0,\ldots,n\\ \epsilon = \pm}} \Val_k^\epsilon(V). 
\end{displaymath}

Let us denote by $\mathbb{P}_+(V^*)$ the set of all oriented lines in $V^*$, or equivalently the set of all co-oriented hyperplanes in $V$. The cosphere bundle of $V$ is defined as $V \times \mathbb P_+(V^*)$. 

An important subspace of $\Val$ is the space $\Val^\infty$ of smooth valuations. They are given by the translation-invariant valuations of the form
\begin{equation} \label{eq_smooth_vals}
 K \mapsto \int_{N(K)} \omega+ c \vol(K),
\end{equation}
where $c \in \C$, $N(K)$ is the normal cycle of $K$ \cite{zaehle86} and $\omega$ is a complex-valued translation-invariant $(n-1)$-form on the cosphere bundle. 

The cosphere bundle $V \times \mathbb P_+(V^*)$ is a contact manifold. A form vanishing on the contact distribution is called vertical. Fixing a choice of a contact form $\alpha$, a form $\omega$ is vertical if and only if $\alpha \wedge \omega=0$ or equivalently $\omega=\alpha \wedge \tau$ for some form $\tau$. We will write $\Omega_V$ for the vertical forms. The normal cycle of a compact convex body is a Legendrian cycle, i.e. it vanishes on vertical forms and also on exact forms. Consequently, the form $\omega$ in the definition of a smooth valuation is not unique, compare \cite{bernig_broecker07} for more information.

Alesker \cite{alesker04_product} has introduced a graded product structure on the space $\Val^\infty$. It is characterized by the property that if $\phi_i(K)=\vol_n(K+A_i)$ with $A_i$ a smooth convex body with positive curvature, then $\phi_1 \cdot \phi_2(K)=\vol_{2n}(\Delta K+A_1 \times A_2)$. Here $\Delta:V \to V \times V$ is the diagonal embedding, $\vol_n$ is any choice of Lebesgue measure on $V$ and $\vol_{2n}$ is the corresponding product measure on $V \times V$.  With respect to the natural Fr\'echet space topology on $\Val^\infty$, the product is continuous. The product is closely related to kinematic formulas and was a key ingredient in the determination of such formulas on hermitian vector spaces, see \cite{bernig_fu06, bernig_fu_hig}. 

The space $\Val_n(V)$ can be identified with the one-dimensional space $\Dens(V)$ of densities on $V$. Then the bilinear map given by the product
\begin{displaymath}
 \Val^\infty_k(V) \times \Val^\infty_{n-k}(V) \to \Val^\infty_n(V)\cong \Dens(V)
\end{displaymath}
is a perfect pairing, i.e. the induced map 
\begin{displaymath}
 \Val_k^\infty(V) \to \Val^\infty_{n-k}(V)^* \otimes \Dens(V)
\end{displaymath}
is injective. 
Moreover , by \cite[Proposition 8.1.2]{alesker_fourier} there is an extension to a continuous product $\Val_k^\infty(V) \otimes \Val_{n-k}(V) \to \Dens(V)$.
Elements of the space 
\begin{displaymath}
 \Val^{-\infty}_k(V) := \Val^\infty_{n-k}(V)^* \otimes \Dens(V)
\end{displaymath}
are called generalized translation-invariant valuations of degree $k$. Equipping $  \Val^{-\infty}_k(V)$ with the weak dual topology, $\Val^{\infty}_k(V)\subset\Val^{-\infty}_k(V)$ is a dense subspace. Generalized translation-invariant valuations were recently related to McMullen's polytope algebra, compare \cite{bernig_faifman}. 

Let us explain two constructions related to even valuations which will be crucial in the following. 

By a result of Klain, even translation-invariant (continuous) valuations can be described in terms of their Klain functions. Let $\phi \in \Val_k^+(V)$. Given a subspace  $E \in \Gr_k(V)$, by a theorem of Hadwiger, the restriction of $\phi$ to $E$ is a multiple of the volume on $E$. We thus get a continuous global section (often called the Klain function of $\phi$) of the Klain bundle $K^{n,k}$ whose fiber over $E \in \Gr_k(V)$ is the one-dimensional space $\Dens(E)$ of densities on $E$. Klain has shown \cite{klain00} that the corresponding map (called Klain map) $\Kl:\Val_k^+(V) \to \Gamma(K^{n,k})$ is injective. 

For even generalized translation-invariant valuations of degree $k$, the Klain map can be extended to an embedding 
\begin{displaymath}
    \Kl_k: \Val_k^{+,-\infty}(V) \to \Gamma^{-\infty}(K^{n,k}),
\end{displaymath}
where the latter space consists of all generalized sections of Klain's bundle. This map is still injective, see \cite[Prop. 4.4]{alesker_faifman}. A generalized translation-invariant valuation of degree $k$ may thus be uniquely described by its (generalized) Klain function. We will do this explicitly for even $\mathrm{SO}^+(p,q)$-invariant valuations in Section \ref{sec_klain_functions}. 

The second description of even valuations is through Crofton measures. Given a Euclidean structure and a smooth measure $m$ on $\Gr_{n-k}(V)$, the following is a  smooth, even, translation-invariant valuation of degree $k$: 
\begin{displaymath}
 K \mapsto \int_{\Gr_{n-k}(V)} \vol_k(\pi_{E^\perp}K)dm(E).
\end{displaymath}
The measure $m$ is called Crofton measure of this valuation.

In $\mathrm{GL}(V)$-equivariant terms, this construction can be described as follows. 

The Crofton bundle $C^{n,k}$ over $\Gr_{n-k}(V)$ is defined as the bundle whose fiber over $E \in \Gr_{n-k}(V)$ is given by 
\begin{displaymath}
 C^{n,k}|_E:=\Dens(V/E) \otimes \Dens(T_E \Gr_{n-k}(V)),
\end{displaymath}
 where $T_E \Gr_{n-k}V$ is the tangent space of the Grassmannian at the plane $E$. 

Given a smooth section $s \in \Gamma^\infty(C^{n,k})$, the smooth even valuation $\mathrm{Cr}_{n-k}(s)$ is defined by 
\begin{displaymath}
 \mathrm{Cr}_{n-k}(s)(K):=\int_{E \in \Gr_{n-k}(V)} s(\mathrm{Pr}_{V/E} K).
\end{displaymath}
The map 
\begin{displaymath}
 \mathrm{Cr}_{n-k}: \Gamma^\infty(C^{n,k}) \to \Val_k^{+,\infty}(V)
\end{displaymath}
is surjective \cite{alesker_bernstein04}.

The composition $T_{n-k,k}=\Kl_k \circ \mathrm{Cr}_{n-k}:\Gamma^\infty(C^{n,k}) \to \Gamma^\infty(K^{n,k})$ is the cosine transform, written in $\mathrm{GL}(V)$-equivariant form \cite{alesker_bernstein04}. In particular, the kernel of $\Cr_{n-k}$ equals the kernel of $T_{n-k,k}$.

The Crofton map $\mathrm{Cr}_{n-k}$ can be extended to generalized Crofton measures and generalized translation-invariant valuations, i.e. we get a map 
\begin{displaymath}
 \mathrm{Cr}_{n-k}: \Gamma^{-\infty}(C^{n,k}) \to \Val_k^{+,-\infty}(V).
\end{displaymath}
Again this map is surjective \cite[Prop. 4.5]{alesker_faifman}. When no confusion can arise, we omit the dimension $k$ in $\Kl_k$, $\Cr_k$.

If a Euclidean structure is fixed on $V$, both line bundles $K^{n,k}$ and $C^{n,k}$ acquire a natural trivialization, and the spaces $\Gr_k(V)$ and $\Gr_{n-k}(V)$ are naturally identified. In this case we will write shortly $T_k:C^\infty(\Gr_k(V)) \to C^\infty(\Gr_k(V))$ for the cosine transform, which is then given by

\begin{displaymath}
 T_{k}(f)(F)=\int_{\Gr_k(V)}f(E)\langle E,F\rangle dE,
\end{displaymath}
where $\langle E,F\rangle$ denotes the cosine of the angle between $E$ and $F$. 

This map is self-adjoint and thus extends to a map $T_k:C^{-\infty}(\Gr_k(V)) \to C^{-\infty}(\Gr_k(V))$. As an example, the cosine transform of the Dirac generalized function $\delta_E, E \in \Gr(V,k)$ is given by the continuous function 
\begin{equation} \label{eq_cosine_delta} 
 T_k \delta_{E}=\langle E,\bullet\rangle.
\end{equation}

Any even generalized translation-invariant valuation can be evaluated at a smooth convex body  $K$ with positive curvature. Indeed, given a generalized Crofton measure, we just apply it to the smooth function $E \mapsto \vol(\pi_EK)$, where $\pi_EK$ is the orthogonal projection onto a $k$-dimensional subspace $E$ (here we use a Euclidean trivialization to keep notations short).

Alesker has introduced in \cite{alesker_fourier} a linear map
\begin{displaymath}
 \mathbb{F}:\Val_k^\infty(V) \to \Val_{n-k}^\infty(V^*) \otimes \Dens(V). 
\end{displaymath}
 usually called the Alesker-Fourier transform.
 
For even valuations, it is characterized in terms of their Klain function by the equation
\begin{displaymath}
 \Kl_{\mathbb F \phi}(E^\perp)=\Kl_\phi(E), \quad E \in \Gr_k(V). 
\end{displaymath}
There is a canonical isomorphism $(E^\perp)^* \cong V/E$, hence the term on the left hand side of the equation is in $\Dens(E^\perp) \otimes \Dens(V) \cong \Dens^*(V/E) \otimes \Dens(V) \cong \Dens(E)$, as is the term on the right hand side. 

Let $Q$ be a non-degenerate bilinear form on $V$ used to identify $V \cong V^*$ and $\Dens(V) \cong \R$. Then $ \mathbb{F}:\Val_k^\infty(V) \to \Val_{n-k}^\infty(V)$. Under the identification \eqref{eq_density_q_dual}, we have the equation
\begin{equation} \label{eq_fourier_even_q}
 \Kl_{\mathbb F \phi}(E^Q)=\Kl_\phi(E), \quad E \in \Gr_k(V).
\end{equation}

\section{Klain-Schneider-continuous valuations}
\label{sec_klain_continuous}

\subsection{Definition of KS-continuous valuations}

Alesker \cite{alesker_fourier} defined pull-back and push-forward under linear maps of continuous translation-invariant valuations. We will extend these operations to the larger class of KS-continuous valuations. 

Let us first recall the different topologies on valuation spaces used in the following. 

The space $\Val(V)$ of translation-invariant valuations has a Banach space structure as follows. Fixing some compact  convex set $B \subset V$ with non-empty interior, the norm of $\phi$ is defined by
\begin{displaymath}
\|\phi\|:=\sup \{|\phi(K)|: K \subset B\}.
\end{displaymath}
Changing $B$ gives a different, but equivalent norm. 

Next, the subspace $\Val^\infty$ of smooth valuations has a natural Fr\'echet space topology. With respect to \eqref{eq_smooth_vals}, this topology may be defined as the quotient topology arising from the natural Fr\'echet space topology on the space of translation-invariant smooth $(n-1)$-forms on the cosphere bundle. 

The space $\Val^{-\infty} = \Val^\infty_{n-k}(V)^* \otimes \Dens(V)$ of generalized translation-invariant valuations is equipped with the weak dual topology. Being the dual of a Fr\'echet space, it is sequentially complete (and even quasi-complete), see \cite{schaefer_wolff}.

Let us now recall the Schneider embedding of continuous translation-invariant valuations. Let $0 \leq k < n$ and $\phi \in \Val_k(V)$. The restriction of $\phi$ to a $(k+1)$-dimensional subspace $E \subset V$ is a continuous, $k$-homogeneous, translation-invariant valuation in $E$, i.e. an element of $\Val_k(E)$. 
Consider the infinite dimensional Banach bundle over $\Gr_{k+1}(V)$, with fiber $\Val_{k}(E)$ over $E$, equipped with a norm defined by $E\cap B$, where $B\subset V$ is some fixed convex set with non-empty interior. We will call it the Schneider bundle.

\begin{Definition}
The Schneider map is given by 
\begin{align*}
	\Sc:\Val_k(V) & \to \Gamma(\Gr_{k+1}(V),\Val_k(E))\\
	\phi & \mapsto [E \mapsto \phi|_E]	
\end{align*}
\end{Definition} 
This map is an injection, \cite{schneider96}. 

Let us define two more norms on translation-invariant valuations, which are related to the Klain and Schneider embeddings.

\begin{Definition}
\begin{enumerate}
\item Define the norm $\|\cdot\|_K$ on $\Val_k^ +(V)$ by
\begin{displaymath}
\|\phi\|_{K}:=\|\Kl_\phi\|_{\infty}
\end{displaymath}
Here the supremum norm arises from the identification $\Gamma(\Gr_k(V), \Dens(E))=C(\Gr_k(V))$ induced by an arbitrary choice of a Euclidean structure on $V$.
\item Fix again an auxiliary Euclidean structure on $V$, and equip all spaces of continuous valuations with the Banach norm associated with the Euclidean unit ball. The space of continuous sections $\Gamma(\Gr_{k+1}(V), \Val_k(E))$ is a Banach space with the supremum norm. Define the norm $\|\bullet\|_{\Sc}$ on $\Val_k(V)$ by 
\begin{displaymath}
\|\phi\|_{\Sc}:=\|\Sc_\phi\|_{\infty}
\end{displaymath}
\end{enumerate}
\end{Definition}
The norms corresponding to different Euclidean structures are equivalent.
\begin{Lemma}
The restriction of $\|\bullet\|_{\Sc}$ to $\Val_k^+(V)$ is equivalent to the norm $\|\bullet\|_K$.
\end{Lemma}

\proof
We have to show that there are constants $c, C>0$ s.t. for $\phi\in \Val_k^+(V)$, \begin{displaymath}
c\|\Kl_\phi\|_\infty \leq \|\Sc_\phi\|_\infty \leq C\|\Kl_\phi\|_\infty.
\end{displaymath}

The left inequality is obvious. For the right inequality, let $B(E)$ be the unit ball inside $E$, $S(E)$ the unit sphere in $E$ and $\sigma_K$ the surface area measure of $K$. Then
\begin{align*} 
\|\Sc_\phi \|_\infty & = \sup_{E \in \Gr_{k+1}(V)} \sup_{K\subset B(E)} |\phi(K)|\\
& =\sup_{ E \in \Gr_{k+1}(V)} \sup_{K \subset B(E)} \left|\int _{\theta \in  S(E)} \Kl_\phi(\theta^\perp) d\sigma_K(\theta)\right| \\
& \leq \|\Kl_\phi\|_\infty \sup_{E\in \Gr_{k+1}(V)} { \sigma_{B(E)}(S(E))}\\
& = \omega_{k+1} \|\Kl_\phi\|_\infty
\end{align*}
where the inequality follows by monotonicity of the surface area, and $\omega_{k+1}$ is the surface area of the unit sphere in a $(k+1)$-dimensional Euclidean vector space. 
\endproof

\begin{Definition}
The space $\Val^{\KS}_k(V)$ of KS-continuous, $k$-homogeneous valuations is the completion of $\Val_k(V)$ in the norm $\|\phi\|_{\Sc}$.
\end{Definition}
We now define a generalization of the Crofton map from Section \ref{sec_basic_definitions} from even to all valuations. Consider the map
\begin{equation*}
\Sc^*: \Gamma^\infty (\Gr_{k+1}(V), \Val_1^{\infty}(E) \otimes \Dens^*(E) \otimes|\omega|) \to  \Val_{k}^{\infty}(V)^*
\end{equation*}
given by 
\begin{displaymath}
\langle \Sc^*(\mu), \phi\rangle =  \int _{\Gr_{k+1(V)}} \langle \mu(E), \Sc_\phi(E)\rangle, \quad \phi \in \Val_k^\infty(V).  
\end{displaymath}

The bracket has to be understood as follows. For a fixed $E$, we have $\mu(E) \in \Val_1^\infty(E) \otimes \Dens^*(E) \otimes |\omega|$ and $\Sc_\phi(E) \in \Val_k^\infty(E)$. The Alesker product gives us a map $\Val_1^\infty(E) \otimes \Val_k^\infty(E) \to \Val_{k+1}^\infty(E) \cong \Dens(E)$, so that altogether we obtain that $\langle \mu(E), \Sc_\phi(E)\rangle \in |\omega|$. We thus get a section of $|\omega|$ which can be integrated over $\Gr_{k+1}(V)$.

The notation $\Sc^*$ is justified by the fact that the dual map $ (\Sc^*)^*$ will be shown to fit into a commutative diagram 
\begin{displaymath}
\xymatrix{ \Val_k(V) \ar@{^{(}->}[d] \ar[r]^-{\Sc} & \Gamma(\Gr_{k+1}(V), \Val_k(E)) \ar@{^{(}->}[d] \\
	\Val_k^{-\infty}(V) \ar[r]^-{(\Sc^*)^*} & \Gamma^\infty (\Gr_{k+1}(V), \Val_1^{\infty}(E) \otimes \Dens^*(E) \otimes|\omega|)^*.
}
\end{displaymath}

Recall that by the Alesker-Poincare duality, there is a natural dense inclusion $\Val_{n-k}^{\infty}(V) \otimes \Dens^*(V) \subset  \Val_k^{\infty}(V)^*$.
\begin{Lemma}
The image of $\Sc^*$ is contained in  $\Val_{n-k}^{\infty}(V) \otimes \Dens^*(V)$, and
\begin{equation}\label{eqn:Sc*_def}
	\Sc^*: \Gamma^\infty (\Gr_{k+1}(V), \Val_1^{\infty}(E) \otimes \Dens^*(E) \otimes|\omega|) \to  \Val_{n-k}^{\infty}(V) \otimes \Dens^*(V)
\end{equation}
is continuous.
\end{Lemma}

\proof
Recall the natural map $i_{k}:\mathcal K(V)\to \Val_{k}(V)\otimes \Dens^*(V)$ given by $i_k(K)\otimes \sigma=V_\sigma(\bullet[k], -K[n-k])$, where we write $V_\sigma$ for the mixed volume induced by $\sigma\in\Dens(V)$.
It then holds \cite[Lemma 3.1]{bernig_faifman} that for $\phi\in \Val_{n-k}^\infty(V)$,  $\phi(K)=\langle \phi, i_k(K)\rangle$, the Alesker-Poincar\'e pairing.

Given $\mu$ in the domain of definition of $\Sc^*$ and $\sigma\in\Dens(V)$, define $\phi_\mu\otimes \sigma \in\Val_{n-k}(V)$ by setting for $K\in\mathcal K(V)$
\[ \phi_\mu\otimes\sigma (K):=  \int _{\Gr_{k+1(V)}} \langle \mu(E), (i_k(K)\otimes\sigma)|_E\rangle\]
Note that the pairing under the integral is between a smooth and a continuous valuation continuously dependent on $E$ in the respective topologies, and is therefore continuous in $E$. 

Considering bodies $K\in\mathcal K(V)$ with smooth support function, we see that $\phi_\mu\otimes \sigma (K)=\langle \Sc^*(\mu)\otimes\sigma, i_k(K)\rangle$, so that $\Sc^*(\mu)=\phi_\mu\in \Val_{n-k}(V)\otimes \Dens^*(V)$, and moreover
\[
\Sc^*: \Gamma^\infty (\Gr_{k+1}(V), \Val_1^{\infty}(E) \otimes \Dens^*(E) \otimes|\omega|) \to   \Val_{n-k}(V)\otimes \Dens^*(V)
\]
is continuous. Since the domain of definition is a smooth $\GL(v)$-module, the $\GL(V)$-equivariance of $\Sc^*$ then implies that the image of $\Sc^*$ consists of smooth vectors, and $\Sc^*$ is a continuous map between the corresponding Fr\'echet spaces.

\endproof

\begin{Proposition}
The image of $\Sc^*$ equals  $\Val_{n-k}^{\infty}(V) \otimes \Dens^*(V)$.
\end{Proposition}
\proof
For a fixed $E \in \Gr_{k+1}(V)$, denote by $\Gr_k^-(E)$ the manifold of cooriented hyperplanes $F\subset E$. For a body $K\in \mathcal{K}(E)$, its equivariantly written support function is $\tilde h_K(F)=\max_{x\in K} \Pr_{E/F}(x) \in \Gamma(\Gr_k^-(E) , \Dens^*(E/F))$ (the maximum is taken using the orientation of $E/F$). It is Minkowski additive, namely $\tilde {h}_{K+L}=\tilde h_K+\tilde h_L$. It relates to the standard version of the support function $h_K(\xi)=\sup_{x \in K} \xi(x), \xi \in V^*$ by $\tilde h_K(\xi^\perp)= h_K(\xi) [\xi^*]$, where $[\xi^*]\in \Dens^*(E/\xi^\perp)=\Dens(\Span(\xi))$ is determined by $\xi$.

Let $\widetilde{\Val}_1(E)$ denote the space of continuous $1$-homogeneous valuations with the property that $\phi(K+x)-\phi(K)$ is a linear functional of $x\in E$, independent of $K\in\mathcal K(E)$. It is naturally a $\GL(E)$-module, and one easily verifies that it fits into an exact sequence 
\begin{displaymath}
0 \to \Val_1(E) \to \widetilde{\Val}_1(E) \to E^*\to 0 
\end{displaymath}
- exactness at $E^*$ follows by lifting $\xi\in E^*$ to $[K \mapsto h_K(\xi)] \in \widetilde{\Val}_1(E)$.

Thus $\widetilde{\Val}_1(E)$ can be considered as the total space of a fibration over $E^*$ with fiber $\Val_1(E)$, with the natural Banach space topology. Considering the smooth elements of this representation, we get the Fr\'echet space $\widetilde\Val^\infty_1(E)$ and the short exact sequence of smooth $\GL(E)$-modules \[0\to \Val^\infty_1(E)\to \widetilde{\Val}^\infty_1(E) \to E^*\to 0 \]
It follows by Alesker's irreducibility theorem that $\widetilde{\Val}^\infty_1(E)$ is admissible and of finite length.

The map
\begin{displaymath}
T_E: \Gamma^\infty(\Gr_k^-(E), \Dens(E/F)\otimes |\omega|_{\Gr^-_k(E)}) \to \widetilde{\Val}_1^{\infty}(E) 
\end{displaymath}	 
given by 
\begin{displaymath}
 T_E(s)(K)=\int_{F\in \Gr^-_k(E)}\langle s(F), \tilde h_K(F)\rangle
\end{displaymath}	
is then a $\GL(E)$-equivariant map of Fr\'echet spaces. 

We claim $T_E$ is an isomorphism. Let us fix a Euclidean structure on $E$. 
Injectivity follows from the fact that any $f\in C^2( S(E))$ can be represented as $f=h_{K}-h_L$ for some $K,L\in \mathcal K(E)$. For surjectivity, we note that $\Im(T_E)$ is a closed invariant subspace (by the Cassleman-Wallach theorem) which contains $\Val_1^\infty(E)$ (e.g. by Alesker's irreducibility theorem), but is strictly larger. Since $E^*$ is irreducible, $\Im(T_E)=\widetilde{\Val}^\infty_1(E)$.
	
Consider the flag manifold of cooriented pairs 
\begin{displaymath}
\Fl_{k+1,k}=\{E\supset F: E\in\Gr_{k+1}(V), F\in\Gr_k^-(E)\}. 
\end{displaymath}

Noting the natural isomorphism 
\begin{displaymath}
\Dens(E/F)\otimes |\omega|_{\Gr^-_k(E)}\otimes \Dens^*(E) \otimes|\omega|_{\Gr_{k+1}} \cong \Dens^*(F)\otimes |\omega|_{\Fl_{k+1,k}},
\end{displaymath}
we obtain a map
\begin{multline*} T:\Gamma^\infty(\Fl_{k+1, k}(V), \Dens^*(F)\otimes |\omega|_{\Fl_{k+1,k}})\\
	\to \Gamma^\infty(\Gr_{k+1}(V),  \widetilde{\Val}_1^{\infty}(E) \otimes \Dens^*(E) \otimes|\omega|_{\Gr_{k+1}}) 
\end{multline*}
given by $T(s)(E)=T_E(s(E,\bullet))$, which is clearly $\GL(V)$-equivariant, continuous and injective. Let us check it is surjective. Take $\psi \in  \Gamma^\infty(\Gr_{k+1}(V),  \widetilde\Val_1^{\infty}(E) \otimes \Dens^*(E) \otimes|\omega|_{\Gr_{k+1}})$, and define $s(E, \bullet )=T_E^{-1}(\psi(E))$, which is clearly a smooth section over $\Fl_{k+1,k}$ since $T_{g^{-1}E}=g^*T_Eg_*$ and so $T_E^{-1}$ depends smoothly on $E$. Thus $T$ is a $\GL(V)$-equivariant isomorphism of Fr\'echet spaces. The domain of $T$ can be considered as the space of smooth vectors of a Banach representation of $\GL(V)$. It follows from \cite[Lemma 11.5.1]{wallach_book2} that the target space of $T$ is a smooth $\GL(V)$-module of moderate growth. The domain of $\Sc^*$ is a closed subspace of the latter, and thus is also a smooth $\GL(V)$-module of moderate growth. 

Now it is shown in \cite{alesker_mcullenconj01} that the target space of $\Sc^*$ is admissible and of finite length. By the Casselman-Wallach theorem \cite{casselman89}, the image of $\Sc^*$ is closed, and obviously it maps sections of even resp. odd valuations, to valuations of the corresponding parity. By Alesker's irreducibility theorem \cite{alesker_mcullenconj01}, $\Sc^*$ is surjective.
\endproof

	 It follows that 
\begin{displaymath}
(\Sc^*)^*:\Val_k^{-\infty}(V)\to \Gamma^\infty (\Gr_{k+1}(V), \Val_1^{\infty}(E) \otimes \Dens^*(E) \otimes|\omega|)^*,
\end{displaymath} 
which will be denoted simply by $\Sc$, is injective. 
	
\begin{Proposition}
The embedding $i:\Val_k(V)\subset \Val_k^{-\infty}(V)$ extends by continuity to a natural embedding $i:\Val_k^{\KS}(V)\subset \Val_k^{-\infty}(V)$.
\end{Proposition}

\proof
Assume $\phi_j \in\Val_k(V)$ is a Cauchy sequence in the $\|\bullet\|_{\Sc}$ norm. Let us show that this sequence is weakly convergent in $\Val_k^{-\infty}(V)$. 

The product $\Val^\infty(V) \otimes \Val^\infty(V) \to \Val^\infty(V)$ extends to a continuous product $\Val^\infty(V) \otimes \Val(V) \to \Val(V)$ \cite[Proposition 8.1.2]{alesker_fourier}. Therefore, for $E\in\Gr_{k+1}(V), \psi \in \Val_1^{\infty}(E) \otimes\Dens^*(E)$  and $\phi \in \Val_k^\infty(V)$, it holds that $|\langle \psi, \phi|_E\rangle |\leq c_\psi\|\phi|_E\|$, where $c_{\psi}$ is a continuous semi-norm on $\psi$. It follows that  
\begin{align*}
|\langle \Sc^*(\mu), \phi_j \rangle | & = \int _{\Gr_{k+1(V)}} \langle \mu(E), \phi_j|_E\rangle  \\
& \leq  \|\phi_j\|_{\Sc}  \int_{\Gr_{k+1}(V)} c_{\mu(E)}\\
& =:C_\mu\|\phi_j\|_{\Sc}.
\end{align*}

Since $\Sc^*$ is surjective, it follows that $\phi_j$ is weakly Cauchy, and by the sequential completeness of $\Val^{-\infty}(V)$, it weakly converges to some generalized valuation $\phi\in \Val_k^{-\infty}(V)$. Denoting $\tilde{\phi}$ the limit of $\phi_j$ in $\Val_k^{\KS}(V)$, we set $i(\tilde \phi):=\phi$. This map obviously does not depend on the choice of $\phi_j$, and extends $i$ to the KS-continuous valuations.
 
It holds that 
\begin{displaymath}
|\langle i(\tilde \phi), \Sc^*(\mu)\rangle | \leq {C_\mu} \|\tilde{\phi}\|_{\Sc}
\end{displaymath}
so that $i$ is continuous. 
To see that $i$ is injective, first note that $i$ is $\GL(V)$-equivariant. Suppose $i$ has a non-trivial kernel in $\Val_k^{\KS}$. Then the dense subspace of smooth vectors in the kernel is a non-trivial subspace of $\mathrm{GL}(V)$-smooth vectors that would then lie in $\Val_k(V)$. Since $i$ has trivial kernel in $\Val_k(V)$, this is a contradiction.

\endproof

For an even valuation, there is a simple way to test it for KS-continuity as follows.

\begin{Proposition}
A generalized translation-invariant valuation $\psi \in \Val_k^{+,-\infty}(V)$ belongs to $\Val_k^{+,\KS}(V)$ if and only if its (generalized) Klain function is continuous. 
\end{Proposition}

\proof
By the extension of the Klain embedding to generalized valuations, we immediately see that if $\phi_j \to \phi$ in $\Val^{+, \KS}$ then $\Kl_\phi$ coincides with the weak limit of $\Kl_{\phi_j}$, which is therefore also a uniform limit. In particular, it is continuous. Conversely, if $\psi \in \Val_k^{+,-\infty}$ has a continuous Klain function, then we may approximate $\psi$ by $\psi_\epsilon=\rho_\epsilon \ast \psi$, where $\rho_\epsilon$ is an approximate identity on $\mathrm{SO}(n)$. Then $\psi_\epsilon \in \Val^{+,\infty}_k$ and $\|\Kl(\psi_\epsilon)- \Kl(\psi)\|_{\infty} \to 0$. 
\endproof

A similar statement also holds for general valuations using the Schneider embedding, although it is not as easy to apply in practice. 

\begin{Proposition}
A generalized translation-invariant valuation $\phi \in \Val_k^{-\infty}(V)$ belongs to $\Val_k^{\KS}(V)$ if and only if $\Sc(\phi) \in  \Gamma(\Gr_{k+1}(V), \Val_k(V))$.
\end{Proposition}

\proof 
Assume $\phi \in \Val^{\KS}$. Then $\phi=\lim \phi_j$ in $\Val^{\KS}$ with smooth $\phi_j$, and it follows that $\Sc \phi=\lim \Sc_{\phi_j} \in \Gamma(\Gr_{k+1}(V), \Val_k(V))$.

In the other direction, if $\Sc(\phi) \in  \Gamma(\Gr_{k+1}(V), \Val_k(V))$ write $\Sc(\phi \ast { \rho_\epsilon})=\Sc\phi \ast { \rho_\epsilon}$ for an approximate identity $\rho_\epsilon \in C^\infty(\mathrm{SO}(V))$.
Then $\phi\ast \rho_\epsilon \in \Val^\infty(V)$, and it converges to $\Sc\phi$ in the norm $\|\bullet\|_{\Sc}$.

\endproof

\subsection{Functorial properties}

Since any linear map is the composition of a monomorphism and an epimorphism, it is enough to prove Theorem \ref{mainthm_klaincont} only for such maps.

\begin{Proposition}  \label{prop_restriction_vals}
Let $V$ be a real vector space of dimension $n$, $j:W \hookrightarrow V$ the inclusion of a subspace. The restriction of continuous translation-invariant valuations $j^*:\Val_k(V)\to \Val_k(W)$ extends by continuity to KS-continuous valuations: $j^*:\Val_k^{\KS}(V)\to \Val_k^{\KS}(W)$. It holds that
\begin{displaymath}
\Sc_{j^*\phi}(E)=\Sc_\phi(jE), \quad E \in \Gr_{k+1}(W).
\end{displaymath}
In the even case, 
\begin{displaymath}
\Kl_{j^*\phi}=j^*\Kl_\phi, \quad \phi \in \Val_k^{+,\KS}(V).
\end{displaymath}
\end{Proposition}

\proof
Choose $\psi\in \Val^{\KS}(V)$. Let us show the existence of a valuation $j^*\psi \in  \Val^{\KS}(W)$ satisfying $\Sc_{j^*\psi}=j^*\Sc_\psi$. Let $(\nu_i)_i$ be a sequence of smooth probability measures on $\GL(V)$, whose supports shrink to the identity element. Then the generalized valuation 
\begin{displaymath}
\psi_i:=\int_{\GL(V)} g\psi d\nu_i(g) \in \Val^{-\infty}_k(V)
\end{displaymath}
is smooth. Since $\Sc_{\psi_i}=\int_{\mathrm{GL}(V)} g^* \Sc_\psi d\nu_i(g)$ and $\Sc_\psi$ is continuous, we have that $\Sc_{\psi_i} \to \Sc_{\psi}$ uniformly. Then $j^*\Sc_{\psi_i}$ converges uniformly to $j^*\Sc_\psi$, so $j^*\psi_i$ converges in $\Val^{\KS}(W)$, and the limit is taken to be $j^*\psi$.

The continuity of $j^*:\Val^{\KS}(V)\to \Val^{\KS}(W)$ is now obvious. The statement in the even case can be shown similarly. 
\endproof

\begin{Proposition}\label{prop:pull_by_epi}
Let $\pi:V \to W$ be a surjective map between vector spaces. Then $\pi^*:\Val_k( W) \to \Val_{k}( V)$ extends by continuity to KS-continuous valuations: $\pi^*:\Val_k^{\KS}( W) \to \Val_{k}^{\KS}( V)$. For $E \in \Gr_{k+1}(V)$ we write $\pi_E:E \to \pi E$ for the restriction of $\pi$; $\pi_E^*:\Val_k(\pi E) \to \Val_k(E)$ for the pull-back. 
It holds for $\psi \in \Val_k^{\KS}(W)$ that 
\begin{displaymath}
\Sc_{\pi^*\psi}(E)=\pi_E^*(\Sc_\psi(\pi E)) \in \Val_k(E), E \in \Gr_{k+1}(V),
\end{displaymath}
whenever $E \cap \Ker \pi=0$, and  $\Sc_{\pi^* \psi}(E)=0$ otherwise.

In the even case,
\begin{displaymath}
\Kl_{\pi^* \psi}(E)=\begin{cases} \pi_E^* \Kl_\psi(\pi E) & \text{if } E \cap \Ker \pi=0 \\ 0 & \text{otherwise} \end{cases}, 
\end{displaymath}
where $E \in \Gr_k(V), \psi \in \Val_k^{+,\KS}(W)$.
\end{Proposition}

\proof
If $E \cap \Ker \pi=0$, we have $\pi(E) \in \Gr_{k+1}(W)$ and the following diagram commutes
\begin{displaymath}
\xymatrix{
	E \ar[r]^{i_E} \ar[d]^{\pi_E} & V \ar[d]^\pi\\
	\pi(E) \ar[r]^{i_{\pi E}} & W}
\end{displaymath}

It follows from \cite[Prop. 3.1.2]{alesker_fourier} that on $\Val_k(W)$ we have 
\begin{displaymath}
i_E^* \circ \pi^*=\pi_E^* \circ i_{\pi E}^*
\end{displaymath}
The rest of the proof follows by approximation as before.
\endproof

\begin{Proposition}
Let $\pi:V \to W$ be a surjective map between vector spaces and set $l:=\dim V-\dim W$. Then $\pi_*:\Val_k(V) \otimes \Dens^*(V)\to \Val_{k-l}(W) \otimes \Dens^*(W)$ extends by continuity to KS-continuous valuations: $\pi_*:\Val_k^{\KS}(V) \otimes \Dens^*(V)\to \Val_{k-l}^{\KS}(W) \otimes \Dens^*(W)$.

For $E \in \Gr_{k-l+1}(W)$, we write $\pi_E:\pi^{-1} E \to E$ for the restriction and $(\pi_E)_*:\Val_k(\pi^{-1}E) \otimes \Dens^*(\pi^{-1}E) \to \Val_{k-l}(E) \otimes \Dens^*(E)$ for the push-forward by $\pi_E$. It holds for $\phi \in\Val_k^{\KS}(V)\otimes\Dens^*(V)$ that 
\begin{displaymath}
\Sc_{\pi_* \phi}(E)=(\pi_E)_*(\Sc_{\phi}(\pi^{-1}E)) \in \Val_{k-l}(E) \otimes \Dens(W^*),\quad E \in \Gr_{k-l+1}(W).
\end{displaymath}
In the even case, 
\begin{displaymath}
\Kl_{\pi_* \psi}(E)=\Kl_\psi(\pi^{-1}E), \quad \psi \in\Val_k^{+,\KS}(V)\otimes\Dens^*(V), E \in \Gr_{k-l}(W). 
\end{displaymath}
\end{Proposition}

\proof
Let us first check that the equations formally make sense. Since $\pi_*\phi \in \Val_{k-l}(W) \otimes \Dens(W^*)$, $\Sc_{\pi_* \phi}(E)$ is an element of $\Val_{k-l}(E) \otimes \Dens(W^*)$. On the other hand, $\pi^{-1}E \in \Gr_{k+1}(V)$ and $\Sc_{\phi}(\pi^{-1}E) \in \Val_k(\pi^{-1}E) \otimes \Dens(V^*)$. Then $(\pi_E)_* \Sc_{\phi}(\pi^{-1}E) \in \Val_{k-l}(E) \otimes \Dens^*(E)/\Dens^*(\pi^{-1}E) \otimes \Dens^*(V)$. But the last factor equals $\Dens^*(W)$. Similarly, for the second equation we use the natural isomorphism $\Dens(E) \otimes \Dens^*(W) \cong \Dens(\pi^{-1}E) \otimes \Dens^*(V)$.

Next, we prove the statement.  For simplicity, let us omit in the following the various twists by densities. Take a sequence of continuous valuations $\phi_i$ such that $\Sc_{\phi_i} \to \Sc_{\phi}$ uniformly. Then for each fixed $E \in \Gr_{k-l+1}(W)$, $j_{\pi^{-1}E}^*\phi_i \to j_{\pi^{-1}E}^*\phi$ uniformly in $E$. It follows from \cite[Thm 3.5.2]{alesker_fourier}, applied to the following diagram:
\begin{displaymath}
\xymatrix{
	\pi^{-1}E \ar[r]^-{j_{\pi^{-1}E}} \ar[d]^{\pi_E} & V \ar[d]^{\pi} \\
	E \ar[r]^{j_E} & W}
\end{displaymath}
that $j_E^*\pi_*\phi_i=(\pi_E)_* j^*_{\pi^{-1}E}\phi $, that is, $\Sc_{\pi_*}\phi(E)=(\pi_E)_*\phi|_{\pi^{-1}(E)}$.

Since the pushforward $(\pi_E)_*$ is continuous in the Banach space topology and has norm that is independent of $E$, we conclude that $\Sc_{\pi_*\phi_i} (E)= j_E^*\pi_*\phi_i$ converges uniformly to $(\pi_E)_* j^*_{\pi^{-1}E}\phi$. We thus define $\pi_*\phi$ as the limit in $\Val_{k-l}^{\KS}(V)$ of $\pi_*\phi_i$.

\endproof

\begin{Proposition}
Let $j:W \to V$ be an inclusion of vector spaces and set $l:=\dim V-\dim W$. Then $j_*:\Val_k(W) \otimes \Dens^*(W)\to \Val_{k+l}(V) \otimes \Dens^*(V)$ extends by continuity to KS-continuous valuations: $j_*:\Val_k^{\KS}(W) \otimes \Dens^*(W)\to \Val_{k+l}^{\KS}(V) \otimes \Dens^*(V)$. For $\phi \in \Val^{\KS}_k(W) \otimes \Dens^*(W)$ we have 
\begin{equation} \label{eq_push_forward_inclusion}
\Sc_{j_*\phi}(E)=(j_E)_* \Sc_{\phi}(E \cap W), \quad E \in \Gr_{k+l+1}(V),
\end{equation}
if $\dim E \cap W=k+1$ and $\Sc_{j_*\phi}(E)=0$ if $\dim E \cap W>k+1$. Here $j_E:E \cap W \to E$ is the inclusion.

In the even case,
\begin{equation}\label{eq_push_forward_inclusion_even}
\Kl_{j_*\phi}(E)=\begin{cases} \Kl_\phi(E \cap W) |\langle E,V/W\rangle| & \text{if } \dim E \cap W=k,\\ 0 & \text{if } \dim E \cap W > k,\end{cases}
\end{equation}
where $E \in \Gr_{k+l}(V)$. If $\dim E \cap W=k$,  then  $E+W=V$ and hence $V/W \cong E/E\cap W$. The element $ |\langle E, V/W\rangle| \in \Dens^*(V/W) \otimes \Dens (E/(E\cap W)) \cong \Dens^*(V/W)\otimes \Dens (V/W)$ is the canonic element in this one-dimensional space. 
\end{Proposition}

\proof
We have the following commuting diagram
\begin{displaymath}
\xymatrix{
E \cap W \ar[r]^-{i_{E \cap W}} \ar[d]^{j_E} & W \ar[d]^j\\
E \ar[r]^{i_E} & V}
\end{displaymath}
It follows that $\Dens(W) \otimes \Dens(E) \cong \Dens(E \cap W) \otimes \Dens(V)$, which may be used to see that \eqref{eq_push_forward_inclusion_even} formally makes sense.
For $\phi \in \Val_k(W) \otimes \Dens^*(W)$ we have 
\begin{displaymath}
\Sc_{j_*\phi}(E)=i_E^*j_*\phi=(j_E)_* i_{E \cap W}^*\phi=(j_E)_* \Sc_{\phi}(E \cap W).
\end{displaymath}
In the last line, we used again \cite[Thm 3.5.2]{alesker_fourier}. The proof is concluded by approximation.
\endproof

The even KS-continuous valuations possess the additional highly useful property of being invariant under the Alesker-Fourier transform.

\begin{Proposition} \label{prop:fourier_klain_continuous}
	Let $\mathbb F: \Val^{+,-\infty}(V) \to \Val^{+,-\infty}(V^*) \otimes \Dens(V)$ be the Alesker-Fourier transform on even generalized translation-invariant valuations \cite[Section 6.2]{alesker_faifman}. Then $\mathbb F$ restricts to an isomorphism of topological vector spaces $\mathbb F:\Val^{+,\KS}(V) \to \Val^{+,\KS}(V^*) \otimes \Dens(V)$.  
\end{Proposition}

\proof
Choosing a Euclidean scalar product to identify $V \cong V^*, \Dens(V) \cong \C$, the (generalized) Klain function of $\mathbb F \psi$ corresponds to the (generalized) Klain function of $\psi$ through the orthogonal complement map $\Gr_k(V)\to\Gr_{n-k}(V)$. It is now obvious that the image of a KS-continuous valuation is KS-continuous, and that the restriction of $\mathbb F$ to the KS-continuous valuations is continuous.
\endproof

\subsection{The Klain function of an even KS-continuous valuation}

Fix $k$ and let $m_\phi \in \Gamma^\infty( C^{n,k})$ be a smooth Crofton measure defining a valuation $\phi \in \Val_k^{+,\infty}$. Since the Klain function of $\phi$ is the cosine transform $T_{n-k,k}m_\phi$, we have, for all $E_0 \in \Gr_k(V)$, by \eqref{eq_cosine_delta},
\begin{displaymath}
\Kl_\phi(E_0)=\left\langle T_{n-k,k}m_\phi,\delta_{E_0}\right\rangle=\left\langle m_\phi,T_k\delta_{E_0}\right\rangle=\left\langle m_\phi, \langle E_0,\bullet\rangle \right\rangle,
\end{displaymath}
i.e. we integrate the continuous function $\langle E_0,\bullet\rangle$ against the smooth measure $m_\phi$.

We will need the following generalization of this equation which applies to KS-continuous, even translation-invariant valuations. 

It will relax the smoothness assumption by considering the singular support, or, more generally, the wavefront of the Crofton measure. For the definitions and basic facts about wavefronts, see for example \cite{grigis_sjostrand,guillemin_sternberg77}. 

We shall now describe the wavefront of the cosine transform.  Let $\langle E, F\rangle\in \Gamma^{-\infty}(\Gr_k(V)\times \Gr_{n-k}(V),\Dens(E)\otimes \Dens^*(V/F))$ denote the Schwartz kernel of the $\mathrm{GL}(V)$-equivariant cosine transform \[T_{n-k,k}:\Gamma^\infty(\Gr_{n-k}(V),C^{n,k})\to \Gamma^\infty(\Gr_k(V),K^{n,k})\] Note that $\langle E,F\rangle$ is continuous everywhere and smooth outside the zero set
\begin{displaymath}
Z:=\{(E,F):  E\cap F\neq \{0\}\} \subset \Gr_k(V) \times \Gr_{n-k}(V).
\end{displaymath} 

Observe that $Z$ is stratified by the $\mathrm{SL}(V)$-orbits on $\Gr_k(V) \times \Gr_{n-k}(V)$, which are locally closed submanifolds classified by $\dim (E\cap F)$. By $N^*Z$ we denote the union of the conormal bundles of the strata. 

\begin{Lemma}\label{lem:cosine_wavefront}
We have $\WF(\langle E, F\rangle)\subset N^*Z$.
\end{Lemma}

\proof
The singular support of $\langle E, F\rangle$ is clearly contained in $Z$, as it is smooth elsewhere. Fix $z=(E, F)\in Z$. As $\langle E,F\rangle$ is $\mathrm{SL}(V)$-invariant while $\mathrm{SL}(V)$ is semisimple, it follows from Corollary \ref{cor:invariant_wavefront} that $\WF(\langle E,F\rangle)\cap T_z^* (\Gr_k(V) \times \Gr_{n-k}(V)) \subset N_z^*Z$, concluding the proof.
\endproof

\begin{Proposition} \label{prop_Klain_continuous_section_computation}
Let $\phi \in \Val_k^{+,\KS}(V)$, and let $m_\phi \in \Gamma^{-\infty}(C^{n,k})$ be a generalized Crofton measure for $\phi$. Let ${E_0} \in  \Gr_k(V)$ be such that ${E_0}^\cap \cap S_\phi=\emptyset$, where ${E_0}^\cap=\{F\in \Gr_{n-k}(V): E_0 \cap F \neq 0\}$ and $S_\phi$ is the singular support of $m_\phi$. 

Then 
\begin{displaymath}
\Kl_\phi({E_0})=\left\langle m_\phi, \langle {E_0},\bullet\rangle \right\rangle.
\end{displaymath}

Moreover, the same conclusion holds under the weaker assumption that the wavefront of $m_\phi$ satisfies $\WF(m_\phi) \cap N^*{E_0}^\cap=\emptyset$. Here $N^*A \subset T^* \Gr_{n-k}(V)$ is the union of the conormal bundles of the strata of a subset $A \subset \Gr_{n-k}(V)$, which is stratified by finitely many locally closed submanifolds.
\end{Proposition}

\proof
We fix for simplicity a Euclidean structure on $V$ to trivialize all the line bundles involved.

Take a sequence of functions $\delta_j \in C^\infty (\Gr_k(V))$ with supports shrinking to ${E_0}$ s.t. $\delta_j \to \delta_{{E_0}}$ as measures.
Then
\begin{align*}
\Kl_\phi({E_0}) & =\langle \Kl_\phi, \delta_{{E_0}} \rangle\\
& =\lim_{j\to\infty} \langle \Kl_\phi,\delta_j\rangle\\
& =\lim_{j\to\infty} \langle T_k(m_\phi), \delta_j\rangle\\
& = \lim_{j\to\infty} \langle m_\phi, T_k(\delta_j)\rangle. 
\end{align*}

Using the assumption of disjoint singular supports, take a smooth function $\rho \in C^\infty(\Gr_{n-k}(V))$ which is identically $1$ in a neighborhood of ${E_0}^\cap$ and identically $0$ in a neighborhood of $S_\phi$. 

Write 
\begin{displaymath}
 \langle m_\phi, T_k(\delta_j)\rangle=\langle \rho m_\phi, T_k(\delta_j)\rangle+\langle m_\phi, (1-\rho)T_k(\delta_j)\rangle.
\end{displaymath}

Since $\rho m_\phi$ is a smooth measure, and since $T_k(\delta_j) \to \langle \bullet,{E_0}\rangle$ weakly, the first summand converges to $\left\langle \rho m_\phi,  \langle {E_0},\bullet\rangle\right\rangle$. 

Next, $(1-\rho)T_k(\delta_j)$ converges by \eqref{eq_cosine_delta} to $(1-\rho)\langle {E_0},\bullet\rangle$ in $C^\infty (\Gr_{n-k}(V))$, so the second summand converges to $\left\langle m_\phi, (1-\rho)\langle {E_0},\bullet\rangle\right\rangle$. 
	
We pass to the general statement concerning disjoint wavefronts. We will write $T^+M$ for $T^*M\setminus 0_M$. For a closed cone $\Gamma \subset T^+\Gr_k(V)$, $C_\Gamma^{-\infty}(\Gr_k(V))$ denotes the space of generalized functions with wavefront contained in $\Gamma$, with the standard locally convex topology. Write $\Gamma_0=T^+_{E_0}\Gr_k(V) \subset T^+ \Gr_k(V)$. We choose $\delta_j \in C^\infty(\Gr_k(V))$ such that $\delta_j \to \delta_{E_0}$ in the topology of $C^{-\infty}_{\Gamma_0}(\Gr_k(V))$. 

Denote $X=\Gr_k(V) \times \Gr_{n-k}(V)$ and $Z=\{(E,F):  E\cap F\neq \{0\}\} \subset X$. The group $\{g\in \mathrm{GL}(n): gE_0=E_0\}$ acts on $E_0^\cap\subset \Gr_{n-k}(V)$ with finitely many orbits classified by $\dim (E_0\cap F)$, which are locally closed submanifolds of $\Gr_{n-k}(V)$. 

By Lemma \ref{lem:cosine_wavefront}, $\WF(T_{k,n-k})=\WF(\langle F,E\rangle )\subset N^*Z$. 
 
Recall that for two cones $C_1 \subset T^*(\Gr_k \times \Gr_{n-k}), C_0\subset T^*(\Gr_k)$, their composition is given by  
\begin{displaymath}
 C_1 \circ C_0:=\{(F, \eta) \in T^+\Gr_{n-k} : \exists (E, \xi) \in C_0 \mbox{ s.t. } (E, \xi, F, \eta) \in C_1 \}.
\end{displaymath}
We will also need the notation
\begin{displaymath}
 C_1'= \{ (E, \xi, F, \eta): (E, -\xi, F, \eta)\in C_1\}.
\end{displaymath}
There is a natural inclusion $E_0^\cap \subset Z, F \mapsto (E_0,F)$. If $(\xi,\eta)\in N^*_{(E_0,F)}Z\subset T^*_{E_0}\Gr_k\oplus T^*_{F}\Gr_{n-k}$ then by restriction to $T_{F}E_0^\cap$ we see that $\eta\in N^*E_0^\cap$. It follows that
\begin{displaymath}
 \Gamma_1:=\WF(\langle F,E\rangle )'\circ \Gamma_0\subset (N^*Z)'\circ  \Gamma_0 \subset N^*E_0^\cap
\end{displaymath}

Finally, since the projection $\pi:Z\to \Gr_{n-k}(V)$ is submersive on each stratum of $Z$, we have 
\begin{displaymath}
(N^*Z)'\circ  0_{\Gr_k}=\emptyset. 
\end{displaymath}

Therefore  (see \cite[Theorem 7.8]{grigis_sjostrand}), 	$T_k$ extends to a sequentially continuous operator $T_k: C^{-\infty}_{\Gamma_0}(\Gr_k(V)) \to C^{-\infty}_{\Gamma_1}(\Gr_{n-k}(V))$.
By the sequential continuity of $T_k$, and since $\WF(m_\phi) \cap \Gamma_1 = \emptyset$,
\begin{displaymath}
\Kl_\phi(E_0)=\lim_{j\to\infty} \langle m_\phi, T_k(\delta_j)\rangle=\langle m_\phi, T_k(\delta_{E_0})\rangle=\langle m_\phi,\langle \bullet,E_0\rangle\rangle.
\end{displaymath}
This concludes the proof.
\endproof

\section{Geometry of the orbits}
\label{sec_orbits}

\subsection{The Grassmannian under the action of $\mathrm{O}(p,q)$}

Let $V$ be a vector space of dimension $n$, $Q$ a non-degenerate quadratic form on $V$ of signature $(p,q)$ and $G:=\mathrm{O}(Q)=\mathrm{O}(p,q)$. The $Q$-orthogonal complement of a subspace $E$ will be denoted by $E^Q$.

Denote by $X^k_{a,b}\subset \Gr_k(V)$ the subset consisting of those subspaces $E$ for which the restriction of $Q$ to $E$ has signature $(a,b)$. It follows by Witt's theorem \cite[Thm. 5.3.]{scharlau85} that whenever $X^k_{a,b}$ is non-empty, it is a $G$-orbit in $\Gr_k(V)$. It is non-empty precisely for those pairs $(a,b)$ for which $a+b \leq k$, $\max(0,k-q) \leq a \leq p$, $\max(0,k-p) \leq b \leq q$ (Lemma \ref{lemma_signatures}). The open orbits are those for which $a+b=k$. There is a unique closed orbit $X^k_0:=X^k_{a_0,b_0}$ with $a_0:=\max(0,k-q), b_0:=\max(0,k-p)$.

\begin{Lemma} \label{lem:antisymmetricextension}
Let $Q$ be a non-degenerate quadratic form on $V$, $E \subset V$ any subspace, and $E_0:=E \cap E^Q$. Let $T:E \to V/E$ be a linear map. Then $Q(Tx,x)=0$ for all $x \in E_0$ if and only if there is a lift $T':V \to V$  of $T$ s.t. $Q(T'x,x)=0$ for all $x \in V$, i.e. $T' \in \mathfrak{so}(Q)$.
\end{Lemma}

\proof
Let $T \in \Hom(E,V/E)$ such that $Q(Tx,x)=0$ for all $x \in E_0$.

Choose any subspaces $E' \subset E$, $E'' \subset E^Q$ with $E=E_0 \oplus E'$, $E^Q=E_0 \oplus E''$. 

Let $r:=\dim E_0$. Then $\dim E'=k-r$ and $\dim E''=(n-k)-r$. Then $E' \oplus E'' \subset V$ is a non-degenerate subspace of dimension $n-2r$, and $E_0$ is an isotropic subspace of the non-degenerate space $W:=(E' \oplus E'')^Q$ of dimension $2r$.  

By Lemma \ref{lemma_isotropic_complement}, one can fix a $Q$-isotropic subspace $F$ s.t. $W=E_0 \oplus F$. Then $V=E' \oplus E'' \oplus E_0 \oplus F$.  

Let $T_1 \in \Hom(E,V)$ be a lift of $T$ such that $T_1(E') \subset E'' \oplus F$, $T_1(E_0) \subset E'' \oplus F \oplus E'$ and 
\begin{displaymath}
Q( \pi'(T_1x),e):=-Q(x,Te), \quad x \in E_0, e \in E',
\end{displaymath}
where $ \pi'(T_1x)$ is the $E'$-component of $T_1x$.

Note that $Q$ gives identifications $E^*=E' \oplus F$, as well as $(E \oplus E'')^*=E' \oplus E'' \oplus F$. We extend $T_1$ to a map $T_2 \in \Hom(E \oplus E'',V)$ by requiring $T_2(E'') \subset E' \oplus F$ and 
\begin{displaymath}
 Q(T_2e'',e):=-Q(e'',T_1e), \quad e'' \in E'', e \in E.
\end{displaymath}

Note that $Q(T_2e_1,e_2)=0$ for all $e_1,e_2 \in  E''$. 

Finally, we extend $T_2$ to a map $T' \in \Hom(V,V)$ by requiring $T'(F) \subset E' \oplus E'' \oplus F$ and  
\begin{displaymath}
 Q(T'f,x)=-Q(f,T_2x), \quad x \in E \oplus E'', f \in F.
\end{displaymath}
Again we have $Q(T'f,f')=0$ for $f,f' \in F$. 

Then $T' \in \mathfrak{so}(Q)$, and $T'$ lifts $T$ as required.
\endproof

\begin{Proposition}\label{prop_EquivariantFormOfNormalBundle}
If $E \in X^k_{a,b}$, $\Stab(E) \subset G$ its stabilizer and $E_0:=E \cap E^Q$, then 
 the normal space $N_E X^k_{a,b}=T_E \Gr_k(V)/T_E X^k_{a,b}$ is $\Stab(E)$-isomorphic to $\Sym^2 E_0^*$. In particular, if $X^k_{a,b}$ is non-empty, then 
\begin{displaymath} 
\dim X^k_{a,b}=k(n-k)-\binom{r+1}{2},
\end{displaymath}
where $r:=\dim E_0=k-a-b$.
\end{Proposition}

\proof
Recall that $T_E\Gr_k(V)\cong\Hom(E,V/E)$.

Consider the $\Stab(E)$-equivariant map $\pi_0: \Hom(E,V/E) \to E_0^* \otimes E_0^*$, given by 
\begin{displaymath}
 \pi_0(T)(u,v)=Q(Tu,v), \quad u,v \in E_0,
\end{displaymath}
which is easily seen to be well-defined and onto.  

Writing $\pi_E: \mathfrak{gl}(V)\to \Hom(E,V/E)$ for the natural projection, it follows from  Lemma \ref{lem:antisymmetricextension} that  
\begin{displaymath}
T_E X_{a,b}^k\cong\pi_E(\mathfrak{so}(Q))=\pi_0^{-1}\left(\largewedge^2 E_0^*\right) .
\end{displaymath}
Thus
\begin{align*}
 N_E X^k_{a,b} &=T_E\Gr_k(V)/T_EX^k_{a,b} \cong \Hom(E,V/E)/\pi_0^{-1}\left(\largewedge^2E_0^*\right) \\&\cong \pi_0\Hom(E,V/E)/\largewedge^2E_0^*
 = E_0^* \otimes E_0^*/\largewedge^2E_0^*  \cong \Sym^2E_0^*,
\end{align*}
where we have used the fact that a linear map $\pi_0:U_1\to U_2$ induces a natural isomorphism $U_1/\pi_0^{-1}(W)\cong \pi_0(U_1)/W$ for any subspace $W\subset \pi_0(U_1)$.
\endproof

\subsection{$\mathrm{O}(p,q)$-invariant sections of the Klain bundle}

\begin{Proposition} \label{prop_invariant_sections_klain}
 Let $V$ be an $n$-dimensional vector space with a non-degenerate bilinear form $Q$ and corresponding orthogonal group $G$. Then the dimension of the space of $G$-invariant generalized sections of the Klain bundle $K^{n,k}$ equals the number of open $G$-orbits in $\Gr_k(V)$. A basis is given by the sections    
\begin{displaymath}
 \kappa_a(E)(v_1 \wedge \ldots \wedge v_k):=\begin{cases} \left|\det (Q(v_i,v_j))_{i,j=1}^k \right|^\frac12 & E=\mathrm{span}\{v_1,\ldots,v_k\} \in X^k_{a,k-a},\\0 & \text{otherwise,}\end{cases}
\end{displaymath}
with $\max(0,k-q) \leq a \leq \min(k,p)$. In particular, each $G$-invariant generalized section of $K^{n,k}$ is continuous. 
\end{Proposition}

\proof
First we show that there are no sections supported on the complement of the open orbits, denoted $ Z \subset \Gr_k(V)$. 
Let $Y=X^k_{a,b}\subset Z$ be any orbit with $r:=k-a-b>0$. Define for every $\alpha \geq 0$ the $G$-module
\begin{displaymath}
 F^\alpha_E:=\Sym^\alpha(N_E Y) \otimes \Dens^*(N_E Y) \otimes K^{n,k}|_E, \quad  E \in Y,
\end{displaymath}
where $N_E Y:=T_E  \Gr_k(V)/T_E Y$ is the normal space of $Y$ at $E$.

By Lemma \ref{lem:LocallyClosedOrbits}, it suffices to check that for all $\alpha \geq 0$, the stabilizer in $G$ of $E$ has no non-trivial invariants in $F^\alpha_E$.

Define $E_0:=E \cap E^Q$. By Proposition \ref{prop_EquivariantFormOfNormalBundle}, there is a $G$-equivariant isomorphism $N_E Y=\Sym^2 E_0^*$, and $\Stab(E)$ acts on $E_0$ as $\GL(E_0)$ (non-faithfully). Note that $K^{n,k}|_E=\Dens(E)=\Dens(E_0) \otimes \Dens(E/E_0)$, and since $E/E_0$ inherits a non-degenerate quadratic form, $G$ acts trivially on $\Dens(E/E_0)$. Taking an element $g \in \Stab(E)$ acting on $E_0$ by the scalar $\lambda \neq 1$, $g$ acts on $F^\alpha|_E$ by $\lambda^{-2\alpha}\lambda^{-2\frac{r(r+1)}{2}}\lambda^{-r}\neq 1$. Therefore, $F^\alpha_E$ admits no $\Stab(E)$-invariants, as required.

Since $K^{n,k}$ is a $1$-dimensional bundle over $\Gr_k(V)$, the space of $G$-invariant sections over any open orbit is at most $1$-dimensional. It is clear that the function $\kappa_a$, as defined in the statement of the proposition, is an invariant continuous section of the Klain bundle. This finishes the proof.  
\endproof

\begin{Corollary}
 Even generalized translation-invariant valuations which are invariant under $\mathrm O(p,q)$ are KS-continuous. 
\end{Corollary}

\begin{Corollary}\label{cor:injectivity_of_restriction}
Fix an isometric inclusion $i:\mathbb R^{p',q'}\to \mathbb R^{p,q}$.
The restriction map $i^*: \Val_k^{+,-\infty}(\mathbb R^{p,q})^{\mathrm O(p,q)}\to \Val_k^{+,-\infty}(\mathbb R^{p',q'})^{\mathrm O(p',q')}$ is injective whenever $\min(k,p)=\min(k,p')$ and $\max(0,k-q)=\max(0,k-q')$. 

\end{Corollary}

\proof
This follows from the injectivity of the Klain map, Proposition \ref{prop_invariant_sections_klain} and Proposition \ref{prop_restriction_vals}.
\endproof


\subsection{The differential of the $\mathrm{O}(p,q)$-action on $\Gr_k(\R^n)$}
 Let us fix a Euclidean form $P$ compatible with $Q$.

\begin{Definition}\label{def:theta}
 Let $k \leq \frac{n}{2}$. 
We define the functions $\theta:\Gr_k(\R^n) \to [0,\frac\pi2]$ by setting
\begin{displaymath}
 \cos 2\theta(E):=\det M(E), 
\end{displaymath}
 where $M(E)=(Q(u_i,u_j))$ for an arbitrary $P$-orthonormal basis $u_j$ of $E$.
\end{Definition} 

\begin{Proposition} \label{prop_jacobian}
For $g \in \mathrm{GL}(n)$  and $E \in \Gr_k(\R^n)$, let 
\begin{displaymath}
 \psi_g(E):=\frac{1}{\Jac(g:E \to gE)^2},
\end{displaymath} 
where $E$ and $gE$ are endowed with the $P$-induced Euclidean scalar product. Clearly this function is smooth and positive. 
\begin{enumerate}
 \item If $E$ is non-degenerate with respect to $Q$ and $g\in \mathrm O(Q)$ then 
 \begin{displaymath}
  \psi_g(E)=\frac{\cos 2 \theta(gE)}{\cos 2 \theta(E)}.
 \end{displaymath}
\item $\psi_g \equiv 1$ for $g \in \mathrm O(n)$. 
\item 
\begin{displaymath}
 \left|\Jac(g:\Gr_k(\R^n) \to \Gr_k(\R^n))|_E\right|=\psi_g(E)^\frac{n}{2}|\det g|^k.
\end{displaymath}
Here the Jacobian is computed with respect to any $\mathrm{O}(n)$-invariant Riemannian metric on $\Gr_k(\R^n)$. 
\end{enumerate}
\end{Proposition}

\proof
\begin{enumerate}
 \item Let $E$ be non-degenerate with respect to $Q$. Let $f_1,\ldots,f_k$ be any basis of $E$. Then 
 \begin{displaymath}
  \cos 2\theta(E)=\frac{\det (Q(f_i,f_j)_{i,j})}{\det (P(f_i,f_j)_{i,j})}, \quad \cos 2\theta(gE)=\frac{\det (Q(gf_i,gf_j)_{i,j})}{\det (P(gf_i,gf_j)_{i,j})}.
 \end{displaymath}
Since $g \in \mathrm{O}(p,q)$, we have $\det (Q(f_i,f_j)_{i,j})=\det (Q(gf_i,gf_j)_{i,j})$ and therefore 
\begin{displaymath}
 \frac{\cos 2\theta(gE)}{\cos 2 \theta(E)}=\frac{\det (P(f_i,f_j)_{i,j})}{\det (P(gf_i,gf_j)_{i,j})}=\psi_g(E).
\end{displaymath}
\item Trivial.
\item By composing with an element of $\mathrm{SO}(n)$ if necessary, we may assume that $gE=E$. 
The tangent space is given by 
\begin{displaymath}
 T_E \Gr_k(\R^n)=\mathrm{Hom}(E,V/E)=E^* \otimes V/E. 
\end{displaymath}
Therefore 
\begin{displaymath}
 \Dens( T_E \Gr_k(\R^n))=\Dens^{n-k}  (E^*)\otimes \Dens^k (V/E)=\Dens^{n}  (E^*)\otimes \Dens^k (V).
\end{displaymath}
By definition, $g$ acts by the scalar $\psi_g(E)^{-\frac12}$ on $\Dens(E)$, and hence by the scalar $\psi_g(E)^{\frac{n}{2}}$ on $\Dens^{n} (E^*)$. Evidently $g$ acts by the scalar $|\det g|^k$ on $\Dens^k( V)$. It follows that $g$ acts by the scalar $|\det g|^k \psi_g(E)^\frac{n}{2}$ on $\Dens(T_E \Gr_k(\R^n))$, as claimed.
\end{enumerate}
\endproof

As $|\det(g)|=1$ for $g \in \mathrm{O}(p,q)$ (see equation \ref{eq_densities_trivial}), we get the following  corollary.

\begin{Corollary} \label{cor_invariant_crofton_klain_sections}
 \begin{enumerate}
  \item Using the Euclidean trivialization, an $\mathrm{O}(p,q)$-invariant section of the Crofton bundle over $\Gr_k(\R^n)$ corresponds to a generalized
function $f\in C^{-\infty}(\Gr_k(\R^n))$ transforming by 
\begin{displaymath}
 g^*(f)=\psi_{g}^{-(n+1)/2} \cdot f
\end{displaymath}
for all $g \in \mathrm{O}(p,q)$. 
\item An $\mathrm{O}(p,q)$-invariant section of the Klain bundle over $\Gr_k(\R^n)$ corresponds to a generalized
function $f \in C^{-\infty}(\Gr_k(\R^n))$ transforming by 
\begin{displaymath}
 g^*(f)=\psi_{g}^{1/2} \cdot f
\end{displaymath}
for all $g \in \mathrm{O}(p,q)$. 
 \end{enumerate}
\end{Corollary}

\subsection{Orbit space on $\mathbb{P}_+(V^*)$ under $\mathrm{SO}^+(p,q)$}

By $\mathbb{P}_+(V^*):=(V^* \setminus \{0\})/\R_{>0}=\Gr_1^+(V^*)$ we denote the space of oriented lines in $V^*$, which is the same as the space of co-oriented hyperplanes in $V$.

\begin{Proposition} \label{prop_orbits_lines}
 \begin{enumerate}
  \item If $\min(p,q)=0$, then the action of $\mathrm{SO}^+(p,q)$ on $\mathbb{P}_+(V^*)$ has one (open) orbit. 
  \item If $\min(p,q)=1$ and $n>2$, then the action of $\mathrm{SO}^+(p,q)$ on $\mathbb{P}_+(V^*)$ has $3$ open orbits $M_1^-, M_2^-, M^+$ and two closed orbits $M_1^0$, $M_2^0$ (the union of which we denote by $M^0$).
  \item If $p=q=1$, then the action of $\mathrm{SO}^+(p,q)$ on $\mathbb{P}_+(V^*)$ has four open orbits $M_1^-,M_2^-,M_1^+,M_2^+$ and four closed orbits $M_{++}^0$, $M_{+-}^0$, $M_{-+}^0$, $M_{--}^0$.
  \item If $\min(p,q) \geq 2$, then the action of $\mathrm{SO}^+(p,q)$ on $\mathbb{P}_+(V^*)$ has two open orbits $M^+,M^-$ and one closed orbit $M^0$. 
 \end{enumerate}
\end{Proposition}

\proof
The first item is trivial. 

Let us prove the last statement. Let $e_1,\ldots,e_p,e_{p+1},\ldots,e_n$ be an orthonormal basis of $V^*$. For $v\in V^*$, $\R v$ has a natural orientation.

We claim that every (unoriented) line in $V^*$ is in the $\mathrm{SO}^+(p,q)$ orbit of $\R e_1, \R e_{p+1}$ or $\R(e_1+e_{p+1})$. Indeed, using elements in $\mathrm{SO}(p) \times \mathrm{SO}(q) \subset \mathrm{SO}^+(p,q)$ we get that every line is in the same orbit as a line given by $L=\R(\lambda_1 e_1+\lambda_{p+1}e_{p+1})$ with $\lambda_1,\lambda_{p+1} \geq 0$. If $\lambda_1=0$ we are done. Otherwise we may assume that $\lambda_1=1$. If $\lambda_{p+1}=1$, $L=\R(e_1+e_{p+1})$ and the claim follows. If $\lambda_{p+1} \neq 1$, we use the boost (hyperbolic rotation) which is the identity on $\Span(e_1, e_{p+1})^Q$ and maps $e_1 \mapsto \cosh \alpha e_1+\sinh \alpha e_{p+1},e_{p+1} \mapsto \sinh \alpha e_1+\cosh \alpha e_{p+1}$ with $\tanh \alpha=-\lambda_{p+1}$ if $0 \leq \lambda_{p+1}<1$ and with $\tanh \alpha=-\frac{1}{\lambda_{p+1}}$ if $\lambda_{p+1}>1$, to move $L$ to $\R e_1$ in the former and to $\R e_{p+1}$ in the latter case.

Using the diagonal matrix with entries $(-1,-1,1,\ldots,1)$ in $\mathrm{SO}(p)$ (here we use $p \geq 2$) and similarly for $\mathrm{SO}(q), q \geq 2$, we may reverse the orientation of each of these lines. Setting $M^+:={ \mathrm{SO}^+(p,q) \cdot} \R e_1, M^-:={\mathrm{SO}^+(p,q) \cdot} \R e_{p+1}, M^0:={ \mathrm{SO}^+(p,q) \cdot} \R(e_1+e_{p+1})$, the statement follows.    

For the second item, suppose that $q=1, p>1$. We argue as above. However, in this case the two orientations of the line $\R e_{p+1}$ define two different open $\mathrm{SO}^+(p,1)$-orbits $M^-_1,M^-_2$, as the set $\{x: Q(x)<0\}$ is disconnected. Similarly the two orientations of $\R(e_1+e_{p+1})$ define two different closed orbits. In the case $p=q=1$, the open orbits are $M_j^\pm:=\SO^+(1,1)(\pm\R e_j)$, $j=1,2$, and the closed orbits are $M_{\pm,\pm}^0:=\R(\pm e_1\pm e_2)$. 
\endproof
In all cases, $M^0$ will be referred to as the light cone.

\subsection{Meromorphic extension of $|\cos 2\theta|^\lambda$.}

For a pair of sets $Y\subset X$ and a function $f:X\to Z$ let us write $f_Y:=f\cdot {1}_{Y}:X\to Z$.
\begin{Proposition} \label{prop_meromorphic_cos2theta}
Let $\mathcal{O} \subset{ \mathbb P(\R^{p,q})}$ be an open orbit of $\SO^+(p,q)$. The function $E \mapsto |\cos 2\theta(E)|_\mathcal{O}^\lambda$ is integrable for $\mathrm{Re} \lambda>-1$ and extends meromorphically to a family of generalized functions with simple poles at $\{-1,-2,\ldots,\}$. The residues are supported on the closed orbit.
\end{Proposition}

\proof
We use a Euclidean trivialization to simplify the proof. We use the double cover $\pi:S^{n-1} \to  \mathbb P(\R^{p,q}), v \mapsto \R v$ and write $S^{n-1}$ as the spherical join of two spheres:
\begin{displaymath}
 \sigma:S^{p-1} \times S^{q-1} \times \left[0,\frac{\pi}{2}\right] \to S^{n-1}, (z_1,z_2,\theta) \mapsto (\cos(\theta)z_1,\sin(\theta)z_2). 
\end{displaymath}
The open orbits are $X_{1,0}=\pi \circ \sigma(S^{p-1} \times S^{q-1} \times [0,\frac{\pi}{4})), X_{0,1}=\pi \circ \sigma(S^{p-1} \times S^{q-1} \times (\frac{\pi}{4},\frac{\pi}{2}])), X_{0,0}=\pi \circ \sigma(S^{p-1} \times S^{q-1} \times \{\frac{\pi}{4}\})$.

Given a function $h$ on $\mathbb P(\R^{p,q})$, set
\begin{displaymath}
 \tilde h(\theta):=\sin(\theta)^{q-1}\cos(\theta)^{p-1} \int_{S^{p-1}} \int_{S^{q-1}}  (h \circ \pi \circ \sigma)(\cdot,\cdot,\theta) \ dz_2 dz_1.
\end{displaymath}

Then for $\mathrm{Re} \lambda>-1$,  
\begin{align*}
\int_{ X_{1,0}} |\cos(2\theta)|^\lambda h d\vol& =\int_{S^{p-1}} \int_{S^{q-1}} \int_0^{\frac{\pi}{ 4}} |\cos(2\theta)|^\lambda \sin(\theta)^{q-1}\cos(\theta)^{p-1} h \circ \pi \circ \sigma \ d\theta dz_2 dz_1\\
& =\int_0^{\frac{\pi}{ 4}} |\cos(2\theta)|^\lambda \tilde h(\theta) \ d\theta\\
& = { \frac12} \int_{ 0}^1 |x|^\lambda { \frac{1}{\sqrt{1-x^2}}} \tilde h\left(\frac12 \arccos x\right) dx
\end{align*}
and similarly 
\begin{displaymath}
 \int_{ X_{0,1}} |\cos(2\theta)|^\lambda h d\vol = { \frac12} \int_{ -1}^0 |x|^\lambda { \frac{1}{\sqrt{1-x^2}}} \tilde h\left(\frac12 \arccos x\right) dx.
\end{displaymath}

Recall \cite[Section 3.2]{hoermander_pde1} that $|x|^\lambda_+:=|x|^\lambda 1_{[0,\infty)}$ and  $|x|^\lambda_-:=|x|^\lambda 1_{(-\infty,0]}$ extend to meromorphic families in $\lambda\in\mathbb C \to C^{-\infty}(\R)$ with simple poles at the negative integers, with the residue at $\lambda=-k$ proportional to the $(k-1)$-st derivative of the delta function at the origin. This concludes the proof.
\endproof

\begin{Remark}
One can also define $|\cos2\theta|^\lambda_{\mathcal O}$ for any open $\mathrm{SO}^+(p,q)$-orbit $\mathcal O$ of the double cover $\mathbb P_+(\mathbb R^{p,q})\to \mathbb P(\mathbb R^{p,q})$ simply by restricting the corresponding pullback. Note also that the meromorphic families $|x|^\lambda=|x|_+^\lambda+|x|_-^\lambda$ and $\sign(x)|x|^\lambda=|x|_+^\lambda-|x|_-^\lambda$ are analytic at even, respectively odd, values $\lambda\in\mathbb Z$. It follows that the generalized functions $|\cos2\theta|^{2m}$ and $\sign(\cos2\theta)|\cos2\theta|^{2m+1}$,  which have full support, are well-defined for all $m\in\mathbb Z$.
\end{Remark}

\section{Dimension of the space of generalized invariant valuations}
\label{sec_dimension_generalized}

This section is devoted to the proof of Theorem \ref{thm_dimensions}. For the remainder of the paper, we will only consider $n\geq 3$.  

In the following, we will work with generalized forms. Let us first describe how generalized forms are related to currents in the sense of geometric measure theory. 

Fix $k,l$ and $\omega \in \Omega^{k,l}(V \times \mathbb{P}_+(V^*))^{tr}$ a translation-invariant form of bidegree $(k,l)$. Given $\phi \in \Omega^{n-k,n-l-1}(V \times \mathbb{P}_+(V^*))^{tr}$, we have $\omega \wedge \phi \in \Omega^{n,n-1}(V \times \mathbb{P}_+(V^*))^{tr}$ and therefore $\pi_*(\omega \wedge \phi) \in \Omega^n(V)^{tr} \cong \Dens(V)$. We thus get a continuous linear functional on the space $\Omega^{n-k,n-l-1}(V \times \mathbb{P}_+(V^*))^{tr}$ with values in $\Dens(V)$. More generally, a translation-invariant generalized form of bidegree $(k,l)$ is by definition a continuous linear functional on $\Omega^{n-k,n-l-1}(V \times \mathbb{P}_+(V^*))^{tr}$ with values in $\Dens(V)$. We will denote this space by $\Omega_{-\infty}^{k,l}(V \times \mathbb{P}_+(V^*))^{tr}$. 
 
By \cite{alesker_bernig}, a generalized translation-invariant valuation on $V$ of degree $1 \leq k \leq n-1$ is uniquely described by a closed, vertical and translation-invariant generalized form on $V \times \mathbb{P}_+(V^*)$ of bidegree $(k,n-k)$.

Given an element $\xi \in \mathbb{P}_+(V^*)$, we denote by $\xi_\perp \subset V$ the annihilator of $\xi$, which is a hyperplane. For a group $G\subset \mathrm{GL}(V)$, $\overline G$ is the group of transformations of $V$ generated by $G$ and parallel translations on $V$.

There is a canonical identification 
 \begin{equation}\label{eq:canonic_identification}
  \Omega_{-\infty}^{k,l}(V \times \mathbb{P}_+(V^*))^{\mathrm{tr}} \cong \Gamma^{-\infty}(\mathbb{P}_+(V^*), D^{k,l}),
 \end{equation}
where $D^{k,l}$ is the vector bundle over $\mathbb{P}_+(V^*)$ with fiber 
\begin{displaymath}
 D^{k,l}|_\xi:=\largewedge^k(V^*) \otimes \largewedge^l(\xi_\perp) \otimes \xi^l. 
\end{displaymath}

For the rest of the section, let $V$ be an $n$-dimensional linear space equipped with a non-degenerate quadratic form $Q$ of signature $(p,q)$. Write $G_0:=\mathrm{SO}^+(Q)=\mathrm{SO}^+(p,q) \subset \mathrm{GL}(V)$. There are $G_0$-equivariant identifications $V \cong V^*$ and $\xi_\perp \cong \xi^Q$.  Thus $D^{k,l}|_\xi\cong \largewedge^k(V)\otimes\largewedge^{l}(\xi^Q)\otimes\xi^{l}$ over $\xi\in\mathbb P_+(V)$.

We will be making frequent use of the following fact.
\begin{Lemma} \label{lemma_dimension_count_linear}
  For $0\leq a,b \leq n$,
 \begin{displaymath}
  (\largewedge^a V \otimes \largewedge^b V)^{\SO^+(Q)} \cong S_{a,b}(V) \oplus D_{a,b}(V),
 \end{displaymath}
with $\dim S_{a,b}(V)=\delta_a^b$ (symplectic part) and $\dim D_{a,b}(V)=\delta_{a+b}^n$ (determinantal part). The spaces $S_{a,a}(V)$ are invariant under the full orthogonal group $\mathrm{O}(Q)$, while $D_{a,n-a}(V)$ equals the sign representation of $\mathrm{O}(Q)$. 
\end{Lemma}

\proof
Let $\mathfrak g:=\mathfrak{so}(Q) \cong \mathfrak{so}(p,q)$ be the Lie algebra of $G_0$ and let $\mathfrak g_\C:=\mathfrak g \otimes_\R \C \cong \mathfrak{so}(n,\C)$ be its complexification. The complexification $V_\C:=V \otimes_\R \C$ is endowed with the $\mathfrak g_\C$-invariant quadratic form 
\begin{displaymath}
 Q_\C(x+iy)=Q(x)-Q(y)+2i Q(x,y).
\end{displaymath}

We have the natural isomorphism 
\begin{displaymath}
 (\largewedge^aV \otimes \largewedge^bV)^{\mathfrak{g}} \otimes_\R \C \cong (\largewedge^aV_\C \otimes \largewedge^bV_\C)^{\mathfrak{g}_\C}.
\end{displaymath}

We may equivariantly identify $\largewedge^aV \otimes \largewedge^bV$ and $\largewedge^aV_\C \otimes \largewedge^bV_\C$ with their respective duals. 

By the first fundamental theorem for $\mathrm{SO}(n,\C)$ \cite{goodman_wallach09}, the algebra of invariant polynomial functions on $V_\C^{a+b}$ is generated by the functions $Q_\C(v_{j_1},v_{j_2})$ and $\det(v_{j_1},\ldots,v_{j_n})$.  

It holds that either $(n-a)+b\leq n$ or $(n-b)+a\leq n$.	
The space $V_{\mathbb C}$ inherits a top form from $V$ which is $\mathrm{SO}(Q_{\mathbb C})$-invariant. Therefore we may $\mathrm{SO}(Q_{\mathbb C})$-equivariantly identify $\largewedge^c V_\C$ with $\largewedge^{n-c}V_\C$. Taking $c$ to be either $a$ or $b$, we may thus assume $a+b\leq n$. Note that this identification interchanges the conditions $a=b$ and $a+b=n$, and at the same time interchanges the trivial and sign representations of $\mathrm{O}(Q)$.

Let $f$ be an element of the subspace $(\largewedge^aV_\C \otimes \largewedge^bV_\C)^{\mathfrak{g}_\C}$. Clearly, $f$ can only contain determinantal-factors if $a+b=n$. In this case, 
\begin{displaymath}
D_{a,b}(V):=\mathrm{span}\{\det(v_1,\ldots,v_a,v_{a+1},\ldots,v_n)\} 
\end{displaymath}
is a direct summand of $(\largewedge^aV_\C \otimes \largewedge^bV_\C)^{\mathfrak{g}_\C}$. 

The other summand $S_{a,b}(V)$ (which is the only summand if $a+b<n$) therefore consists of $\mathrm{O}(n,\C)$-invariant polynomials, and this space is one-dimensional if $a=b$ and trivial otherwise.  This concludes the proof.

\endproof

We may describe the space $S_{a,a}(V)$ more precisely as follows. Using the identification $V \cong V^*$ induced by the quadratic form, the canonical symplectic form on $V \oplus V^*$ gives rise to an element in $\largewedge^{1,1}(V \oplus V)^{\mathrm{O}(Q)} \cong (V \otimes V)^{\mathrm{O}(Q)}$, and its $a$-th exterior power generates $S_{a,a}(V)$.

\subsection{Open orbits}

Let us define certain natural differential forms $\alpha, \beta_j$, where $0 \leq j\leq 2n-1$. The form $\alpha \in \Omega^{1,0} \left(V\times(\mathbb P_+(V^*)\setminus M^0)\right)^{\overline{\mathrm O(Q)}}$ is the $\mathrm O(Q)$-invariant contact form, which can be defined as the restriction of the canonic $1$-form on $V \times V^*$ to $V \times \{x\in V^*: Q(x)^2=1\}$. We subsequently identify $\{x\in V^*: Q(x)^2=1\}=\mathbb P_+(V^*) \setminus M^0$. It corresponds to the section $\alpha\in \Gamma(\mathbb P_+(V^*)\setminus M^0, D^{1,0})$ given by $\alpha(\xi)=\xi_Q$, where $\xi_Q\in \xi$ is the unique positively oriented vector with $|Q(\xi_Q)|=1$.
Then for $\epsilon=0,1$ and $0 \leq k \leq n-1$ define 
\begin{equation}\label{eq:beta_def}
\beta_{2k+\epsilon}:=(d\alpha)^k \wedge \alpha^\epsilon \in \Omega^{k+\epsilon,k} \left(V \times \left(\mathbb P_+(V^*) \setminus M^0\right)\right)^{\overline{\mathrm O(Q)}}.  
\end{equation}
Note that $d\beta_j=0$ if $j$ is even and $d\beta_j=\beta_{j+1} \neq 0$ if $j<2n-1$ is odd.

A generalized $\mathrm{SO}^+(Q)$-invariant and translation invariant differential form restricted to an open orbit is necessarily smooth. We will prove the following:
\begin{Proposition} \label{prop_forms_openorbit}
Let $1 \leq k \leq n-1$ and let $\mathcal O \subset \mathbb{P}_+(V^*)$ be an open $\mathrm{SO}^+(Q)$-orbit. There exists a unique (up to scale) translation-invariant,  $\mathrm{SO}^+(Q)$-invariant, vertical and closed form on $V \times \mathcal O$ of bidegree $(k,n-k)$, i.e. 
\begin{displaymath}
 \dim \Omega^{k,n-k}_V(V \times \mathcal O)^{ \overline{G_0}} \cap \ker(d)=1.
\end{displaymath}
Moreover, if $g\in\OO(Q)$ stabilizes $\mathcal O$ then $g^*$ acts as the scalar $\det g$ on this space.
If $k\neq \frac{n+1}{2}$, then 
\begin{displaymath}
 \dim \Omega^{k,n-k}_V(V \times \mathcal O)^{\overline{G_0}}=1,
\end{displaymath}
i.e. every invariant vertical form of bidegree $(k,n-k)$ is automatically closed. 
\end{Proposition}

\proof
Since $G_0=\mathrm{SO}^+(Q)$ acts transitively on $\mathcal O$, the group $\overline{G_0}$ acts transitively on $V \times \mathcal O$. We therefore have a natural isomorphism
\begin{displaymath}
 \Omega^{k,n-k}(V \times \mathcal O)^{\overline{G_0}} \cong (D^{k,n-k}|_\xi)^H\cong \left(\largewedge^k V^* \otimes \largewedge^{n-k} \xi^Q\right)^H,
\end{displaymath}
where $\xi \in \mathcal O$ is an arbitrary element, and $H$ is the stabilizer in $G_0$ of $\xi$. From $V^* \cong \xi \oplus \xi^Q$ we infer the more precise splitting 
\begin{equation} \label{eq_dec_forms}
 \Omega^{k,n-k}(V \times \mathcal O)^{\overline{G_0}} \cong \left(\largewedge^k \xi^Q \otimes \largewedge^{n-k} \xi^Q\right)^H \oplus \left(\largewedge^{k-1} \xi^Q \otimes \largewedge^{n-k} \xi^Q\right)^H,
\end{equation}
where the second summand corresponds to vertical forms. Hence 
\begin{displaymath}
 \Omega_V^{k,n-k}(V \times \mathcal O)^{\overline{G_0}} \cong \left(\largewedge^{k-1} \xi^Q \oplus \largewedge^{n-k} \xi^Q\right)^H.
\end{displaymath}

By Lemma \ref{lemma_dimension_count_linear}, applied to $W:=\xi^Q$ of dimension $(n-1)$, we obtain that 
\begin{displaymath}
 \dim \Omega_V^{k,n-k}(V \times \mathcal O)^{\overline{G_0}} = \begin{cases} 1 & k \neq \frac{n+1}{2}\\ 2 & k=\frac{n+1}{2}.\end{cases}
\end{displaymath}

Let $\omega_k \in \Omega_V^{k,n-k}(V \times \mathcal O)^{\overline{G_0}}$ be a representative corresponding to the determinantal part. If $n$ is odd and $k=\frac{n+1}{2}$, $\beta_{n} \in \Omega_V^{k,n-k}(V \times \mathcal O)^{\overline{G_0}}$ from eq. \eqref{eq:beta_def} is a representative of the symplectic part. We claim that $\omega_k$ is closed. To prove the claim, we first argue as for \eqref{eq_dec_forms} that 
\begin{displaymath}
 \dim \Omega^{k,n-k+1}(V \times \mathcal O)^{\overline{G_0}} = \begin{cases}  0 & k \neq \frac{n+1}{2}, \frac{n+2}{2}\\ 1 & k=\frac{n+1}{2}, \frac{n+2}{2}.\end{cases}
\end{displaymath}
In the cases $k=\frac{n+1}{2}, \frac{n+2}{2}$, the space $\Omega^{k,n-k+1}(V \times \mathcal O)^{\overline{G_0}}$ is spanned by $\beta_{n+1}$. 

It follows that $\omega_k$ must be closed if $k \neq \frac{n+1}{2}, \frac{n+2}{2}$. 
 
If $k \in \left\{\frac{n+1}{2}, \frac{n+2}{2}\right\}$, then $d\omega_k$ is a multiple of $\beta_{n+1}$, which is invariant under the subgroup of $\OO(Q)$ stabilizing $\mathcal O$. On the other hand, one may choose $g\in \OO(Q)$ stabilizing $\mathcal O$ with $\det g=-1$. Then $g^*\omega_k=-\omega_k$ by Lemma \ref{lemma_dimension_count_linear}, implying $g^*d\omega_k=-d\omega_k$. It follows that $d\omega_k=0$.
\endproof

\subsection{Global extensions}
\label{subsec_global_extension}

Fix an orientation on $V$. Then the natural volume form associated with $Q$ is preserved by $G_0$ and therefore there are $G_0$-isomorphisms $\largewedge^kV^* \cong \largewedge^nV^* \otimes \largewedge^{n-k}V \cong \largewedge^{n-k}V$. It follows that there is a $G_0$-isomorphism of vector bundles $D^{k,n-k} \cong D_0^{k}$, where 
\begin{displaymath}
D_0^{k}|_\xi = \largewedge^{n-k}(V) \otimes \largewedge^{n-k}(\xi^Q) \otimes \xi^{n-k}. 
\end{displaymath}

According to  Lemma \ref{lemma_dimension_count_linear} and Proposition \ref{prop_forms_openorbit}, there exists a $G_0$-invariant section $E_{(p,q)}^k$ of $D_0^k$ over $\mathbb P_+(V^*) \setminus M^0$ corresponding to a closed form $\omega_k$. 

Since the isomorphism $D^{k,n-k} \cong D_0^{k}$ interchanges the symplectic and the determinantal invariant subspaces, we may write down this section explicitly as
\begin{equation}\label{eq:def_of_E}
E_{(p,q)}^k(\xi) =  \sum_{I \subset \{1,\ldots,n-1\},\# I=n-k} \epsilon_I v_I \otimes v_I \otimes \xi_Q^{n-k}.
\end{equation}
Here $v_1,\ldots,v_{n-1}$ is a $Q$-orthonormal basis of $\xi^Q$, $v_I=v_{i_1} \wedge \ldots \wedge v_{i_{n-k}}$, $\xi_Q\in\xi$ the unique positively oriented vector with $|Q(\xi_Q)|=1$, and $\epsilon_I=Q(v_I,v_I)$. We will now rewrite this section in Euclidean terms using a compatible Euclidean form $P$ in the sense of Definition \ref{def_compatible_forms}. We will subsequently use this form to extend across the light cone certain linear combinations of the restrictions $E^k_{(p,q)}|_{\mathcal O}$ of $E_{(p,q)}^k$ to the different open orbits $\mathcal O$.

\begin{Lemma}
Let $S: V\to V$ be the map satisfying $Q(u,v)=P(Su,v)$. Let $\mathfrak U \subset \mathbb P_+(V^*)$ be the open dense  set of oriented lines $\xi$ s.t. $\xi \neq \pm S^*\xi$. Then the light cone $\{Q(\xi)=0\}$ is a subset of $\mathfrak U$.
\end{Lemma}

\proof
For $v \in V$ non-zero we have $Q(Sv,v)=P(v,v)>0$. If $S^*\xi=\pm \xi$, then $\xi_\perp$ is invariant under $S$ and hence $Q|_{\xi_\perp}$ is non-degenerate, as $S|_{\xi_\perp}$ is invertible. Therefore, $\xi$ is not in the light cone. 
\endproof

In the following, let $\xi \in \mathfrak U$. Set 
\begin{displaymath}
W_\xi:=\xi_\perp \cap S\xi_\perp=(\xi_\perp\cap V_p) \oplus (\xi_\perp \cap V_q), 
\end{displaymath}
which is a non-degenerate $(n-2)$-dimensional subspace of $V$. 

By Lemma \ref{Lem:TwoFormsOrthogonalComplements}, $W_\xi^P=W_\xi^Q$, and we define the line $\eta:=\xi_\perp\cap W_\xi^P$. Thus we have a decomposition
\begin{displaymath}
 \xi_\perp=(\xi_\perp \cap V_p) \oplus (\xi_\perp \cap V_q) \oplus \eta,
\end{displaymath}
which is both $P$- and $Q$- orthogonal. 

Let $\cos2\theta: \mathbb P_+(V^*)\to [-1,1]$ be defined by the identification of $\mathbb P(V^*)$ with $\Gr_{n-1}(V)$ through $\xi \mapsto \xi_\perp$. 
 We will write $\xi_P, \xi_Q\in\xi$ for the unique positively oriented covectors with $|P(\xi_P)|=|Q(\xi_Q)|=1$, so that $\xi_Q=|\cos2\theta|^{-1/2}\xi_P$. Define similarly $\eta_P, \eta_Q\in\eta$ with $\eta_Q=|\cos2\theta|^{-1/2}\eta_P$. Let us choose $Q$-orthonormal bases $\{v_1,\ldots,v_{p-1}\}$ of $\xi_\perp \cap V_p$, $\{v_p,\ldots,v_{n-2}\}$ of $\xi_\perp \cap V_q$ and $v_{n-1}:=\eta_Q \in \eta$. 

For any natural $l$ define $Q^l_\xi \in \largewedge^l(W_\xi) \otimes \largewedge^l(W_\xi)$ for $\xi \in \mathfrak U$ by taking it to be dual to the $Q$-product induced on $\largewedge^l(W_\xi)$. In terms of the given orthonormal basis, this may be written as 
\begin{displaymath}
 Q^l_\xi=\sum_{I \subset \{1,\ldots,n-2\},\# I=l} \epsilon_I v_I \otimes v_I,
\end{displaymath}
where $\epsilon_I=(-1)^{ \# (I \cap \{p,\ldots,n-2\})}$.

Over $\mathfrak U \setminus M^0$, 
setting as before $\epsilon_I=(-1)^{ \# (I \cap \{p,\ldots,n-2\})}$ if $n-1 \notin I$ and $\epsilon_I=(-1)^{ \# (I \cap \{p,\ldots,n-2\})} \cdot \mathrm{sign}(\cos 2 \theta)$ if $n-1 \in I$, we can split $E_{(p,q)}^k$ as follows. 
\begin{align*}
E_{(p,q)}^k(\xi) & = \sum_{\substack{\# I=n-k\\ n-1 \not\in I}} \epsilon_I v_I \otimes v_I \otimes \xi_Q^{n-k}+\sum_{\substack{\# I=n-k-1\\ n-1 \not\in I}} \epsilon_{I \cup \{n-1\}} (v_I \wedge v_{n-1}) \otimes (v_I \wedge v_{n-1}) \otimes \xi_Q^{n-k}\\
 & = \sum_{\substack{I \subset \{1,\ldots,n-2\}\\ \# I=n-k}} \epsilon_I v_I \otimes v_I \otimes \xi_Q^{n-k}+  \mathrm{sign} (\cos 2 \theta) \sum_{\substack{I \subset \{1,\ldots,n-2\}\\ \# I=n-k-1}} \epsilon_I v_I \otimes v_I \otimes \eta_Q^2\otimes \xi_Q^{n-k} \\
& = |\cos 2 \theta|^{-\frac{n-k}{2}} Q^{n-k}_\xi \otimes \xi_P^{n-k} +\mathrm{sign}(\cos 2 \theta) |\cos 2 \theta|^{-\frac{n-k}{2}-1} Q^{n-k-1}_\xi \otimes \eta_P^2\otimes \xi_P^{n-k} .
\end{align*}

Setting for $\xi \in \mathfrak U \setminus M^0$,
\begin{equation}\label{eq:def_of_AB}
A(\xi):=Q^{n-k}_\xi \otimes \xi_P^{n-k}, B(\xi):=Q^{n-k-1}_\xi \otimes \eta_P^2 \otimes \xi_P^{n-k},
\end{equation}
we may rewrite this as 
\begin{equation}\label{eq:Ekpq}
 E_{(p,q)}^k(\xi)=|\cos 2 \theta|^{-\frac{n-k}{2}} A(\xi)+\mathrm{sign}(\cos 2\theta) |\cos 2\theta|^{-\frac{n-k}{2}-1} B(\xi).
\end{equation}
 
Observe that $\xi_P$ is in fact well-defined for any $\xi\in \mathbb P_+(V^*)$, while $Q_\xi^{n-k}, Q_\xi^{n-k-1}, \eta_P$ are well-defined and smooth for $\xi\in\mathfrak U$, so that $A(\xi)$, $B(\xi)$ are well-defined and smooth over $\xi\in\mathfrak U$. 

Fix a non-negative cut-off function $\rho\in C^\infty[-1,1]$ which is identically 1 near 0 and identically 0 near $\pm1$. Write for $\xi\notin M^0$
\begin{equation}\label{eq:E_cutoff_decomposition}
 E_{(p,q)}^k(\xi)=\rho(\cos2\theta(\xi))E_{(p,q)}^k(\xi)+(1-\rho(\cos2\theta(\xi)))E_{(p,q)}^k(\xi).
\end{equation}
The function $1-\rho(\cos2\theta(\xi))$ vanishes in a neighborhood of $M^0$.

Hence the second summand of equation (\ref{eq:E_cutoff_decomposition}) is well-defined and smooth for all $\xi\in\mathbb P_+(V^*)$. 

We will now extend the first summand across the light cone. In order to do so, we replace $-\frac{n-k}{2}$ by a complex parameter $\lambda$ and use meromorphic extension (depending on the situation this amounts to analytic continuation or to taking residues) to define certain generalized sections $E^k_{(p,q)},F^k_{(p,q)}$ of $D_0^k$. Since $\rho(\cos2\theta(\xi))$ is compactly supported in $\mathfrak U$, we may use eq. \eqref{eq:Ekpq}. 

For every open orbit $\mathcal O$, define a meromorphic in $\lambda \in \mathbb C$ family of global sections 
\begin{align} \label{eq:def_of_E2}
E^k_{(p,q), \mathcal O, \lambda} & = \rho(\cos2\theta(\xi)) \left(|\cos2\theta|^{\lambda}A(\xi)+\sign(\cos 2\theta)|\cos2\theta|^{\lambda-1}B(\xi)\right)\\
& \quad   + \left(1-\rho(\cos2\theta(\xi))\right)|\cos2\theta|^{\lambda+\frac{n-k}{2}}E_{(p,q)}^k(\xi)
\end{align}
 supported on the closure of $\mathcal O$.

Note that when $\Re\lambda>1$, $E^k_{(p,q), \mathcal O, \lambda}$ is a continuous global section, vanishing on $M^0$ and coinciding with a multiple of $|\cos2\theta|^{\lambda+\frac{n-k}{2}} E^k_{(p,q)}$ on every open orbit. It follows that for $g\in G_0$,
\begin{equation} \label{eq_equivariance_a_b}
g^*E^k_{(p,q), \mathcal O, \lambda}=\psi_g(\xi)^{\lambda+\frac{n-k}{2}} E^k_{(p,q), \mathcal O, \lambda}.
\end{equation}

As both sides are meromorphic, we arrive at the following conclusion:
\begin{Corollary}\label{cor:invariant_sections}
Equation (\ref{eq_equivariance_a_b}) holds for arbitrary values of $\lambda$ in the domain of analyticity, as well as for the residues at the simple poles, for any linear combination of $E^k_{(p,q), \mathcal O, \lambda}$ over the open orbits.
\end{Corollary}

According to Proposition \ref{prop_orbits_lines}, the orbit space structure is different in the cases $\min(p,q)=0,1, >1$, and will require slightly different treatments. The case $\min(p,q)=0$ is of course the case of the classical Hadwiger theorem, so we are left with two new cases. 

When $\min(p,q)=1$, by Proposition \ref{prop_orbits_lines} the open $G_0$-orbits are $M_1^-,M_2^-,M^+$, and the closed orbits $M^0_1$, $M^0_2$.

When $\min(p,q)\geq 2$ there are two open $G_0$-orbits $M^+,M^-$, and one closed orbit $M^0$. In this case we will sometimes write $M_1^-=M^-$. We will denote by $N_c\in\{1,2\}$ the number of closed orbits, 

If $n-k$ is odd, the generalized functions $|\cos 2\theta|_\mathcal{O}^{-\frac{n-k}{2}}$ and $|\cos 2\theta|_\mathcal{O}^{-\frac{n-k}{2}-1}$ supported on the closure of each orbit $\mathcal O$ are well defined, yielding an invariant global section $E^k_{(p,q),\mathcal{O}}:=E^k_{(p,q),\mathcal{O}, -\frac{n-k}{2}}$ supported on $\overline{\mathcal O}$.

When $n-k$ is even, define the invariant global section $E_{(p,q),0}^{k}$ as the value at $\lambda=-\frac{n-k}{2}$ of
\begin{displaymath}
E^k_{(p,q), M^+, \lambda}+ (-1)^{\frac{n-k}{2}} \sum_{j=1}^{N_c}E^k_{(p,q), M^-_j, \lambda}
\end{displaymath}
The fact that this particular linear combination is in fact analytic at $\lambda$ is implied by Proposition \ref{prop_meromorphic_cos2theta} and the remark thereafter, using equation \eqref{eq:Ekpq}.
Still assuming $n-k$ is even, define the invariant global sections $F^k_{(p,q),j}, 1\leq j\leq N_c$, respectively supported on $M_j^0$, as the residue of $E^k_{(p,q), M_j^-, \lambda}$ at $\lambda=-\frac{n-k}{2}$.

It follows from Corollary \ref{cor:invariant_sections} that all the constructed sections are $G_0$-invariant.

We use the identification \eqref{eq:canonic_identification} to define the generalized $(k,n-k)$-forms $\omega_{k,0},\omega_{k,2}$ by
\begin{align} \label{eq:omegak0}
 \omega_{k,0} & :=\begin{cases} E^k_{(p,q),M^+}-\sum_{j=1}^{N_c}E^k_{(p,q),M^-_j} & n-k \equiv 1,3 \mod 4,\\
                E^k_{(p,q),0} & n-k \equiv 2 \mod 4,\\
                \sum_{j=1}^{N_c} F^k_{(p,q),j} & n-k \equiv 0 \mod 4;
               \end{cases}\\
\label{eq:omegak2}  \omega_{k,2} & :=\begin{cases} E^k_{(p,q),M^+}+\sum_{j=1}^{N_c}E^k_{(p,q),M^-_1} & n-k \equiv 1,3 \mod 4,\\
                \sum_{j=1}^{N_c} F^k_{(p,q),j} & n-k \equiv 2 \mod 4,\\
                E^k_{(p,q),0} & n-k \equiv 0 \mod 4.
               \end{cases}
\end{align}

If $\min(p,q)=1$ we also define $\omega_{k,1}$ by
\begin{align} \label{eq:omegak1}
 \omega_{k,1} & :=\begin{cases} E^k_{(p,q),M_1^-}-E^k_{(p,q),M_2^-} & n-k \equiv 1 \mod 2,\\
 	F^k_{(p,q),1}-F^k_{(p,q),2} & n-k \equiv 0 \mod 2;
 \end{cases}
\end{align}

We summarize the results above in the following proposition.

\begin{Proposition}\label{prop:open_orbits_forms_summary} 
For every $1\leq k\leq n-1$ there are linearly independent forms $\omega_{k,j}\in  \Omega_{-\infty,V}^{k,n-k}(V \times \mathbb P(V^*))^{\overline{G_0}}$, with $j=0,2$ if $\min(p,q)\geq 2$ and $j=0,1,2$ if $\min(p,q)= 1$. These forms are closed on the open orbits.
It holds that $g^*\omega_{k,j}=\det(g) \omega_{k,j}$ for $g\in \mathrm O(Q)$ and $j=0,2$. When $\min(p,q)=1$ and $g\in \mathrm O(Q)$, $g^*\omega_{k,1}=\det(g)\epsilon_0(g)\omega_{k,1}$, where $\epsilon_0(g)=\pm1$ according to whether $g$ stabilizes the connected components of the complement of the light cone. 
\end{Proposition}

We will now show that there are no other extensions of the vertical closed $\overline{G_0}$-invariant differential forms on the open orbits across the light cone, up to forms supported on the light cone.

\begin{Lemma}\label{lem:no_other_extensions}
Assume $n-k\equiv 0\mod 2$. Let $\omega\in \Omega_{-\infty}^{k,n-k}(V \times \mathbb P_+(V^*))^{\overline {G_0}}$ restrict to 
\begin{displaymath}
 \omega|_{V\times(\mathbb P_+(V^*)\setminus M^0)}=c^+E^k_{(p,q)}|_{M^+}+c^{-}E^k_{(p,q)}|_{M^-}
\end{displaymath}
if $\min(p,q)\geq 2$, and to
\begin{displaymath}
\omega|_{V\times(\mathbb P_+(V^*)\setminus M^0)}=c^+E^k_{(p,q)}|_{M^+}+c^{-}_1E^k_{(p,q)}|_{M_1^-}+c^{-}_2E^k_{(p,q)}|_{M_2^-}
\end{displaymath}
if $\min(p,q)=1$. Then $c^+=(-1)^{\frac{n-k}{2}}c^-$ in the first case, and $c^+=(-1)^{\frac{n-k}{2}}c_j^-$ for $j=1,2$ in the latter case.
\end{Lemma}

\proof We will only consider the case of $n-k\equiv 0 \mod 4$, as the other case can be treated similarly.
We will use the identification given by equation \eqref{eq:canonic_identification}.
Setting  $e:=\frac{n-k}{2} \in \mathbb N$, we may represent 
$\omega=f(\xi) A(\xi)+h(\xi)B(\xi)$ for some global extensions $f,h\in C^{-\infty}(\mathbb P_+(V^*)\setminus\mathfrak U)$ of the functions 
$c^+|\cos2\theta|^{-e}1_{M^+}+c^-|\cos2\theta|^{-e}1_{M^-}$ and respectively $c^+|\cos2\theta|^{-e-1}1_{M^+}-c^-|\cos2\theta|^{-e-1}1_{M^-}$ (with the obvious modification when $\min(p,q)=1$).

Use $Q$ to identify $V\cong V^*$. Fix a $Q$-orthogonal, $S$-invariant decomposition $V=\mathbb R^{1,1}\oplus \mathbb R^{p-1,q-1}$, and write $L:=\mathbb P_+(\mathbb R^{1,1})\subset \mathbb P_+(V)$. 

Take $g_\alpha \in \mathrm{SO}^+(Q)$ fixing $\mathbb R^{p-1,q-1}$ and acting by an $\alpha$-boost on $\mathbb R^{1,1}$. Recall that $A(\xi)=Q^{n-k}_\xi \otimes \xi_P^{n-k}, B(\xi):=Q^{n-k-1}_\xi \otimes \eta_P^2 \otimes \xi_P^{n-k}$ for $\xi \in \mathfrak U$. For $\xi \in L$ it holds that $W_\xi=\mathbb R^{p-1,q-1}$, and therefore also $Q_\xi^l$ is independent of $\xi$. Thus $A|_{L \cap \mathfrak U}$, $B|_{L \cap \mathfrak U}$ may be extended as smooth sections over $L$.
 	
We may write  $g_\alpha^*\omega(\xi)=g_\alpha^*f(\xi) \cdot g_\alpha^*A(\xi)+g_\alpha^*h(\xi)\cdot g_\alpha^*B(\xi)$, 
where by definition $g^*s(\xi)=g^{-1}s(g\xi)$ for a section $s$ of any $G_0$-vector bundle over $\mathbb P_+(V^*)$. 

Since the connected components of $M^0$ are in fact orbits of the maximal compact subgroup of $G_0$, it follows from the $G_0$-invariance of $\omega$ by Lemma \ref{lem:invariant_wavefront} that $\WF(\omega)\subset N^*(M^0)$.  We may therefore restrict $\omega$ (as a global section of $D^k_0$) to $L$. Then $f|_L,h|_L$ extend to generalized functions on $L$, smooth outside $M^0$.

Now for all $l$ and all $\xi \in L$, $g_\alpha^*Q_\xi^{l}|_L=Q_\xi^{l}|_L$ and $g_\alpha^* \xi_P|_L=\psi_{g_\alpha}(\xi)^{\frac12} \cdot \xi_P|_L$.
Since $g_\alpha^*\omega=\omega$, it follows that $g_\alpha^*f|_L=\psi_{g_\alpha}^{-e}f|_L$.

The following claim then completes the proof.

{\it Claim}. Let $f\in C^{-\infty}(\mathbb P_+(\mathbb R^{1,1}))$ satisfy $g_\alpha^*f=\psi_{g_\alpha}^{-e}f$ for all $g_\alpha\in \mathrm{SO}^+(1,1)$. Then $f$ is a linear combination of $|\cos2\theta|^{-e}$ and a certain $f_0$ supported on the light cone. 

To see this, let $\theta$ be the Euclidean angle on $S^1=\mathbb P_+(\mathbb R^{1,1})$, and restrict $f$ to $\theta \in \left(0,\frac{\pi}{2}\right)$. Consider only $\alpha>0$ so that $g_\alpha\left(0,\frac \pi 2\right) \subset \left(0,\frac \pi 2\right)$. Make the change of coordinates $x=\cos 2\theta$, $F(x) =f(\theta)$, $F\in C^{-\infty}(-1,1)$. 

Write $g_\alpha(\cos \theta e_1+\sin \theta e_2)=c(\cos \theta' e_1+\sin \theta' e_2)$, where $e_1,e_2$ is an orthonormal basis of $\R^{1,1}$, $c \in (0,\infty)$ and $x'=\cos 2\theta'$. 

Then
\begin{displaymath}
x=\cos^2\theta-\sin^2\theta=c^2(\cos^2\theta'-\sin^2\theta')=c^2 x'
\end{displaymath}
and 
\begin{displaymath}
c^2=(\cosh\alpha \cos\theta+\sinh\alpha\sin\theta)^2+(\sinh\alpha \cos\theta+\cosh\alpha\sin\theta)^2=\cosh2\alpha+\sinh2\alpha\sin2\theta.
\end{displaymath}
It follows that 
\begin{displaymath}
x'=g_\alpha(x):=\frac{x}{\cosh2\alpha+\sinh2\alpha\sqrt{1-x^2}}.
\end{displaymath}

Hence $F$ satisfies the functional equation $g_\alpha^*F= \psi_{g_\alpha}^{-e}F$
and therefore, by Proposition \ref{prop_jacobian}, 
\begin{displaymath}
\Psi_{g_\alpha}(x) =\frac{g_\alpha(x)}{x}=\frac{1}{\cosh2\alpha+\sinh2\alpha\sqrt{1-x^2}}. 
\end{displaymath}

Differentiating the functional equation at $\alpha=0$ we get $L_{\underline{X}}F=2e\sqrt{1-x^2}F$, where 
\begin{displaymath}
 \underline{X}=\left.\frac{d}{d\alpha}\right|_{\alpha=0}g_\alpha(x)=-2x\sqrt{1-x^2}{\frac{\partial}{\partial x}}.
\end{displaymath}

Thus $xF'=-eF$, implying that $F$ is $(-e)$-homogeneous, that is, $F(\lambda x)=\lambda^{-e} F(x)$ for $0<\lambda<1$, where the equation has to be understood in the appropriate sense for generalized functions.
But the space of $(-e)$-homogeneous functions is $2$-dimensional, spanned by $|x|^{-e}$ and $\delta_0^{(e-1)}$. Repeating the argument around the other points of the light cone in $\mathbb P_+(\mathbb R^{1,1})$ completes the proof of the claim and the lemma.
\endproof

\subsection{Generalized $n$-forms supported on the closed orbit}

Next, we study the space of closed vertical invariant generalized forms supported on the closed orbit $M^0$. 

\subsubsection{First steps} \label{subsubsec_first_steps}

Fix $\xi \in M^{0}$ with stabilizer $H:= \Stab(\xi) \subset G_0$, $\alpha$, $0 \leq k \leq n$, $0\leq l\leq n-1$. Denote by $H_+:=\Stab_+(\xi)=\{g\in \Stab(\xi): \det (g|_\xi)>0\}$. 
Define the vector bundle $F_{k,l}^{\alpha}$ over $M^{0}$ with fiber 
\begin{displaymath}
F_{k,l}^\alpha|_\xi=\largewedge^k V^* \otimes \largewedge^l \xi_{\perp} \otimes \xi^{l} \otimes \Sym^{\alpha}(N_{\xi}M^{0}) \otimes \Dens^*(N_\xi M^{0}). 
\end{displaymath}

Those are the auxiliary vector bundles that appear in Lemma \ref{lem:LocallyClosedOrbits} associated with the vector bundle $D^{k,l}$. We denote by $\Omega_{(-\infty,\alpha), M^0}$ the space of generalized differential forms supported on $M^0$ of differential order normal to $M^0$ not greater than $\alpha$, see Appendix \ref{appendix}.

Set
\begin{displaymath}
T^{k,l,\beta}:=\largewedge^kV \otimes \largewedge^l \xi^Q \otimes \xi^\beta.
\end{displaymath}

\begin{Lemma}\label{lem:value_of_alpha}
There is a $\Stab_+(\xi)$-equivariant isomorphism,
\begin{displaymath}
 F_{k,l}^\alpha|_\xi \cong T^{k,l,l-2\alpha-2}.
\end{displaymath}
\end{Lemma}

\proof
By Proposition \ref{prop_EquivariantFormOfNormalBundle}, $N_\xi M^0$ is $\Stab(\xi)$-equivariantly isomorphic to $(\xi^*)^2$. It follows that $\Dens^*(N_\xi M^{0})=\Dens^*((\xi^*)^2)=\Dens(\xi)^2=(\xi^*)^2$, where the last equality is $\Stab_+(\xi)$-equivariant.  Recall that $Q$ defines natural isomorphisms $V \cong V^*$, $\xi_\perp \cong \xi^Q$ and $V^*/\xi^Q\cong V/\xi_\perp \cong \xi^*$. We therefore get
\begin{displaymath}
 F_{k,l}^\alpha|_\xi \cong \largewedge^k V \otimes \largewedge^l \xi^Q \otimes \xi^l \otimes (\xi^*)^{2\alpha} \otimes (\xi^*)^2 \cong \largewedge^kV \otimes \largewedge^l \xi^Q \otimes \xi^{l-2\alpha-2}\cong T^{k,l,l-2\alpha-2}.
\end{displaymath}
\endproof
For any $k,l,\beta\geq 0$, the subspace 
\begin{displaymath}
 U:=\largewedge^k\xi^Q \otimes \largewedge^l \xi^Q \otimes \xi^\beta \subset T^{k,l,\beta}
\end{displaymath}
is $H_+$-invariant. 

By Lemma \ref{lem:long_quotient}, there is a $\Stab(\xi)$-equivariant isomorphism 
\begin{displaymath}
\largewedge^k V/\largewedge^k \xi^Q\cong\largewedge^{k-1} \xi^Q \otimes (V/\xi^Q) \cong \largewedge^{k-1} \xi^Q \otimes \xi^*,
\end{displaymath}
and so the quotient $W:=T^{k,l,\beta}/U$ is $\Stab(\xi)$-equivariantly isomorphic to 
\begin{displaymath}
 W \cong \largewedge^{k-1} \xi^Q \otimes \largewedge^l \xi^Q \otimes \xi^{\beta-1}.
\end{displaymath}

If there is a non-zero $G_0$-invariant section of $F_{k,l}^\alpha$, then it defines an element in $T^{k,l,\beta}$ which is non-zero and $H_+$-invariant. Then either its image in $W$ is non-zero (and $H_+$-invariant), or the element belongs to $U$ and is $H_+$-invariant. We are therefore led to study $H_+$-invariant elements in $U$ and $W$. 

This procedure will be iterated several times. The notation will be as following: $U$ stands for a subspace, $W$
for the corresponding quotient space. At each new level, the order of the letters is preserved. Thus $U_{uw}$ is a subspace of $U_{w}$ and the quotient space is $U_{w}/U_{uw}=W_{uw}$. If the value of $\beta$ isn't clear
from the context, we indicate $\beta$ explicitly, e.g. $U_{u,\beta}$. 

The kernel of the quotient map $\largewedge^{k-1} \xi^Q \to \largewedge^{k-1}(\xi^Q/\xi)$ is $H_+$-isomorphic to $\largewedge^{k-2}(\xi^Q/\xi) \otimes \xi$. We therefore can define  
\begin{align*}
U_w & :=(\largewedge^{k-2}(\xi^Q/\xi) \otimes \xi) \otimes \largewedge^l \xi^Q \otimes \xi^{\beta-1}=\largewedge^{k-2}(\xi^Q/\xi) \otimes \largewedge^l \xi^Q \otimes \xi^\beta \subset W, \\
W_w & := W/U_w \cong \largewedge^{k-1}(\xi^Q/\xi) \otimes \largewedge^l \xi^Q \otimes \xi^{\beta-1}.
\end{align*}

Setting 
\begin{displaymath}
 A_{k,l,\beta}:=\largewedge^k(\xi^Q/\xi) \otimes \largewedge^l (\xi^Q/\xi) \otimes \xi^\beta,
\end{displaymath}
we get the spaces
\begin{align*}
 U_{ww} & :=A_{k-1,l-1,\beta}  \subset W_w,\\
 W_{ww} & :=W_w/U_{ww} \cong A_{k-1,l,\beta-1}, \\
 U_{uw} & :=A_{k-2,l-1,\beta+1} \subset U_w,\\
 W_{uw} & := U_w/U_{uw} \cong A_{k-2,l,\beta}.
\end{align*}

Starting with $U$ instead of $W$, we define the following spaces.
\begin{align*}
 U_u & :=\largewedge^{k-1}(\xi^Q/\xi) \otimes \largewedge^l \xi^Q \otimes \xi^{\beta+1} \subset U,\\
 W_u & := U/U_u \cong \largewedge^k(\xi^Q/\xi) \otimes \largewedge^l \xi^Q \otimes \xi^\beta,\\
 U_{wu} & :=A_{k,l-1,\beta+1} \subset W_u,\\
 W_{wu} & := W_{u}/U_{wu} \cong A_{k,l,\beta},\\
 U_{uu} & :=A_{k-1,l-1,\beta+2} \subset U_u,\\
 W_{uu} & :=U_u/U_{uu} \cong A_{k-1,l,\beta+1}.
\end{align*}

The dimension of the space of $H_+$-invariant elements in $A_{k,l,\beta}$ can be computed with Lemma \ref{lemma_dimension_count_linear}. The space $\xi^Q/\xi$ is of dimension $n-2$ and inherits a quadratic form of signature $(p-1,q-1)$. Both $\Stab(\xi)$ and $\Stab_+(\xi)$ act as $\mathrm{SO}^+(p-1,q-1)$ on $\xi^Q/\xi$.

If $k+l=n-2=\dim \xi^Q/\xi$ and $\beta=0$, there exists a one-dimensional subspace in $A_{k,l,\beta}$ of invariants corresponding to the determinantal part. If $k=l$ and $\beta=0$, there exists another one-dimensional space corresponding to the symplectic part. 

Let us first study the case $l=n-k$. The conditions for invariant subspaces can be summarized in the following table
\begin{equation}\label{tbl1}
 \begin{array}{c| c | c}
  & \text{determinantal} & \text{symplectic}\\ \hline
  U_{ww} & \beta=0 & n=2k, \beta=0\\
  W_{ww} & - & n=2k-1, \beta=1\\
  U_{uw} & - & n=2k-1,\beta=-1\\
  W_{uw} & \beta=0 & n=2k-2,\beta=0\\
  U_{wu} & - & n=2k+1,\beta=-1\\
  W_{wu} & - & n=2k, \beta=0\\
  U_{uu} & \beta=-2 & n=2k, \beta=-2\\
  W_{uu} & - & n=2k-1,\beta=-1
 \end{array}
\end{equation}

\begin{Corollary}
 If $n-k$ is odd and $n \neq 2k-1,2k+1$, there exists no invariant section of $F_{k,n-k}^\alpha$. 
\end{Corollary}

\proof
By Lemma \ref{lem:value_of_alpha}, such an invariant section defines an invariant element in $T^{k,n-k, \beta}$ with $\beta=n-k-2\alpha-2$ and odd integer. However, given the assumption on $n$ and $k$, none of the admissible spaces in the table above contains an invariant element.
\endproof

\subsubsection{Vertical forms}
The next step is to study the space of sections lifting to vertical generalized forms. 
\begin{Lemma} \label{lemma_vertical_lifting}
Let $s \in F_{k,n-k}^\alpha|_\xi$ lift to a vertical generalized form $\omega \in \Omega^{k,n-k}_{(-\infty,\alpha),M^0}(V \times \mathbb P_+(V^*))^{tr}$. Set $\beta=n-k-2\alpha-2$. Then 
\begin{enumerate} 
\item  $s \in \xi \otimes \largewedge^{k-1}(V/\xi) \otimes \largewedge^{n-k}\xi^Q \otimes \xi^{\beta}$. 
\item If $s \notin U$, then $0 \neq \mathrm{Pr}_W(s) \in U_w$.
\item If $s \in U$, then $s \in U_u$. 
\end{enumerate}
\end{Lemma}

\proof
\begin{enumerate}
\item 
The space of smooth translation-invariant forms of bidegree $(k,n-k)$ is given by 
\begin{displaymath}
 \Omega^{k,n-k}(V \times \mathbb{P}_+(V^*))^{tr} \cong \largewedge^k V^* \otimes \Omega^{n-k}(\mathbb{P}_+(V^*)).
\end{displaymath}
An element $\omega \in \Omega^{k,n-k}(V \times \mathbb{P}_+(V^*))^{tr}$ is vertical if and only if for each $\xi \in \mathbb{P}_+(V^*)$ we have $\omega|_\xi \in \xi \otimes \largewedge^{k-1}(V/\xi) \otimes \largewedge^{n-k} T^*_\xi \mathbb{P}_+(V^*)$. Correspondingly, an element $\omega \in \Omega^{k,n-k}_{(-\infty,\alpha),M^0}(V \times \mathbb{P}_+(V^*))^{tr}$ is vertical only if for each $\xi \in M^0$ its image in $F_{k,n-k}^\alpha|_\xi$ belongs to $\xi \otimes \largewedge^{k-1}(V/\xi) \otimes \largewedge^{n-k}\xi^Q \otimes \xi^{\beta}$. 

\item The image of $\xi \otimes \largewedge^{k-1}(V/\xi) \subset \largewedge^k V$ under the quotient map $\largewedge^k V \to \largewedge^k V / \largewedge^k \xi^Q$, which is essentially the quotient map $T^{k,n-k,\beta}\to W$, equals $\xi\otimes \largewedge^{k-2} (\xi^Q/\xi)\otimes V/\xi^Q\cong \largewedge^{k-2} (\xi^Q/\xi)$. It follows by i) that the image of $s$ in $W=(\largewedge^kV/\largewedge^k\xi^Q) \otimes \largewedge^{n-k}\xi^Q \otimes \xi^{\beta}$ is contained in $\largewedge^{k-2}( \xi^Q/\xi) \otimes \largewedge^{n-k}\xi^Q \otimes \xi^{\beta}$, which is just $U_w$. Since $s \notin U$, $\mathrm{Pr}_W(s)\neq0$. 
\item
Immediate from i) and the definition of $U_u$.
\end{enumerate}
\endproof

\begin{Proposition} \label{prop_dim_vertical_forms}
Assume $1\leq k\leq \frac n 2$. The space of $G_0$-invariant sections of $F_{k,n-k}^\alpha$ lifting to a vertical generalized form is at most one-dimensional if $n-k$ is even, $n \neq 2k$ and $\alpha=\frac{n-k}{2}$; it is at most two-dimensional if $n=2k$ with $k$ even and $\alpha=\frac{k}{2}$. These invariants belong to $U_{uu}$. In all other cases there are no such $G_0$-invariants. 
\end{Proposition}

\proof
Set $\beta=n-k-2\alpha-2$. We claim that there is no non-zero $H_+$-invariant element $s \in F^\alpha_{k,n-k}$ lifting to a vertical generalized form if $\beta=0$. 

Let us first show how the claim implies the statement of the proposition. 

Let $s \in F^\alpha_{k,n-k}$ be a non-zero invariant section lifting to a vertical generalized form. If $s \notin U$, then $0 \neq \mathrm{Pr}_W(s) \in U_w$ by Lemma \ref{lemma_vertical_lifting}. In this case, $U_{uw}$ or $W_{uw}$ must contain a non-zero invariant. But this is not the case by the assumption $k\leq n/2$ or respectively by the claim. 

Hence $s \in U$, and Lemma \ref{lemma_vertical_lifting} implies that in fact $s \in U_u$. It follows that $U_{uu}$ or $W_{uu}$ contains a non-zero invariant. The latter space cannot contain an invariant element by the assumption $k\leq n/2$, hence the former space does. 

Since $U_{uu}$ is a subspace of $T^{k,l,\beta}$, the $H_+$-invariant elements in $U_{uu}$ can be considered as $H_+$-invariant elements in $T^{k,l,\beta}$. Now $U_{uu}$ contains an invariant line corresponding to the determinantal part if $\alpha=\frac{n-k}{2}$; and another invariant line if $n=2k$ with $k$ even, $\alpha=\frac{k}{2}$ corresponding to the symplectic part. 
\newline\newline
It remains to prove the claim. 

Suppose that $0 \neq s \in F^\alpha_{k,n-k}$ lifts to a vertical generalized form, and $\beta=0$.  

If $s \in U_{u,\beta=0}$, then $U_{uu,\beta=0}$ or $W_{uu,\beta=0}$ would contain an invariant element, which is not the case. By Lemma \ref{lemma_vertical_lifting} again, $s \notin U$ and $0 \neq \mathrm{Pr}_W(s) \in U_{w,\beta=0}=\largewedge^{k-2}(\xi^Q/\xi)\otimes \largewedge^{n-k}\xi^Q$. 

Fix a subspace $Y \subset \xi^Q$ complementary to $\xi$ and let $H_Y \subset H_+$ be the subgroup stabilizing $Y$. Under the action of $H_Y$ we can decompose
\begin{displaymath}
 U_{w,\beta=0}=\largewedge^{k-2}(\xi^Q/\xi) \otimes \largewedge^{n-k}Y \oplus \largewedge^{k-2}(\xi^Q/\xi) \otimes \largewedge^{n-k-1}Y \otimes \xi.
\end{displaymath}
The second summand does not contain any non-zero $H_Y$-invariant element, hence $\mathrm{Pr}_W(s)$ lies in the first summand. However, under the action of the larger group $H_+$, $\mathrm{Pr}_W(s)$ does not remain in the first summand, which is a contradiction. 

\endproof

\begin{Proposition} \label{prop_vertical_n_forms}
Assume $1\leq k\leq \frac n 2$. The space of vertical generalized translation and $G_0$-invariant $(k,n-k)$-forms supported on the closed orbit $M^0$ has the following dimension:
 \begin{displaymath}
  \dim \Omega_{-\infty,V,M^0}^{k,n-k}(V \times \mathbb P(V^*))^{\overline{G_0}}=\begin{cases} 0 & k \not\equiv n \mod 2,\\ 1 & k \equiv n \mod 2, \min(p,q)\geq 2, \\ 2 & k \equiv n \mod 2,\min(p,q)=1.\end{cases} 
 \end{displaymath}
\end{Proposition}

\proof

Let $N_c$ denote the number of connected components of $M^0$. Thus $N_c=2$ if $\min(p,q)=1$, and $N_c=1$ otherwise.

If $n-k$ is odd, then by Proposition \ref{prop_dim_vertical_forms} there is no non-zero invariant generalized vertical translation-invariant $(k,n-k)$-form supported on $M^0$.

Let us assume that $n-k$ is even and $n \neq 2k$. Again by Proposition \ref{prop_dim_vertical_forms}, applied to each connected component of $M^0$, we have $\dim \Omega_{-\infty,V,M^0}^{k,n-k}(V \times \mathbb P(V^*))^{\overline{G_0}} \leq N_c$. Since we already constructed $N_c$ linearly independent elements in this space (see Subsection \ref{subsec_global_extension}), we must have equality. More precisely, if $n-k \equiv 0 \mod 4$, then $\omega_{k,2} \in \Omega_{-\infty,V,M^0}^{k,n-k}(V \times \mathbb P(V^*))^{\overline{G_0}}$. If $n-k \equiv 2 \mod 4$, then $\omega_{k,0} \in \Omega_{-\infty,V,M^0}^{k,n-k}(V \times \mathbb P(V^*))^{\overline{G_0}}$. If $\min(p,q)=1$, then in both cases also $\omega_{k,1}\in \Omega_{-\infty,V,M^0}^{k,n-k}(V \times \mathbb P(V^*))^{\overline{G_0}}$.

It remains to study the case $n=2k$ with $k$ even. In this case, Proposition \ref{prop_dim_vertical_forms} implies that $\dim \Omega_{-\infty,V,M^0}^{k,k}(V \times \mathbb P(V^*))^{\overline{G_0}} \leq 2N_c$. 

Identify 
\begin{displaymath}
 \Omega^{k,k}_{-\infty,M^0}(V \times \mathbb P_+(V^*))^{\overline{G_0}}=\Gamma^{-\infty}_{M^0}(\mathbb P_+(V^*), D^{k,k})^{\overline{G_0}}.
\end{displaymath}
Under this identification, vertical forms correspond to sections of $\xi\otimes \largewedge^{k-1}(V^*/\xi)\otimes \largewedge^{k}\xi_\perp \otimes \xi^k$. 

We already constructed $N_c$ linearly independent elements $\gamma_\pm$ in this space:  If $n-k \equiv j \mod 4$, take $\gamma_+=\omega_{k,j}$, where $j=0,2$. If $\min(p,q)=1$ we also have $\gamma_-=\omega_{k,1}$.

By the proof of Proposition \ref{prop_dim_vertical_forms}, as $0 \neq \gamma_\pm \in \Omega_{-\infty,V,M^0}^{k,k}(V \times \mathbb P(V^*))^{\overline{G_0}}$, its image in $\Gamma^{-\infty}_{M^0}(\mathbb P_+(V^*),D^{k,k})$ has  differential order normal to $M^0$ equal to $\alpha:=\frac{k}{2}$. 

Recall the definition of $A(\xi), B(\xi)$ for $\xi \in \mathfrak U$ in equation \eqref{eq:def_of_AB}. We denote the residue at $\lambda=-\frac{k}{2}$ of 	
\begin{displaymath}
\tilde E_{(p,q)}^k(\xi)=|\cos 2 \theta|^{\lambda} A(\xi)+\mathrm{sign}(\cos 2\theta) |\cos 2\theta|^{\lambda-1} B(\xi)
\end{displaymath}
by $\zeta_+$ if $\min(p,q)\geq 2$. If $\min(p,q)=1$, this residue is supported on two connected components, and we denote the corresponding generalized forms by $\zeta_\pm$.
Note that $\gamma_\pm$ is given by the same formula, except that $\Omega^{k,k}(V \times \mathbb P(V^*))^{tr}$ is identified with sections of $D^k_0$. Since $Q_\xi^k$ belongs to $\largewedge^kW_\xi \otimes \largewedge^kW_\xi \subset \largewedge^k \xi_\perp \otimes \largewedge^k \xi_\perp$ and is non-zero (since $k \leq \dim W_\xi=n-2$), $\zeta_\pm$ is not vertical.

Observe also that $\zeta_\pm$ has differential order normal to $M^0$ equal to $\alpha$. The space $S_\zeta$ spanned by $\zeta_\pm$ is $N_c$-dimensional and contains only non-vertical forms (except for the trivial combination). Let 
\begin{displaymath}
 S_V:=\Omega_{-\infty,V,M^0}^{k,k}(V \times \mathbb P(V^*))^{\overline{G_0}} \subset \Omega^{k,k}_{-\infty,M^0}(V \times \mathbb P_+(V^*))^{\overline{G_0}}
\end{displaymath}
be the space spanned by vertical invariant forms, so that $S_V\cap S_\zeta=\{0\}$. It follows from Table \ref{tbl1} with $\beta=n-k-2\alpha-2=-2$ that $\dim \Gamma^\infty(M^0,F^\alpha)^{G_0}=2N_c$.  
Since the map
\begin{displaymath}
 \Omega^{k,k}_{(-\infty,\alpha),M^0}(V \times \mathbb P_+(V^*))^{\overline{G_0}}/\Omega^{k,k}_{(-\infty,\alpha-1),M^0}(V \times \mathbb P_+(V^*))^{\overline{G_0}} \to \Gamma^\infty(M^0,F^\alpha)^{G_0}
\end{displaymath}
is injective by  \cite[Proposition 4.9]{alesker_faifman}, we have $2N_c \geq \dim (S_\zeta\oplus S_V)=\dim S_\zeta+\dim S_V=N_c+\dim S_V$ and hence $\dim S_V \leq N_c$ as claimed.  
\endproof

For a group $G$, a multiplicative character $\chi:G \to \mathbb R$, and a $G$-module $X$ write $X^{G,\chi}=\{\omega \in X: g\omega=\chi(g)\omega , g\in G \}$. We will use the characters $\det:\mathrm O(Q) \to \{\pm 1\}$ and also $\epsilon_0:\mathrm{O}(p,q) \to \{\pm 1\}$ for $\min(p,q)=1$, where the latter is defined according to whether $g$ stabilizes the connected components of the complement of the light cone. 

\begin{Corollary}\label{cor:vertical_forms_summary}  
Assume $1\leq k\leq \frac{n}{2}$.
\begin{displaymath}
\Omega_0:=\Omega_{-\infty,V}^{k,n-k}(V \times \mathbb P(V^*))^{\overline{G_0}}=\begin{cases} \Span\{\omega_{k,0}, \omega_{k,2}\} & \min(p,q)\geq 2, \\  \Span\{\omega_{k,0},\omega_{k,1}, \omega_{k,2}\} &\min(p,q)=1 .\end{cases} 
\end{displaymath}
Moreover, $\Omega_0$ is an $\mathrm O(Q)$-module. We have 
\begin{align*}
 \Omega_0^{\mathrm O(Q),\det} & = \Span\{\omega_{k,0}, \omega_{k,2}\} \text{ if } p,q \geq 1,\\
 \Omega_0^{\mathrm O(Q),\det\cdot\epsilon_0} & =\Span\{\omega_{k,1}\} \text{ if } \min(p,q)=1.
\end{align*}
\end{Corollary}

\proof
The first statement follows immediately from  Proposition \ref{prop:open_orbits_forms_summary}, Proposition \ref{prop_vertical_n_forms} and Lemma \ref{lem:no_other_extensions}. The second statement is simply a restatement of the last part of Proposition \ref{prop:open_orbits_forms_summary}.
\endproof

\subsubsection{Closed forms}

Our next aim is to study which of the forms from Proposition \ref{prop_vertical_n_forms} are closed. 

\begin{Proposition} \label{prop:case_n-k_odd}
Let $1 \leq k \leq \frac{n}{2}$ with $(n-k)$ odd. Then the space $\Omega_{-\infty, V}^{k,n-k}(V\times\mathbb P_+(V^*))^{\overline{G_0}}$ consists of closed forms. 
\end{Proposition}

\proof
Let $\omega \in \Omega_{-\infty, V}^{k,n-k}(V\times\mathbb P_+(V^*))^{\overline{G_0}}$ be non-zero. By Proposition \ref{prop:open_orbits_forms_summary}, $\omega$ is closed on the open orbits, hence $d\omega \in \Omega_{-\infty,M^0}^{k,n-k+1}(V \times \mathbb P(V^*))^{\overline{G_0}}$. To prove that $d\omega=0$, we argue as in Subsection \ref{subsubsec_first_steps}. 

Consider the bundle $F^\alpha_{k,n-k+1}$ corresponding to $D^{k,n-k+1}$. Take $\xi\in M^0$, and define the spaces $U,W,\ldots$ as before. By Lemma \ref{lem:value_of_alpha}, 
$F^\alpha_{k,n-k+1}\cong T^{k,n-k+1,\beta}$ with $\beta=n-k-2\alpha-1$. The corresponding table is as follows
\begin{displaymath}
\begin{array}{c| c | c}
  & \text{determinantal} & \text{symplectic}\\ \hline
  U_{ww} & - & n=2k-1, \beta=0\\
  W_{ww} & - & n=2k-2, \beta=1\\
  U_{uw} & \beta=-1 & n=2k-2,\beta=-1\\
  W_{uw} & - & n=2k-3,\beta=0\\
  U_{wu} & - & n=2k,\beta=-1\\
  W_{wu} & - & n=2k-1, \beta=0\\
  U_{uu} & - & n=2k-1, \beta=-2\\
  W_{uu} & - & n=2k-2,\beta=-1.
 \end{array}
\end{displaymath}

None of the conditions in the table can hold. Indeed under the assumptions on $n$ and $k$, $\beta=(n-k)-2\alpha-1$ must be even, while $n\geq 2k$. Therefore, $d\omega=0$.
\endproof
 
\subsubsection{Proof of Theorem \ref{thm_dimensions}}

\proof
The statement in the cases $k=0, n$ follows from the fact that $\Val_0^{-\infty}(V)=\Val_0(V) = \mathbb C \chi$ and $\Val_n^{-\infty}(V)=\Val_n(V) = \mathbb C \vol$, noting that $\chi$ and $\vol$ are invariant under $\mathrm O(Q)$.

Moreover, using the Alesker-Fourier transform, it is enough to prove the theorem under the assumption $1 \leq k \leq \frac{n}{2}$. 

By Corollary \ref{cor:vertical_forms_summary} we have the following upper bounds. 

\begin{align*}
 \dim\Val_k^{-\infty}(V)^{\mathrm{SO}^+(Q)} & =\dim\Val_k^{-\infty}(V)^{\mathrm{O}(Q)} \leq 2 \text { if } \min(p,q) \geq 2,\\
\dim\Val_k^{+,-\infty}(V)^{\mathrm{SO}^+(Q)} & =\dim\Val_k^{-\infty}(V)^{\mathrm O(Q)} \leq 2 \text { if } \min(p,q) =1,\\
\dim\Val_k^{-,-\infty}(V)^{\mathrm{SO}^+(Q)} & \leq 1 \text { if } \min(p,q)=1.
\end{align*}

If $(n-k)$ is odd, we have equalities by Corollary \ref{cor:vertical_forms_summary} and Proposition \ref{prop:case_n-k_odd}. 

Next, suppose that $(n-k)$ is even. Since $k \leq \frac{n}{2}$, we either have $k \leq p$ or $k \leq q$. In the first case,  the restriction map
$i^*: \Val_k^{-\infty}(\mathbb R^{p+1,q})^{\mathrm O(p+1,q)}\to \Val_k^{-\infty}(\mathbb R^{p,q})^{\mathrm O(p,q)}$ is injective by Corollary \ref{cor:injectivity_of_restriction}. Similarly, the restriction map $i^*: \Val_k^{-\infty}(\mathbb R^{p,q+1})^{\mathrm O(p,q+1)}\to \Val_k^{-\infty}(\mathbb R^{p,q})^{\mathrm O(p,q)}$ is injective in the second case. 

Since $(n+1-k)$ is odd, by the case treated above we have 
\begin{displaymath}
 \dim \Val_k^{-\infty}(\mathbb R^{p+1,q})^{\mathrm O(p+1,q)}=\dim \Val_k^{-\infty}(\mathbb R^{p,q+1})^{\mathrm O(p,q+1)}=2,
\end{displaymath}
and in both cases we conclude that 
\begin{displaymath}
 \dim\Val_k^{-\infty}(\mathbb R^{p,q})^{\mathrm O(p,q)} \geq 2.
\end{displaymath}

This proves the statement in the case of even valuations. Let us finally consider odd valuations, which appear only if $\min(p,q)=1$. The case where $n-k$ is odd was already treated above, so assume $n-k$ is even. Recall from Corollary \ref{cor:vertical_forms_summary}  that both $\omega_{k,0}$, $\omega_{k,2}$ define the sign representation of $\mathrm O(Q)$, and span the full space of invariant differential forms with this property. 

 Since  $\dim\Val_k^{-\infty}(\mathbb R^{p,q})^{\mathrm O(p,q)} =2$, $\omega_{k,0}$ and $\omega_{k,2}$ must be closed. Denote $\omega_+=\omega_{k,j}$ if $(n-k)\equiv j \mod 4$, $j=0,2$. Thus $\omega_+=\omega_{1}+\omega_{2}$ and $\omega_{k,1}=\omega_{1}-\omega_{2}$ with $\omega_{j}:=F^k_{(n-1,1),j} \in \Omega_{-\infty,V,M^0_j}(V\times\mathbb P_+(V^*))^{\overline{G_0}}$ - see equations \eqref{eq:omegak0} - \eqref{eq:omegak1}. 
 
 Since $d\omega_+=0$ and $\omega_1,\omega_2$ have disjoint supports, it follows that $d\omega_1=d\omega_2=0$. Thus also $d\omega_{k,1}=0$, and hence $\omega_{k,1}$ corresponds to a non-zero element in $\Val_k^{-,-\infty}(V)^{\mathrm{SO}^+(Q)}$, which concludes the proof.
\endproof

\section{The image of the Klain map}
\label{sec_klain_functions}

In Proposition \ref{prop_invariant_sections_klain} we determined a basis $\{\kappa_a\}, a=\max(0,k-q),\ldots,\min(k,p)$ of the space of $\mathrm O(Q)$-invariant sections of the Klain bundle $K^{n,k}$. We next determine which invariant sections are Klain functions of even generalized translation-invariant and $\mathrm{SO}^+(p,q)$-invariant valuations, by an inductive argument. 

Let us formulate Theorem \ref{main_thm_klain} in invariant terms. 

\begin{Theorem} \label{thm_Klain_description}
The invariant section $\sum_{a=\max(0,k-q)}^{\min(k,p)} c_a \kappa_a$ of the Klain bundle is in the image of the Klain map $\Kl:\Val_k^{+,-\infty}(\R^{p,q}) \to \Gamma^{-\infty}(K^{n,k})$ if and only if 
\begin{equation} \label{eq_image_klain}
 c_{a+1}+c_{a-1}=0, \quad \max(0,k-q)<a<\min(k,p).
\end{equation}
\end{Theorem}

\proof
We prove the statement by induction on $p+q$. If $k=0$ or $k=n$, the statement is trivial, so let us assume $0<k<n$.

If $\min(p,q) = 0$, the condition is empty and the statement says that the unique invariant section on the Klain bundle is in the image of the Klain map. Since the Klain function of the $k$-th intrinsic volume is a non-zero invariant section, this statement is in fact trivial.   

If $\min(p,q) = 1$, the condition is also empty. By \cite[Thm. 1.3]{alesker_faifman}, the space $\Val_k^{+,-\infty}(\R^{p,q})$ is of dimension $2$. Since there are two open orbits, the statement follows by injectivity of the Klain map.

The induction start is $p=q=2$. Here the statement follows from Proposition \ref{prop_zero_coeffs_of_klain} and Proposition \ref{prop_klain_mu00} in Section \ref{sec_r22} below. The proof of these propositions will be independent of the current section.  

Let us now suppose that $p,q \geq 2$ and that the statement holds for all $(p',q')$ with $p'+q'<p+q$.

Using the quadratic form $Q$ to identify $V \cong V^*$ and $\Dens(V) \cong \C$, the Alesker-Fourier transform is an isomorphism 
\begin{displaymath}
\mathbb F:\Val_k^{+,-\infty}(\R^{p,q})^{\mathrm{SO}^+(p,q)} \to \Val_{n-k}^{+,-\infty}(\R^{p,q})^{\mathrm{SO}^+(p,q)}. 
\end{displaymath}
The Klain function of $\mathbb F \phi$ is the same as the Klain function of $\phi$ composed with the involution $E \to E^Q$, compare \eqref{eq_fourier_even_q} for a precise statement. Under this map, the $\mathrm{SO}^+(p,q)$-orbit $X_{a,b}^k$ is mapped to the $\mathrm{SO}^+(p,q)$-orbit $X_{p-a,q-b}^{n-k}$ and the section $\kappa_a$ of $K^{n,k}$ is mapped to the section $\kappa_{p-a}$ of $K^{n,n-k}$. Hence it is enough to prove \eqref{eq_image_klain} in the case $2k \leq n$.     

Let $\phi \in \Val_k^{+,-\infty}(\R^{p,q})$ and let $\sum c_a \kappa_a$ be its Klain function. By Proposition \ref{prop_restriction_vals}, we may restrict $\phi$ to the space $\R^{p,q-1}$. Each open orbit $X^k_{a,k-a}$ of $\Gr_k(\R^{p,q-1})$ is a subset of the corresponding orbit in $\Gr_k(\R^{p,q})$. The induction hypothesis implies that \eqref{eq_image_klain} holds for all  $\max(0,k-q)+1<a<\min(k,p)$. 

On the other hand, restriction to $\R^{p-1,q}$ and the induction hypothesis imply that \eqref{eq_image_klain} holds for all  $\max(0,k-q)<a<\min(k,p)-1$. Taken together, this yields \eqref{eq_image_klain} for all $\max(0,k-q)<a<\min(k,p)$, which implies the ``only if''-part. 

On the other hand, the Klain map is injective, hence the Klain functions of the invariant valuations span a $2$-dimensional space by Theorem \ref{thm_dimensions}. This proves the ``if'' part. 
\endproof

Recall the involution $j:V\to V$ appearing in the split case.

\begin{Corollary} \label{cor_splitting_split_case}
In the split case $p=q$, there is a canonical direct sum decomposition into one-dimensional spaces: for $1 \leq k \leq n-1$, 
\begin{displaymath}
 \Val_k^{-\infty}(\mathbb R^{p,p})^{\mathrm O(Q)}=\Val_k^{-\infty}(\mathbb R^{p,p})^{\mathrm O(Q),j}\oplus \Val_k^{-\infty}(\mathbb R^{p,p})^{\mathrm O(Q),-j}.
\end{displaymath}
\end{Corollary}

\section{Continuity of invariant valuations}

\label{sec_classification_continuous}

\subsection{Classification of invariant continuous valuations}

\proof[Proof of Theorem \ref{mainthm_continuous_case}]
If $\min(p,q) = 1$, then the statement was shown in \cite[Theorem 1.1]{alesker_faifman}. The idea of the proof in the general case is to use suitable restrictions or push-forwards to reduce to this case. 

Let us assume that $2 \leq q \leq p$.  By Theorem \ref{thm_dimensions} we know that $\Val_k(\mathbb R^{p,q})^{\mathrm{SO}^+(p,q)}=\Val_k(\mathbb R^{p,q})^{\mathrm{O}(p,q)}$. In particular, only even valuations appear and we may apply the results from the previous section.

Let $\phi \in \Val_k(\mathbb R^{p,q})^{\mathrm{O}(p,q)}$ be a continuous valuation. Let $e_j$ be the standard basis, and $\mathbb R^{p,1}=\Span\{e_1,\ldots,e_{p+1}\}$. The subgroup $\mathrm{O}(p,1) \subset \mathrm{O}(p,q)$ acts on $\mathbb R^{p,1}$ and fixes $e_j$ for $j\geq p+2$. Let $\iota:\R^{p,1} \to \R^{p,q}$ and $\pi:\R^{p,q} \to \R^{p,1}$ denote the inclusion and orthogonal projection map respectively. 

If $1 \leq k \leq p-1$, the restriction $\phi_{p,1}:=\iota^*\phi \in \Val_k^{\mathrm{O}(p,1)}(\R^{p,1})$ of $\phi$ to $\mathbb R^{p,1}$ is trivial by the above mentioned result from \cite{alesker_faifman}.  Since $\Kl_{\phi}(\mathbb R^{k,0})=\Kl_{\phi_{p,1}}(\mathbb R^{k,0})=0$, and similarly 
$\Kl_{\phi}(\mathbb R^{k-1,1})=\Kl_{\phi_{p,1}}(\mathbb R^{k-1,1})=0$, it follows by Theorem \ref{thm_Klain_description} that $\Kl(\phi)=0$, and so $\phi=0$.

If $p \leq k \leq n-2$, let $\tilde \phi_{p,1}:=\pi_* \phi \in \Val^{\mathrm{O}(p,1)}_{k-(q-1)}(\R^{p,1})$ be the push-forward under the projection. Using $k-(q-1) \geq p-q+1 \geq 1$, $k-(q-1) \leq (n-2)-(q-1)=p-1$ and \cite{alesker_faifman} again, we deduce that $\tilde \phi_{p,1}=0$. It follows that $\Kl_{\phi}(\R^{k-q+1,q-1})=\Kl_{\tilde \phi_{p,1}}(\mathbb R^{k-q+1,0})=0$ and  $\Kl_{\phi}(\mathbb R^{k-q,q})=\Kl_{\tilde \phi_{p,1}}(\mathbb R^{k-q,1})=0$. By Theorem \ref{thm_Klain_description} we conclude that $\phi=0$.
\endproof

\subsection{KS-continuity of the odd invariant valuations} 

 Let us recall that odd invariant valuations only appear in the Lorentz case, and only for the corresponding connected orthogonal group. We thus assume that $p=n-1,q=1$. By Proposition \ref{prop_orbits_lines}, the open orbits of $\mathrm{SO}^+(Q)$ on $\mathbb{P}_+(V^*)$ are given by $M_1^-,M_2^-$ and $M^+$, and the antipodal map $a$ interchanges $M_1^-$ and $M_2^-$. 

Let $\phi^+_k, \psi_k^+ \in \Val_k^{+,-\infty}(V)^{\mathrm{O}(Q)}$ be the unique valuations with Klain sections $\kappa_{k -1}$, respectively $\kappa_{k}$. Recall that $\Gr^-_{n-1}(V)$ denotes the Grassmann manifold of cooriented hyperplanes in $V$. Fix a positive time direction in $V$.
\begin{Definition}
 The section $\mu_Q\in \Gamma^-(\Gr^-_{n-1}(V), \Dens(E))$ is given by $\mu_Q(E)= \epsilon(E)\vol_Q|_E$, where $\epsilon(E)\in\{\pm 1\}$ measures the time orientation of $V/E$ if $Q|_E>0$, and $\epsilon(E)=0$ otherwise.
	Define $\phi^-_{n-1}\in \Val_{n-1}^{-}(V)^{\mathrm{SO}^+(Q)}$ for a convex body $K$ with smooth boundary by $\phi^-_{n-1}(K)=\int_{\partial K} \mu_Q(T_x\partial K)$, where the tangent planes are cooriented outwards from $K$. For a non-degenerate subspace $E \subset V$ with $\dim E=k+1$, define $s_k^-(E):=(\pi_E)_*(\phi^-_{n-1})\in \Val_k^-(E)$, where $\pi_E:V\to E$ is the projection corresponding to $V=E \oplus E^Q$. 
\end{Definition}

We remark that $\phi^-_{n-1}$ was first considered in \cite{alesker_faifman}. By McMullen's description of $(n-1)$-homogeneous continuous valuations, $\phi^-_{n-1}$ is well-defined. We can write it explicitly for the standard Lorentz structure in $\R^{n-1,1}$ in Euclidean terms as follows:
$\phi_{n-1}^-(K)=\int_{S^{n-1}} \sign(\cos\theta) |\cos2\theta|^{\frac12} d\sigma_K$, where $\cos\theta=x_n$.

It holds that $\phi^-_{n-1}(C)=\phi^+_{n-1}(C)$ for any cone $C$ bounded by the positive part of the light cone and an arbitrary hyperplane. By the base-change theorem, \cite[Theorem 3.5.2]{alesker_fourier}, $s_{k}^-(E)$ is $\mathrm{SO}^+(Q|_E)$-invariant.

We will need the following simple observation.

\begin{Lemma} \label{lemma_push_forward_positive}
Let $\pi:V\to W$ be a linear surjection with $m:=\dim V-\dim W$. Fix positive Lebesgue measures $\sigma_V$, $\sigma_W$ on $V$ and $W$ which allow to identify $\Dens(V^*) \cong \R, \Dens(W^*) \cong \R$. Then $\pi_*:\Val_k(V) \to \Val_{k-m}(W)$ maps positive valuations to positive valuations. 
\end{Lemma}

\proof
The push-forward of a valuation $\tau \in \Val_k(V)$ is defined as follows (compare \cite{alesker_barcelona}). Write $V=W \oplus M$, and take the Lebesgue measure $\sigma_M$ on $M$ such that $\sigma_V=\sigma_W \times \sigma_M$. Let $S$ be a compact convex subset of $M$ of volume $1$. Then 
\begin{displaymath}
\pi_*\tau(K)=\frac{1}{m!} \left.\frac{d^m}{d\epsilon^m}\right|_{\epsilon=0} \tau(K+\epsilon S).
\end{displaymath}
In other words, $\pi_*\tau(K)$ is the coefficient of $\epsilon^m$ in the polynomial $\tau(K+\epsilon S)$. The degree of this polynomial is at most $m$, since $S \mapsto \tau(K+S)$ is an element of $\Val(M)$ for fixed $K$. If $\tau$ is positive, then this polynomial is positive for $\epsilon \geq 0$, which implies that its highest coefficient is positive. 
\endproof

\begin{Lemma}\label{lem:s_k_continuous}
$s_k^-$ extends to a continuous section: $s_k^-\in\Gamma(\Gr_{k+1}(V), \Val_k(E))^{\mathrm{SO}^+(Q)}$.
\end{Lemma}

\proof
Let $X_c\subset \Gr_{k+1}(V)$ denote the $Q$-degenerate subspaces. Extend $s_k^-$ by zero to $X_c$. The base-change theorem \cite[Theorem 3.5.2]{alesker_fourier} implies the $\mathrm{SO}^+(Q)$-invariance of $s_k^-$. It then follows from the locally transitive action of $\mathrm{SO}^+(Q)$ on $\Gr_{k+1}(V)\setminus X_c$ that $s_k^-$ is continuous outside $X_c$. For the global continuity of $s_k^-$, we argue as follows. 
 
Fix a Euclidean structure on $V$. It holds that  $\phi_{n-1}^-(K)=\int_{S^{n-1}} f(\theta)d\sigma_K$, while $\phi_{n-1}^+(K)=\int |f(\theta)|d\sigma_K$ where $f$ is an odd function, and $|f|$ is the Klain section of $\phi^+_{n-1}$. In particular, $\phi^+_{n-1}$ is positive, and for any convex body $K$, $|\phi^-_{n-1}(K)| \leq \phi^+_{n-1}(K)$. By Lemma \ref{lemma_push_forward_positive}, applied to $\tau:=\phi_{n-1}^+ \pm \phi_{n-1}^-$ it follows that for any projection $\pi: V \to W$
\begin{displaymath}
\|\pi_*\phi^-_{n-1}\|\leq \|\pi_*\phi^+_{n-1}\|.
\end{displaymath}
	
Recall that $\Val^+(V)^{\mathrm{O}(Q)}$ consists of KS-continuous valuations.
We next show the equality $(\pi_E)_*\phi_{n-1}^+=i_E^*\phi_k^+$ for non-degenerate $E$,  where $i_E:E \to V$ is the inclusion map. From $\psi_1^+=\mathbb F\phi_{n-1}^+$ and using \cite[Theorem 6.2.1]{alesker_fourier} we deduce that $(\pi_E)_*\phi_{n-1}^+=(\pi_E)_*\mathbb F\psi_1^+=\mathbb F i_E^* \psi_1^+$, and the desired equality follows by comparing the Klain sections of $\mathbb F i_E^* \psi_1^+$ and $i_E^*\phi_k^+$.
	
We define $s_k^+:=(i_E)^*\phi_{k}^+\in \Gamma(\Gr_{k+1}(V), \Val_k({ E}))$, and note that for non-degenerate $E$ one has $\|s_k^-(E)\|=\|(\pi_E)_* (\phi^-_{n-1})\|\leq  \|(\pi_E)_*(\phi_{n-1}^+)\|=  \|s_k^+(E) \|$.
Since $s_k^+$ vanishes at $E\in X_c$, we conclude that $s_k^-$ is continuous.
\endproof

\subsubsection{The case of $k=1$}
We present a short proof that only applies for $k=1$, but contains some of the ingredients that appear in the general case. Define \begin{displaymath}
\phi_1^-:=\mathbb F \phi_{n-1}^-\in \Val_1^{-,-\infty}(V)^{\mathrm{SO}^+(Q)}.
\end{displaymath}
 
\begin{Proposition}
The elements of $\Val_1^{-,-\infty}(V)^{\mathrm{SO}^+(Q)}$ are KS-continuous.
\end{Proposition}

\proof
We claim that 
\begin{displaymath}
s_1^-\in \mathrm{Im}\left(\Sc:\Val_1^{\KS}(V)\to \Gamma(\Gr_{2}, \Val_1(E))\right).
\end{displaymath}

Take an approximate identity $\rho_\epsilon$ on $\mathrm{SO}(n)$.
Recall that the square of the Fourier transform acts as minus the identity on odd valuations. 
Consider the family of projections $\pi_{V^*/H}:V^* \to V^*/H$ for $H \in \Gr_{n-2}(V^*)$ and the dual inclusion $i_{H^\perp}:H^\perp\to V$. Then

\begin{displaymath}
(\pi_{V^*/H})_*\phi_{n-1}^- \ast { \rho}_\epsilon=(\pi_{V^*/H})_* (\phi_{n-1}^- \ast { \rho}_\epsilon) =-\mathbb F_{H^\perp} i_{H^\perp}^*\mathbb F (\phi_{n-1}^- \ast { \rho}_\epsilon) = -\mathbb F_{H^\perp} i_{H^\perp}^* (\phi_1^- \ast { \rho}_\epsilon).
\end{displaymath}

It follows that 
\begin{displaymath}
\mathbb F_{V^*/H} \left((\pi_{V^*/H})_*\phi_{n-1}^- \ast { \rho}_\epsilon\right) =  i_{H^\perp}^*(\phi_1^-\ast { \rho}_\epsilon)
\end{displaymath}
lies in the image of $\Sc$.
Now $(\pi_{V^*/H})_*\phi_{n-1}^-$ is an $\mathrm {SO}^+(Q)$-invariant section of the bundle over $\Gr_{n-2}(V^*)$ with fiber $\Val_1^-(V^*/H)$. The fiberwise Fourier transform $\mathbb F_{V^*/H}$ yields a continuous $\mathrm{SL}(n)$-equivariant map between the following spaces of continuous sections:
\[\Gamma (\Gr_{n-2}(V^*), \Val_1^{-,-\infty}(V^*/H)) \to\Gamma (\Gr_2(V), \Val_1^{-,-\infty}(E))\]
It follows that $s:=\mathbb F_{V^*/H} \left((\pi_{V^*/H})_*\phi_{n-1}^-\right)\in \Gamma (\Gr_2(V), \Val_1^{-,-\infty}(E))^{\mathrm{SO^+(Q)}}$. By Hadwiger's theorem for $\SO(2)$, this section vanishes over subspaces of definite signature. By the Lorentzian analogue \cite{alesker_faifman} of Hadwiger's theorem, it coincides with a multiple of $s_1^-$ over subspaces of signature $(1,1)$. By continuity and equivariance it follows that $s$ is a multiple of $s_1^-$. 
We conclude that for some constant $c\neq 0$,

\begin{displaymath}
c s_1^-\ast { \rho}_\epsilon=  i_{H^\perp}^*(\phi_1^-\ast { \rho}_\epsilon),
\end{displaymath}
belongs to the image of $\Sc$ and  converges in the $\Sc$-norm to a multiple of $s_1^-$. This proves the claim. 

We may thus write $s_1^-=\Sc(\phi)$ for some $KS$-continuous valuations $\phi$. By the injectivity of $\Sc$, $\phi$ is $\mathrm{SO}^+(Q)$-invariant. Since $\dim \Val_1^{-,-\infty}(V)^{\mathrm{SO}^+(Q)}=1$, this completes the proof.
\endproof

\subsubsection{General $k$}
First, let us sketch the plan of the proof. As in the case of $k=1$, it remains to show that $s_k^-$ is in the image of the Schneider map. Of course one expects it to equal the restriction of $\phi_k^-$, namely $\Sc(\phi_k^-)$; however restriction to any fixed subspace is only defined for KS-continuous valuations, which is not known a-priori about $\phi_k^-$; we therefore approximate it by convolving with an approximate identity. 

We represent valuations by closed, vertical generalized forms as in section \ref{sec_dimension_generalized}, where we are able to fit our form of interest into a meromorphic family as was done during its construction. The constructed meromorphic family $\omega_\lambda^{k,n-k}$ is shown to be essentially characterized by the way it transforms under the action of $\mathrm{SO}^+(Q)$. 

Note that restrictions of valuations essentially correspond to push-forwards of generalized forms, while push-forwards of valuations by a projection correspond to restrictions of generalized forms. As restriction of forms is easier to carry out, we use the Alesker-Fourier transform to interchange the two operations. We then carry out the restriction for those values of the complex parameter $\lambda$ where the forms are continuous. We rely on the aforementioned classification through the $\SO^+(Q)$ action description to deduce the value of the restriction without explicit computation (this method fails for push-forward). 

We now give a detailed proof. Let us construct a meromorphic in $\lambda \in \mathbb C$ family of forms $\omega^{k,n-k}_\lambda$ in $\Omega_{-\infty}^{k,n-k} (V \times \mathbb P_+(V^*))^{tr}$ with the property that for any  $\lambda$ which is not a pole, $\omega^{k,n-k}_\lambda$ is vertical and has eigenvalue $(-1)^{n-1}$ under the antipodal map  (for a closed vertical form this is equivalent to saying that the corresponding valuation is odd)  and for $g\in \mathrm{SO}^+(Q)$ satisfies $g^*\omega^{k,n-k}_\lambda=\psi_g(\xi)^{\lambda+\frac{n-k}{2}} \omega^{k,n-k}_\lambda$. We will also classify all such forms for almost every $\lambda$. 

We use the constructions from Subsection \ref{subsec_global_extension}, and take $\omega^{k,n-k}_\lambda$ to correspond to $E^k_{(n-1,1), M_1^-,\lambda}-E^k_{(n-1,1), M_2^-, \lambda}$, as defined in \eqref{eq:def_of_E2}.
By eq. \eqref{eq_equivariance_a_b} and Corollary \ref{cor:invariant_sections} we have   
\begin{displaymath}
g^* \omega_\lambda^{k,n-k}=\psi_g(\xi)^{\lambda+\frac{n-k}{2}} \omega_\lambda^{k,n-k}, \quad g\in \mathrm{SO}^+(Q).
\end{displaymath}
The remaining stated properties are immediate.

Note that while $\omega_\lambda^{n-k,k}$ is not closed in general, the value or residue at $\lambda=-\frac{n-k}{2}$ is closed by Proposition \ref{prop_forms_openorbit}. 

\begin{Proposition}\label{prop:classify_equiv_sections}
For $k \neq \frac{n+1}{2}$ and almost every $\lambda$, the space of vertical forms  in $\Omega_{-\infty}^{k,n-k} (V \times \mathbb P_+(V^*))^{tr}$  that satisfy   \begin{equation}\label{eq:omega_equivariance} g^*\omega=\psi_g(\xi)^{\lambda +\frac{n-k}{2}}\omega \textbf{ }\forall g \in \mathrm{SO}^+(Q) \end{equation} 
and $a^*\omega=(-1)^{n-1}\omega$ is  one-dimensional. For large $\Re\lambda$, this space consists of continuous forms.
\end{Proposition}

\proof
Consider an open $G_0$-orbit $\mathcal O$. As in the proof of Proposition \ref{prop_forms_openorbit}, the space of  translation-invariant forms on $ V \times \mathcal O$ satisfying equation (\ref{eq:omega_equivariance}) is one-dimensional.

Next, since $a^*\omega=(-1)^{n-1}\omega$, we conclude that the restriction of $\omega$ to $ V \times \mathbb P_+(V^*)\setminus M^0$ is unique up to scaling. More precisely, on each of the three open orbits, the space of $\mathrm{SO}^+(Q)$-invariant forms is $1$-dimensional. Since $a$ maps $M_1^-$ to $M_2^-$, the value of $\omega$ on any of these two orbits is determined by the other. On the other hand, $M^+$ is $a$-invariant. Note that $\omega$ is $\mathrm{SO}(n-1)$-invariant by equation \eqref{eq:omega_equivariance}. 

Fix $\xi\in M^+$ and $g\in\mathrm{SO}(n-1)$ such that $g(\xi)=a(\xi)$. Then $(a^*\omega)(0,\xi)=(a^*g^*\omega)(0,\xi)=(g\circ a)^*(\omega (0,\xi))=(-1)^n\omega(0,\xi)$, so that $a^*\omega|_{M^+}=(-1)^n\omega|_{M^+}$. Since $a^*\omega=(-1)^{n-1}\omega$, we conclude that $\omega|_{M^+}=0$. 

For $\lambda$ outside the poles, $\omega^{k,n-k}_\lambda$ constitutes a global extension for $\omega|_{ V \times \mathbb P_+(V^*)\setminus M^0}$. For $\Re\lambda>1$, this form is continuous - this follows from eq. \eqref{eq:def_of_E2}.

Let us finally show that for almost every $\lambda$, there are no global forms with the desired properties that are supported on $M^0$.  This will conclude the proof. 

Write $\lambda'=\lambda+\frac{n-k}{2}$. Fix $\xi\in M^0$, and note that for $g\in \mathrm{Stab}(\xi)$, $\psi_g(\xi)^{\lambda'}=(g|_\xi)^{-2\lambda'}$. Consider $F_{k,l,\lambda'}^\alpha|_{\xi}= F_{k,l}^\alpha|_{\xi}\otimes \xi^{-2\lambda'}$. 
A form as above would define a non-trivial $\mathrm{Stab}(\xi)$-invariant element of $F_{k,l,\lambda'}^\alpha|_\xi$ for some non-negative integer $\alpha$. It follows at once from table (\ref{tbl1}), applied to $T^{k,l,\beta}$ with $l=n-k$ and $\beta=l-2\alpha-2-2\lambda'$ that for almost every $\lambda'$, such an element does not exist.
\endproof

Recall that for $E \in \Gr_{k+1}(V)$, the Alesker-Fourier transform is $\mathbb F_E:\Val_k^{-\infty}(E) \to \Val_1^{-\infty}(E^*) \otimes \Dens(E)$. An element in the image can be identified with an element $\tau \in \Omega^{1,k}_{-\infty,V}(E^* \times \mathbb P_+(E))^{tr} \otimes \largewedge^{k+1}(E)$. Since $a$ acts by multiplication by $(-1)^{k+1}$ on the second factor, odd valuations correspond to those $\tau$ satisfying $a^*\tau=-\tau$.

\begin{Proposition}\label{prop:classify_equiv_sections2}
There is a family of continuous (with respect to the weak topology on the fiber) sections 
\begin{displaymath}
\tau_\lambda\in \Gamma (\Gr_{k+1}(V), \Omega^{1,k}_{-\infty,V}(E^* \times \mathbb P_+(E))^{tr} \otimes \largewedge^{k+1}(E)), \lambda \in \C 
\end{displaymath}
 with the following properties:
\begin{enumerate}
 \item $\tau_\lambda$ is meromorphic in $\lambda$;
 \item $g^*\tau_\lambda=\psi_g^{\lambda+\frac{k}{2}}\tau_\lambda$
 \item $a^*\tau_\lambda =-\tau_\lambda$
\end{enumerate}
Moreover, for almost every $\lambda$, a section satisfying ii)-iii) is unique up to scaling.
\end{Proposition}

\begin{Remark}
Here we consider $\psi_g$ as an element of $\Gamma(\Gr_{k+1}(V), C^\infty(\mathbb P_+(E)))$, which is obtained from $\psi_g(\xi)\in C^\infty(\mathbb P_+(V))$ by restrictions to $\mathbb P_+(E)$ as $E\in \Gr_{k+1}(V)$. The space of sections in the proposition forms a module over this ring.
\end{Remark}

\proof
Note that the space of sections of the bundle under consideration is naturally a module over $\Gamma^\infty(\Gr_{k+1}(V), C^\infty(\mathbb P_+(E)))$. Note also that there is a natural inclusion $C^\infty(\mathbb P_+(V))\to \Gamma^\infty(\Gr_{k+1}(V), C^\infty(\mathbb P_+(E)))$ given by restrictions.

Fix $E_0$ of signature $(k,1)$. We take $\tau_\lambda(E_0)$ to be the form $\omega^{1,k}_\lambda$ on $E_0^*\times\mathbb P_+(E_0)$, twisted by the $Q$-defined density on $E_0$.  Extend $\tau_\lambda$ to a global section by equivariance   (i.e. by $g^*\tau_\lambda=\psi_g(\xi)^{\lambda+\frac{k}{2}} \tau_\lambda$), by requiring that $a^*\tau_\lambda=-\tau_\lambda$, and that $\tau_\lambda$ vanishes on the planes of signature different from $(k,1)$. 

It is straightforward to check that for any fixed $\lambda$, this is a well-defined assignment of an element of the fiber over all $E\in\Gr_{k+1}(V)$. It will be later shown to form a continuous section. The restriction of $\tau_\lambda$ to any of the open $\mathrm{SO}^+(Q)$-orbits is a smooth section.
To see that $\tau_\lambda$ is continuous and meromorphic in $\lambda$, decompose
\begin{displaymath}
\tau_\lambda=\rho(\cos2\theta(\xi))\tau_\lambda+(1-\rho(\cos2\theta(\xi)))\tau_\lambda
\end{displaymath} 
with $\rho$ as in eq. \eqref{eq:E_cutoff_decomposition}.

The second summand is simply $(1-\rho(\cos2\theta(\xi)))|\cos2\theta(\xi)|^\lambda\tau_0$, which is evidently an analytic family of smooth sections. The first summand we may decompose as $ |\cos2\theta|_-^{\lambda}\alpha_E(\xi)- |\cos2\theta|_-^{\lambda-1}\beta_E(\xi)$, corresponding to the $A$ and $B$ summands in the definition of $\omega^{1,k}_\lambda$ as in \eqref{eq:Ekpq}. Thus $\alpha_E$ and $\beta_E$ are smooth sections of our bundle, multiplied by a meromorphic in $\lambda$ family of sections $\Gamma(\Gr_{k+1}(V), C^{-\infty}(\mathbb P_+(E)))$. This shows that $\tau_\lambda$ is a meromorphic family of continuous sections.
Finally, applying Proposition \ref{prop:classify_equiv_sections} with $k$ replaced by $1$ and $V$ replaced by $E^*$,  uniqueness for almost every $\lambda$ follows. 
\endproof

\begin{Proposition}
 The valuations in $\Val_k^{-,-\infty}(V)^{\mathrm{SO}^+(Q)}$ are KS-continuous. Moreover, the image of the Schneider embedding $\Sc(\Val_k^{-,\KS}(V)^{\mathrm{SO}^+(Q)})$
is spanned by $s_k^-$.
\end{Proposition}

\proof
Let $\Omega_{cont}^{n-k,k}(V^*\times \mathbb P_+(V))^{tr}$ denote the space of continuous forms.

Consider the $\mathrm{GL}(V)$-equivariant map 
\begin{displaymath}
 j^*:\Omega_{cont}^{n-k,k}(V^*\times \mathbb P_+(V))^{tr}\otimes \largewedge^{top}(V)\to \Gamma (\Gr_{k+1}(V), \Omega_{cont}^{1,k}(E^*\times \mathbb P_+(E))^{tr}\otimes\largewedge^{k+1}(E)).
\end{displaymath}

It is given by the Gelfand transform $j^*\omega(E):=(\pi_E)_*i_E^*\omega$ with respect to the double fibration 
\begin{displaymath}
\xymatrix{&  V^* \times \mathbb P_+(E)\ar[dl]_{\tilde i_E} \ar[dr]^{\tilde\pi_E}&\\ V^* \times \mathbb P_+(V) & &E^* \times \mathbb P_+(E) }
\end{displaymath}
induced by the natural projection $\pi_E:V^* \to E^*$ and inclusion $i_E:\mathbb P_+(E) \to \mathbb P_+(V)$. 

This operation is clearly well defined for continuous forms. Moreover, vertical forms are mapped to sections of vertical forms - this is easiest to see by considering the closely related double fibration \begin{displaymath}
\xymatrix{&  V^* \times E\ar[dl]_{\hat i_E} \ar[dr]^{\hat\pi_E}&\\ V^* \times V & &E^* \times E }
\end{displaymath}

Denoting the canonical $1$-forms $\alpha_E$ and $\alpha_V$ on $E^*\times E$ and $V^*\times V$, respectively, it clearly holds that $\hat i_E^*\alpha_V=\hat{\pi}_E^*\alpha_E$. Thus one can choose contact forms on $V^*\times\mathbb P_+(V)$ and $E^* \times \mathbb P_+(E)$ with the same property, and we fix such a choice of forms. It is then easy to see that for arbitrary $\eta$, 

\begin{displaymath}
(\tilde \pi_E)_* \tilde i_E^* (\alpha_V\wedge \eta)=(\tilde\pi_E)_*(\tilde \pi_E^*\alpha_E\wedge \tilde i_E^*\eta)=\alpha_E\wedge (\tilde\pi_E)_* \tilde i_E^*\eta
\end{displaymath}
is vertical. 

For closed vertical forms, $j^*$ corresponds to the push-forward of the associated valuation by $\pi_E$. For a precise statement and proof, see \cite[Proposition 3.2.3]{alesker_intgeo}.

Fix an orientation on $V$, thus identifying $\largewedge^{top}V=\Dens(V)$, and recall that $\vol_Q$ is the density on $V$ induced by $Q$. The contracted restriction $j^*(\omega_\lambda^{n-k, k}\otimes \vol_Q)$ is well defined for sufficiently large values of $\Re\lambda$, as $\omega_\lambda^{n-k,k}$ is continuous for those values. 

Using the $\mathrm{SO}^+(Q)$-equivariance and sign-reversal by the antipodal map, it follows from Proposition \ref{prop:classify_equiv_sections2} that for almost every $\lambda$ with $\Re\lambda$ sufficiently large, $j^*(\omega_\lambda^{n-k, k}\otimes \vol_Q)$ is a multiple of $\tau_\lambda$. Let us write $j^*(\omega_\lambda^{n-k,k}\otimes \vol_Q)=\tau'_\lambda=c_\lambda\tau_\lambda$. 

Choose an approximate identity ${ \rho}_\epsilon$ on $\mathrm{SL}(V)$. Since $j^*$ is $ \mathrm{SL}(V)$-equivariant, it follows that for almost every $\mathrm{Re} \lambda$ sufficiently large  
\begin{equation} \label{eq_push_forward_form}
j^*(\omega_\lambda^{n-k,k} \otimes \vol_Q\ast { \rho}_\epsilon )=c_\lambda\tau_\lambda \ast { \rho}_\epsilon.
\end{equation}
As the left hand side is meromorphic in $\lambda\in\mathbb C$, we may use this equality to extend $c_\lambda$ meromorphically to $\lambda\in\mathbb C$, see Lemma \ref{lemma_meromorphic_extension} below.
It follows that this identity also holds for $\lambda= -\frac{k}{2}$ when $k \equiv 1 \mod 2$, or for the residues at $\lambda=-\frac{k}{2}$ when $k \equiv 0 \mod 2$. Also, as $\tau_\lambda$ and $\omega_\lambda^{n-k,k}$ are either both analytic or have a simple pole at $\lambda=-\frac{k}{2}$, we deduce that $c_\lambda$ is analytic at $\lambda=-\frac{k}{2}$.

One has the natural inclusion 
\begin{displaymath}
 \Val^{-\infty}_{n-k}(V^*) \otimes \Dens(V) \subset \Omega_{-\infty}^{n-k,k}(V^*\times \mathbb P_+(V))^{tr} \otimes \largewedge^{top}(V),
\end{displaymath}
corresponding to the closed vertical forms (compare the beginning of Section \ref{sec_dimension_generalized}). Under this correspondence, either the value or residue of $\omega_\lambda^{n-k,k}$ at $\lambda=-\frac{k}{2}$ corresponds to some valuation $\psi^-_{n-k}\in \Val_{n-k}^{-,-\infty}(V^*)^{\mathrm{SO}^+(Q)}$, while $\tau_\lambda(E)$ corresponds to some $\psi_{1,E}^-\in \Val_1^{-,-\infty}(E^*)^{\mathrm{Stab}(E^*)}$. Equation \eqref{eq_push_forward_form} then yields an identity on push-forwards of valuations, namely 
$(\pi_{E^*})_*(\psi^-_{n-k}\ast { \rho}_\epsilon)= c_0 \psi^-_{1,E}\ast { \rho}_\epsilon$.

Applying the Alesker-Fourier transform on both sides yields 
\begin{displaymath}
(i_E)^*(\mathbb F\psi_{n-k}^-\ast{ \rho}_\epsilon)=c_0 \cdot (s^-_{k}\ast { \rho}_\epsilon)(E). 
\end{displaymath}
By $\mathrm{SO}^+(Q)$-invariance of both sides, $c_0$ is in fact independent of $E$. 

This proves that $\Sc(\mathbb F\psi_{n-k}^-)=c_0s_k^-$, and since $s_k^-$ is continuous by Lemma \ref{lem:s_k_continuous}, this shows that 
\begin{displaymath}
\mathbb F \psi_{n-k}^- \in \Val_{k}^{-,\KS}(V)^{\mathrm{SO}^+(Q)}.  
\end{displaymath}

It remains to recall that $\mathbb F:\Val_{n-k}^{-,-\infty}(V^*)^{\mathrm{SO}^+(Q)} \to \Val_{k}^{-,-\infty}(V)^{\mathrm{SO}^+(Q)}$ is an isomorphism.
\endproof

Let us define $\phi_k^-\in\Val_{k}^{-,-\infty}(V)^{\mathrm{SO}^+(Q)}$ as the unique valuation with $\Sc(\phi_k^-)=s_k^-$. It is then evident  using Proposition \ref{prop_restriction_vals} that given an inclusion of Lorentz spaces $i:U\to V$, one has the identity $i^*\phi_{V,k}^{-}=\phi^-_{U,k}$.

\begin{Lemma} \label{lemma_meromorphic_extension}
Let $f, g:\mathbb C \to V$ be meromorphic functions with values in some quasi-complete locally convex topological vector space $V$. Assume that $g \neq 0$ and that for an open set $D\subset\mathbb C$ it holds that $f(s)=c(s)g(s)$ for some function $c:D\to \mathbb C$ (which is necessarily meromorphic). Then $c$ extends to a meromorphic function on $\C$ and the above identity holds everywhere.
\end{Lemma}

\proof
Fix $s_0$ where $g$ is holomorphic and non-zero. Take $\phi\in V^*$ such that $\phi(g(s_0)) \neq 0$. 

Then 
\begin{displaymath}
c(s):=\frac{\phi(f(s))}{\phi(g(s))}
\end{displaymath} 
defines a meromorphic extension for $c$. Now the identity $f=cg$ is between two globally meromorphic sections on an open set $D$, so it holds everywhere.
\endproof

\section{Invariant Crofton measures for $\R^{2,2}$}
\label{sec_r22}

In this section, we study in detail the case of $\R^{2,2}$. We will write down explicitly the invariant generalized sections of the Crofton bundle and compute the Klain functions of the corresponding generalized valuations. This will also complete the proof of Theorem \ref{thm_Klain_description}. 

Let $V:=\R^{2,2} \cong \R^4$ with its standard basis $e_1,\ldots,e_4$. We fix the standard Euclidean form $P$ and use the indefinite form $Q(x_1e_1+x_2e_2+x_3e_3+x_4e_4)=x_1^2+x_2^2-x_3^2-x_4^2$.  

\subsection{Construction of invariant Crofton measures of degree $1$ and $3$}

Let us construct an invariant generalized measure on $\Gr_1(\R^{2,2})$. The case of $\Gr_3(\R^{2,2})$ is similar or can be directly reduced to this one by applying the Alesker-Fourier transform. 

The orbits of $\mathrm{O}(2,2)$ on $\Gr_1(\R^{2,2})$ are $X_{1,0}=\{\cos(2\theta)>0\}, X_{0,1}=\{\cos(2\theta)<0\}$ (open) and $X_{0,0}=\{\cos(2\theta)=0\}$ (closed). 

According to Proposition \ref{prop_meromorphic_cos2theta}, there is a meromorphic family of generalized measures $|\cos(2\theta)|^{\lambda}$ on $\Gr_1(\R^{2,2})$. Using the Euclidean trivialization as in Corollary \ref{cor_invariant_crofton_klain_sections}, we obtain two meromorphic families of generalized measures $m_{1,0}^{(\lambda)}, m_{0,1}^{(\lambda)}$ on the open orbits.

We set $m_{1,0}:=m_{1,0}^{\left(-\frac52\right)}$ and $m_{0,1}:=m_{0,1}^{\left(-\frac52\right)}$. By Corollary \ref{cor_invariant_crofton_klain_sections}, $m_{1,0}$ and $m_{0,1}$ are invariant generalized Crofton measures. We denote the corresponding invariant valuations by $\phi_{1,0},\phi_{0,1} \in \Val_3^{-\infty}(\R^{2,2})$.

\begin{Proposition}
 The Klain functions of $\phi_{1,0}, \phi_{0,1}$ are given by 
\begin{align*}
 \Kl_{\phi_{1,0}} & =-\frac{8}{3} \pi \kappa_1, \\                                                                                           
 \Kl_{\phi_{0,1}} & = -\frac{8}{3} \pi \kappa_2.
\end{align*}
\end{Proposition}

\proof
We will use the notation from Proposition \ref{prop_meromorphic_cos2theta} and its proof in the special case $p=q=2$. 

We know that $\Kl_{\phi_{1,0}}=a \kappa_1+b \kappa_2$, where $\kappa_1,\kappa_2$ are defined in Proposition \ref{prop_invariant_sections_klain}. In our situation, $\kappa_1,\kappa_2$ are just the restrictions of $|\cos(2\theta)|^{\frac12}$ to the two open orbits. To compute $a,b$, we evaluate the valuations at the $3$-planes $e_1^P$ and $e_3^P$, noting that $\kappa_1(e_1^P)=\kappa_2(e_3^P)=1, \kappa_1(e_3^P)=\kappa_2(e_1^P)=0$.  

Let us first take $E:=e_1^P$. The cosine of the angle between the line given by $(z_1,z_2,\theta)$ and $E$ is $h=|\cos(\theta) \mathrm{Re}(z_1)|=\cos(\theta) |\cos(\tau_1)|$, where $z_j=\cos \tau_j+i \sin \tau_j, j=1,2$. We thus get 
\begin{align*}
  \langle m_{1,0}^{(\lambda)},\langle e_1^P,\cdot\rangle\rangle & = \int_0^{2\pi} \int_0^{2\pi} \int_0^{\frac{\pi}{4}} |\cos(2\theta)|^\lambda \sin(\theta)\cos^2(\theta) |\cos(\tau_1)| d\theta d\tau_2 d\tau_1\\
  & = 8\pi \int_0^{\frac{\pi}{4}} |\cos(2\theta)|^\lambda \sin(\theta)\cos^2(\theta) d\theta,
\end{align*}
and it follows (using the lemma below) that 
\begin{displaymath}
a=\Kl_{\phi_{1,0}}(e_1^P)=\langle m_{1,0},\langle e_1^P,\cdot\rangle \rangle= -\frac{8}{3} \pi.
\end{displaymath}
Similarly, 
\begin{align*}
  \langle m_{1,0}^{(\lambda)},\langle e_3^P,\cdot\rangle\rangle & = \int_0^{2\pi} \int_0^{2\pi} \int_0^{\frac{\pi}{4}} |\cos(2\theta)|^\lambda \sin^2(\theta)\cos(\theta) |\cos(\tau_2)| d\theta d\tau_2 d\tau_1\\
  & = 8\pi \int_0^{\frac{\pi}{4}} |\cos(2\theta)|^\lambda \sin(\theta)^2 \cos(\theta) d\theta,
\end{align*}
and hence 
\begin{displaymath}
b=\Kl_{\phi_{1,0}}(e_3^P)=\langle m_{1,0},\langle e_3^P,\cdot\rangle \rangle= 0.
\end{displaymath}

The computation for $\phi_{0,1}$ is similar. 
\endproof

In the proof, we used the following lemma. 

\begin{Lemma}
The meromorphic extensions of 
 \begin{displaymath}
  \int_0^{\frac{\pi}{4}} |\cos(2\theta)|^\lambda \sin \theta \cos^2 \theta d\theta \quad \text{and} \quad \int_{\frac{\pi}{4}}^{\frac{\pi}{2}} |\cos(2\theta)|^\lambda \sin^2 \theta \cos \theta d\theta
 \end{displaymath}
at $\lambda=-\frac52$ equal $-\frac13$. 

The meromorphic extensions of 
 \begin{displaymath}
  \int_0^{\frac{\pi}{4}} |\cos(2\theta)|^\lambda \sin^2 \theta \cos \theta d\theta \quad \text{and} \quad \int_{\frac{\pi}{4}}^{\frac{\pi}{2}} |\cos(2\theta)|^\lambda \sin \theta \cos^2 \theta d\theta
 \end{displaymath}
at $\lambda=-\frac52$ equal $0$.
\end{Lemma}

\proof
Set for $\mathrm{Re} \lambda>1$,
\begin{displaymath}
 I(\lambda):=\int_0^1 s^\lambda \sqrt{1+s} ds.
\end{displaymath}

Substituting $s:=\cos(2\theta)$ we obtain
\begin{displaymath}
  \int_0^{\frac{\pi}{4}} \cos(2\theta)^\lambda \sin \theta \cos^2 \theta d\theta = \frac{1}{4\sqrt{2}} \int_0^1 s^\lambda \sqrt{1+s} ds=\frac{1}{4\sqrt{2}} I(\lambda).
\end{displaymath}

Using integration by parts, one obtains the recurrence relation 
\begin{displaymath}
 (2\lambda+5)I(\lambda+1)+(2\lambda+2) I(\lambda)=4\sqrt{2}. 
\end{displaymath}
Since $I(\lambda)$ is analytic at $\lambda=-\frac32$, the result follows by plugging in $\lambda:=-\frac52$. 

The other values are obtained in a similar way.
\endproof

\subsection{The geometry of $\Gr_2(\mathbb R^{2,2})$}\label{sub:n=4}
\label{subsec_gr2_r22}

It is well-known that there is a double cover $S^2 \times S^2 \to \Gr_2 \R^4$. Explicitly, let $(w_1,w_2,w_3),(z_1,z_2,z_3) \in S^2$. Set 
\begin{align*}
 x_{12} & :=\frac{w_1+z_1}{2} \quad x_{34}:=\frac{w_1-z_1}{2} \quad x_{13}:=-\frac{w_2+z_2}{2}\\
x_{24} & :=\frac{w_2-z_2}{2} \quad x_{14}:=\frac{w_3+z_3}{2} \quad x_{23}:=\frac{w_3-z_3}{2}. 
\end{align*}
Then the $2$-vector 
\begin{equation} \label{eq_2_vector}
 \tau:=\sum_{1 \leq i<j\leq 4} x_{ij} e_i \wedge e_j
\end{equation}
is simple and defines the $2$-plane $E:=\{v \in \R^4: v \wedge \tau=0\}$. We will be identifying a subspace with one of its lifts to $S^2\times S^2$ whenever this cannot lead to ambiguity.

Under the action of $\mathrm{O}(2,2)$, the orbits on $\Gr_2(\R^4)$ lift as follows: $X_{2,0}=\{z_1^2+w_1^2>1,z_1w_1>0\}$, $X_{0,2}=\{z_1^2+w_1^2>1,z_1w_1<0\}$, $X_{1,1}=\{z_1^2+w_1^2<1\}$, $X^2_{1,0}=\{z_1^2+w_1^2=1,z_1w_1>0\}$, $X^2_{0,1}=\{z_1^2+w_1^2=1,z_1w_1<0\}$, $X^2_{0,0}=\{z_1=\pm 1, w_1=0\}\cup \{z_1=0, w_1=\pm 1\}$.

Let $s_i$ be the map which changes the sign of the $i$-th coordinate. The notation will be the same for $s_i$ acting on $\R^N$ for different values of $N$ as well as for the restriction on the unit spheres. 

Note that both $P$ and $Q$ induce a quadratic form on $\largewedge^2 \R^4$, denoted by the same letter. The vectors $e_i \wedge e_j, 1 \leq i<j\leq 4$ form an orthonormal basis for both $P$ and $Q$, with $Q(e_1 \wedge e_2)=Q(e_3 \wedge e_4)=1$ and $Q(e_1 \wedge e_3)=Q(e_1 \wedge e_4)=Q(e_2 \wedge e_3)=Q(e_2 \wedge e_4)=-1$. The Euclidean orthogonal complement lifts to $(z,w)^P=(-z,w)$, while the $Q$-orthogonal complement is $(z,w)^Q=(s_1z, s_2s_3w)$.

We identify $\R^{2,2}=\C^2$ by $(x_1,x_2,x_3,x_4) \mapsto (x_1+ix_3,x_2+ix_4)$. Multiplication by the complex unit $i$ corresponds to $(z,w) \mapsto (s_1s_3z,w)$. In particular, the complex Grassmannian $\Gr^\C_1(\C^2) \subset \Gr^\R_2(\R^4)$ is covered by $(0,\pm 1,0) \times S^2$. Note that multiplication by $i$ is not in $\mathrm{O}(Q)=\mathrm{O}(2,2)$ but satisfies $i^*Q=-Q$.

Let $u,v$ be an orthonormal basis of $E \in \Gr^\R_2(\C^2)$. The K\"ahler angle $\theta_\C$ of $E$ is defined by $\cos^2 \theta_\C(E)=P(u,iv)^2, \theta_\C(E) \in [0,\frac{\pi}{2}]$. 

\begin{Lemma} \label{lemma_sigma_on_s2s2}
Let $E \in \Gr_2(\R^4)$ be represented by $(z,w)$. Then 
\begin{align}
  \cos 2\theta(E) & =z_1^2+w_1^2-1 \nonumber\\
\cos^2 \theta_\C(E) & =z_2^2. \label{eq_kaehler_angle}
\end{align}
\end{Lemma}

\proof
Write $\tau=\sum_{1 \leq i<j\leq 4} x_{ij} e_i \wedge e_j=u \wedge v$ with $u,v \in V$ a $P$-orthonormal basis of $E$. Then 
\begin{align*}
 \cos 2\theta(E) & =\det \left(\begin{matrix} Q(u,u) & Q(u,v)\\ Q(v,u) & Q(v,v)\end{matrix}\right)\\
& =  Q(u \wedge v, u \wedge v)\\
& = x_{12}^2+x_{34}^2-x_{13}^2-x_{14}^2-x_{23}^2-x_{24}^2\\
& = z_1^2+w_1^2-1
\end{align*}
and
\begin{align*}
\cos^2 \theta_\C(E) & =\det \left(\begin{matrix} 0 & \pm \cos \theta_\C(E) \\ \mp \cos \theta_\C(E) & 0\end{matrix}\right)\\
& =\det \left(\begin{matrix} P(u,iu) & P(u,iv) \\ P(v,iu) & P(v,iv)\end{matrix}\right)\\
& =P (u \wedge v,iu \wedge iv)\\
& =P(\tau,i\tau)\\
& =x_{13}^2+x_{24}^2+2x_{12}x_{34}+2x_{14}x_{23}\\
& =z_2^2.
\end{align*}
\endproof

\subsection{Construction of invariant Crofton measures of degree $2$}

Since we are in the split case $p=q=2$, we can use the notion of $j$-even and $j$-odd elements introduced in Definition \ref{def_split}. By Corollary \ref{cor_splitting_split_case}, there is a non-trivial decomposition
\begin{displaymath}
\Val_2^{-\infty}(\R^{2,2})^{\mathrm O(Q)}=\Val_2^{-\infty}(\R^{2,2})^{\mathrm O(Q),j} \oplus \Val_2^{-\infty}(\R^{2,2})^{\mathrm O(Q),-j}.
\end{displaymath}

Similarly, the action of $j$ on the Grassmannian $\Gr_2(\R^{2,2})$ is an involution and we can speak of $j$-even and $j$-odd invariant (generalized) Crofton measures.

Let $x_+^\mu$ be the well-known generalized function on the real line \cite{hoermander_pde1}. This is a meromorphic function in $\mu$ with poles at $\mu=-1,-2,\ldots$ and corresponding residues 
\begin{displaymath}
 \Res_{\mu=-m}x_+^\mu=\frac{(-1)^{m-1}}{(m-1)!} \delta^{(m-1)}.
\end{displaymath}

In the following, we will be making use of generalized measures on a manifold with corners, namely, $\mathcal M^{-\infty}(X)$ with $X\subset\mathbb R^d$ a parallelepiped. 
In the literature such generalized measures are sometimes called supported distributions. By definition, those are generalized measures given in some neighborhood $U\subset \mathbb R^d$ of $X$ and supported inside $X$.

\begin{Lemma} \label{lemma_generalized_tan}
There exists a meromorphic family of generalized measures 
\begin{displaymath}
\tan^{2\lambda+2}(t)dt \in \mathcal M^{-\infty}\left[0,\frac{\pi}4\right].
\end{displaymath}
It has simple poles at $\lambda=-\frac m2$, $m\geq 3$, and 
\begin{displaymath} 
\langle \Res_{\lambda=-\frac m2} \tan ^{2\lambda+2}(t)dt,\psi\rangle=\frac{ (-1)^{m-1}}{2(m-3)!} \left.\frac {d^{m-3}}{dx^{m-3}}\right|_{x=0}\frac{\psi(\arctan x)}{1+x^2}.
\end{displaymath} 
\end{Lemma}

\proof
It suffices to set  
\begin{displaymath}
\int_0^{\pi/4}\tan^{2\lambda+2}(t) \psi(t)dt := \int_0^{1} x_+^{2\lambda+2} \frac{\psi(\arctan x)}{1+x^2}dx, 
\end{displaymath}
using the generalized function $x_+^{2\lambda+2}$ \cite[Section 3.2]{hoermander_pde1}. 
\endproof

Define for $\Re(\lambda )>0$ and $(a,b) \in \{(2,0), (1,1), (0,2)\}$ the generalized measures $m_{a,b}^{(\lambda)}$ on $\Gr_2(\R^{2,2})$ by 
 \begin{displaymath}
  \phi \mapsto \int_{X_{a,b}^2} \phi(E) |w_1^2+z_1^2-1|^\lambda dE.
 \end{displaymath}
 
\begin{Proposition} \label{prop_construction_crofton_measures}
\begin{enumerate}
 \item $m_{a,b}^{(\lambda)}$ admits a meromorphic extension to $\lambda\in\mathbb C$ with simple poles at $\lambda=-\frac{m}{2}, m\geq 2$ if $(a,b) \in \{(2,0),(0,2)\}$ and with simple poles at $\lambda=-1,-2,\ldots$ if $(a,b)=(1,1)$. The measure $m_-^{(\lambda)}=m_{2,0}^{(\lambda)}-m_{0,2}^{(\lambda)}$ has poles at $\lambda=-m$ with $m\geq1$.
\item Denote $m_{0,0}:=\Res_{\lambda=-\frac52} m_{2,0}^{(\lambda)}$. Then $m_{0,0}$ is an $\mathrm{O}(2,2)$-invariant and  $j$-even generalized Crofton measure supported on $X_{0,0}^2$.
\item Denote $m_+:=m_{1,1}^{\left(-\frac52\right)}$. Then $m_+$ is an $\mathrm{O}(2,2)$-invariant, $j$-even generalized Crofton measure supported on $\overline{X_{1,1}^2}$.
\item Denote $m_-:=m_-^{\left(-\frac52\right)}$. Then $m_-$ is an $\mathrm{O}(2,2)$-invariant,  $j$-odd  generalized Crofton measure supported on $\overline{X_{2,0}^2} \cup \overline{X_{0,2}^2}$.
\end{enumerate}
\end{Proposition}

We denote the corresponding generalized translation and $\mathrm{O}(2,2)$-invariant valuations by $\phi_{0,0},\phi_+,\phi_-$. 

\begin{Remark}
One can moreover show that the space of invariant generalized Crofton measures of degree $2$ is spanned by $m_+,m_-,m_{0,0}$. Since we do not need this stronger statement in the following, we omit the (rather technical) details.  

We also refer to \cite{faifman_crofton} for a more general  study of Crofton measures for $\OO(p,q)$.
\end{Remark}

\proof[Proof of the proposition]
Define two operators 
\begin{align*}
Q :C^\infty([-1,1]^2 \setminus \{0\}) & \to C^\infty([-1,1]^2 \setminus \{0\}),\\
 \Phi & \mapsto  \frac{z_1\frac{\partial\Phi}{\partial z_1}+w_1\frac{\partial\Phi}{\partial w_1}}{z_1^2+w_1^2},\\
 R  : C^\infty([-1,1]^2 \setminus \{0\}) & \to C\left(\left[0,\frac{\pi}{4}\right]\right),\\
 \Phi & \mapsto \Phi(1,\tan t)+\Phi(\tan t,1)+\Phi(-1,-\tan t)+\Phi(-\tan t,-1).
\end{align*}

 For $\phi\in C^\infty (\Gr_2(\R^{2,2}))$, denote by $\tilde \phi$ its lift to $S^2 \times S^2$. Then put $\Phi(z_1,w_1)=\int \int _{\mathrm{SO}(2) \times \mathrm{SO}(2)} \tilde \phi(g_1(z), g_2(w))dg_1dg_2$, where the left resp. right copy of $\mathrm{SO}(2)$ is the stabilizer in $\mathrm{SO}(3)$ of $z_1$ resp. $w_1$, and $dg$ is the invariant probability measure. The orbit $X^2_{2,0}$ corresponds to $A_+:=\{z_1^2+w_1^2>1,\max(|z_1|,|w_1|)\leq 1, z_1w_1>0\}$, and 
  we have
\begin{displaymath}
\int_{\Gr_2(\R^{2,2})} \phi dm_{2,0}^{(\lambda)}= \frac{1}{4} \left( \int_{A_+}(z_1^2+w_1^2-1)^\lambda \Phi (z_1,w_1)dz_1dw_1\right). 
\end{displaymath}

Thus it suffices to define the meromorphic extension of the generalized measures $\tilde m_{2,0}^{(\lambda)}= (z_1^2+w_1^2-1)^\lambda dz_1dw_1$ on   $A_+$. Those are well-defined for $\Re(\lambda)>0$, and we will show that a meromorphic extension to $\Re(\lambda)>\sigma_0$ implies a meromorphic extension to $\Re(\lambda)>\sigma_0-1$. 

Set 
\begin{displaymath}
 I:=\int_{A_+} \Phi d\tilde m_{2,0}^{(\lambda)}.
\end{displaymath}

By using polar coordinates $z_1:=\cos(t) \cosh(s), w_1:=\sin(t) \cosh(s)$, we obtain that  
\begin{align*}
I & = \int_{0}^\frac{\pi}{4}\int_0^{\arcosh\frac{1}{\cos t}} \sinh^{2\lambda+1} s \cosh s \Phi(\cosh s\cos t,\cosh s\sin t)dsdt \\
& \quad + \int_\frac{\pi}{4}^\frac{\pi }{2}\int_0^{\arcosh\frac{1}{\sin t}} \sinh^{2\lambda+1} s \cosh s \Phi(\cosh s\cos t,\cosh s\sin t)dsdt\\
& \quad +\int_\pi^\frac{5 \pi }{4} \int_0^{\arcosh\frac{-1}{\cos t}} \sinh^{2\lambda+1} s \cosh s \Phi(\cosh s\cos t,\cosh s\sin t)dsdt \\
& \quad + \int_\frac{5 \pi}{4}^\frac{3 \pi }{2}\int_0^{\arcosh\frac{-1}{\sin t}} \sinh^{2\lambda+1} s \cosh s \Phi(\cosh s\cos t,\cosh s\sin t)dsdt.
\end{align*}

Integrating by parts of the inner integrals gives $I=I_1+I_2$ with 
\begin{align*}
 I_1 & = \frac{1}{2\lambda+2} \Bigg[\int_0^\frac{\pi}{4} \tan(t)^{2\lambda+2} \Phi(1,\tan t) dt  + \int_\frac{\pi}{4}^\frac{\pi}{2} \cot(t)^{2\lambda+2} \Phi(\cot t,1) dt\\
& \quad \int_\pi^\frac{5 \pi}{4} \tan(t)^{2\lambda+2} \Phi(-1,-\tan t) dt+\int_\frac{5 \pi}{4}^\frac{3 \pi}{2} \cot(t)^{2\lambda+2} \Phi(-\cot t,-1) dt\Bigg]\\
& = \frac{1}{2\lambda+2} \int_0^\frac{\pi}{4} \tan(t)^{2\lambda+2} [\Phi(1,\tan t)+\Phi(\tan t,1)+\Phi(-1,-\tan t)+\Phi(-\tan t,-1)] dt\\
& = \frac{1}{2\lambda+2} \int_0^\frac{\pi}{4} \tan(t)^{2\lambda+2} (R\Phi)(t) dt,
\end{align*}
 and 
\begin{align*}
 I_2 & = - \frac {1}{2\lambda+2} \Bigg[\int_{0}^{\pi/4} \int_0^{\arcosh \frac{1}{\cos t}}\sinh ^{2\lambda +3}s\left(\frac{\partial\Phi}{\partial z_1}\cos t+\frac{\partial\Phi}{\partial w_1}\sin t\right) dsdt\\
 & \quad + \int_\frac{\pi}{4}^\frac{\pi }{2} \int_0^{\arcosh \frac{1}{\sin t}}\sinh ^{2\lambda +3}s\left(\frac{\partial\Phi}{\partial z_1}\cos t+\frac{\partial\Phi}{\partial w_1}\sin t\right)\\
 & \quad + \int_\pi^\frac{5 \pi }{4} \int_0^{\arcosh \frac{-1}{\cos t}}\sinh ^{2\lambda +3}s\left(\frac{\partial\Phi}{\partial z_1}\cos t+\frac{\partial\Phi}{\partial w_1}\sin t\right)\\
 & \quad + \int_\frac{5 \pi}{4}^\frac{3 \pi }{2}\int_0^{\arcosh \frac{-1}{\sin t}}\sinh ^{2\lambda +3}s\left(\frac{\partial\Phi}{\partial z_1}\cos t+\frac{\partial\Phi}{\partial w_1}\sin t\right)\Bigg]\\
 & = - \frac {1}{2\lambda+2} \int_{A_+} (Q\Phi) d\tilde m_{2,0}^{(\lambda+1)}.
\end{align*}

The term $I_1$ is meromorphic on $\C$ by Lemma \ref{lemma_generalized_tan}, while $I_2$ is meromorphic in $\Re(\lambda)>\sigma_0-1$ by the hypothesis.

We thus have shown the first claim and the equation  
\begin{equation}
 \label{eqn:recursion}\int_{A_+} \Phi d\tilde m_{2,0}^{(\lambda)}=\frac{1}{2\lambda+2} \int_0^\frac{\pi}{4} \tan(t)^{2\lambda+2} (R\Phi)(t)dt-\frac{1}{2\lambda+2} \int_{A_+} (Q\Phi) d\tilde m_{2,0}^{(\lambda+1)}. 
\end{equation}

The orbit $X_{0,2}^2$ corresponds to $A_-:=\{z_1^2+w_1^2>1,\max(|z_1|,|w_1|)\leq 1, z_1w_1<0\}$. A similar computation as for $A_+$ gives 
\begin{displaymath}
 \int_{A_-} \Phi d\tilde m_{0,2}^{(\lambda)}=\frac{1}{2\lambda+2} \int_{-\frac{\pi}{4}}^0 \tan(t)^{2\lambda+2} (R\Phi)(t)dt-\frac{1}{2\lambda+2} \int_{A_-} (Q\Phi) d\tilde m_{0,2}^{(\lambda+1)}
\end{displaymath}
and therefore for $\lambda=-\frac{2m+1}{2}$,
\begin{align*}
 \int_{A_+} \Phi d\tilde m_{2,0}^{(-\frac{2m+1}{2})}-\int_{A_-} & \Phi d\tilde m_{0,2}^{(-\frac{2m+1}{2})} =\frac{1}{2\lambda+2} \int_{-\frac{\pi}{4}}^\frac{\pi}{4} \tan(t)^{-2m+1} (R\Phi)(t)dt\\
 & \quad -\frac{1}{2\lambda+2} \left(\int_{A_+} (Q\Phi) d\tilde m_{2,0}^{(\lambda+1)}-\int_{A_-} (Q\Phi) d \tilde m_{0,2}^{(\lambda+1)}\right). 
\end{align*}

Moreover, since the generalized measure $|\tan t|^{2\lambda+2} dt$ considered on the two intervals $[0,\pi/4]$ and $[-\pi/4,0]$  has equal residues of opposite signs at $\lambda=-(2m+1)/2$ for $m\geq 1$, it follows that one has the generalized measure $(\tan t)^{-2m+1}dt$ on $[-\pi/4,\pi/4]$. The claim for $m_-^{(\lambda)}$ follows.

Iteration of the formula above gives us 
\begin{align*}
  \int_{A_+} \Phi d\tilde m_{2,0}^{(\lambda)} & =\frac{1}{2\lambda+2} \int_0^\frac{\pi}{4} \tan(t)^{2\lambda+2} (R\Phi)(t)dt\\
  & \quad -\frac{1}{(2\lambda+2)(2\lambda+4)}  \int_0^\frac{\pi}{4} \tan(t)^{2\lambda+4} (RQ\Phi)(t)dt \\
  & \quad +\frac{1}{(2\lambda+2)(2\lambda+4)}  \int_{A_+} (Q^2\Phi) d \tilde m_{2,0}^{(\lambda+2)}. 
\end{align*}

Using Lemma \ref{lemma_generalized_tan} and the fact that $\tilde m_{2,0}^{(\lambda)}$ has no pole at $-\frac12$ we get 
\begin{align*} 
\int_{\Gr_2(\R^{2,2})} \phi dm_{0,0} & =\frac14 \Res_{\lambda=-\frac52} \int_{A_+} \Phi d\tilde m_{2,0}^{(\lambda)} \\
& = \frac{-1}{48} \left.\frac{d^2}{dx^2}\right|_{x=0} \frac{R\Phi(\arctan x)}{1+x^2}-\frac{1}{24}  RQ\Phi(0) \\ & =-\frac{1}{48} (R\Phi)''(0)+\frac{1}{24}  (R\Phi(0)-RQ\Phi(0)).
\end{align*}

Let us denote by $\Phi_4$ the average of $\Phi$ over the $4$ rotational symmetries of the square $[-1,1] \times [-1,1]$. Then we may rewrite the above equation as 

\begin{equation} \label{eq_int_m00}
\int_{\Gr_2(\R^{2,2})} \phi dm_{0,0}=\left[-\frac{1}{12}\frac{\partial ^2\Phi_4}{\partial w_1^2}+\frac{1}{6}  \left(\Phi_4-\frac{\partial \Phi_4}{\partial z_1}\right)\right]_{(1,0)}.
\end{equation}

An alternative way to write this equation is 
\begin{displaymath}
 \int_{\Gr_2(\R^{2,2})} \phi dm_{0,0}=\left[ \frac16 \pi^*\Phi_4+\frac1{12}\left(\frac{\partial^2 \pi^*\Phi_4}{\partial z_2^2}+\frac{\partial^2 \pi^*\Phi_4}{\partial z_3^2}-\frac{\partial^2 \pi^*\Phi_4}{\partial w_1^2}\right)\right]_{(1,0,0)\times(0,1,0)},
\end{displaymath}
where $\pi: S^2\times S^2\to [-1,1]^2$ is given by $\pi(z,w)=(z_1,w_1)$.

Let us finally study the case $(a,b)=(1,1)$. Let $B:=\{z_1^2+w_1^2<1\}$. Defining $\Phi$ as above, we have for $\Re(\lambda)>0$
\begin{displaymath}
 \int_{\Gr_2(\R^{2,2})} \phi d\tilde m_{1,1}^{(\lambda)}=\frac14 \int_B (1-z_1^2-w_1^2)^\lambda \Phi(z_1,w_1) dz_1 dw_1.
\end{displaymath}

Write this in polar coordinates:
\begin{align}
\int_{\Gr_2(\R^{2,2})} \phi d\tilde m_{1,1}^{(\lambda)} & =\frac14 \int_0^1 (1-r^2)^\lambda r \int_0^{2\pi}\Phi(r\cos t,r\sin t)dt dr \nonumber \\
& = \frac18 \int_0^1 (1-s)^\lambda \int_0^{2\pi}\Phi(\sqrt s\cos t,\sqrt s\sin t)dt ds \nonumber\\
& = \frac18 \int_0^1 x^\lambda \int_0^{2\pi}\Phi(\sqrt{1-x}\cos t,\sqrt{1-x}\sin t)dt dx.  \label{eq_measure_m_+}
\end{align}

Note that $\Psi(x):=\int_0^{2\pi}\Phi(\sqrt{1-x}\cos t,\sqrt{1-x}\sin t)dt$ is continuous in $x$ and smooth away from $x=1$.
Therefore, one may apply the generalized measure $x^\lambda dx$ to $\Psi(x)$, which is meromorphic in $\lambda$ with poles in the negative integers.

It remains to prove $\mathrm{O}(2,2)$-invariance of the Crofton measures just constructed.
We identify the generalized measures $m_{a,b}^{(\lambda)}$ constructed above with generalized functions over subsets of the Grassmannian using the Euclidean trivialization. It then follows that for $\lambda>0$ and $g\in \mathrm{O}(2,2)$ one has 
\begin{displaymath}
g^*m_{a,b}^{(\lambda)}=\psi_g^\lambda \cdot m_{(a,b)}^{(\lambda)}, 
\end{displaymath}
where $\psi_g$ is the function from Proposition \ref{prop_jacobian}. Note that 
\begin{displaymath}
\psi_g^\lambda:\mathbb C \to C^\infty (\Gr_2(\R^{2,2}))
\end{displaymath}
is an entire function of $\lambda$. Both sides are meromorphic in $\lambda$, and so this identity holds for all $\lambda$ where $m_{a,b}^{(\lambda)}$ is analytic, as well as for the residues at the simple poles. Taking $\lambda=-\frac52$, the statement follows from Corollary \ref{cor_invariant_crofton_klain_sections}.
\endproof 

\subsection{Computation of Klain functions}

We next will compute the Klain functions of the valuations $\phi_{0,0},\phi_+,\phi_-$. 

\begin{Proposition} \label{prop_zero_coeffs_of_klain}
\begin{enumerate}
 \item The Klain functions corresponding to $\phi_+$ and $\phi_{0,0}$ vanish on $X^2_{2,0} \cup X^2_{0,2}$. 
 \item The Klain function corresponding to $\phi_-$ is $\frac13(\kappa_2-\kappa_0)$. 
\end{enumerate}
\end{Proposition}

\proof
We evaluate those valuations on the unit square $S_{12}$ in $\R^{2,0} \subset \R^{2,2}$. 

The projection function of this square is easily seen to be equal $|x_{12}|=\frac{|z_1+w_1|}{2}$. We thus apply our formulas to $\Phi(z_1,w_1):=\frac{|z_1+w_1|}{2}$. Note that $\Phi$ is $1$-homogeneous, and so $Q\Phi=\frac{z_1+w_1}{2(z_1^2+w_1^2)}$, $Q^2 \Phi=-\frac{z_1+w_1}{2(z_1^2+w_1^2)^2}$.

The support of $m_{0,0}$ is disjoint from the singular support of the projection function of $S_{12}$, which is given by $M_\Phi=\{z_1+w_1=0\}$. We may thus use Proposition \ref{prop_Klain_continuous_section_computation}, and by \eqref{eq_int_m00} we obtain 
\begin{displaymath}
 \phi_{0,0}(S_{12}) = -\frac{1}{48} \left.\frac{d^2}{dx^2}\right|_{x=0} \frac{2(x+1)}{x^2+1}-\frac{1}{12}=0.
\end{displaymath}

We claim that the wavefronts of $m_+$ and $m_-$ are disjoint from the conormal bundle of $M_\Phi$: The singular support of $m_\pm$ is $\{z_1^2+w_1^2=1\}$, which intersects $M_\Phi$ at a union of two tori, $\pm P=\{z_1=\pm\frac1{\sqrt 2}, z_2=\mp\frac1{\sqrt 2} \}$ (which is one  torus on the unoriented Grassmannian). 
The function $\sigma_2:=z_1^2+w_1^2-1$ is regular at a neighborhood of $P$. Therefore, one can choose a neighborhood $U$ of $P$ and a diffeomorphism $\Psi_{\pm}: U\to U_0$ where $U_0\subset \mathbb R^4$ is some open subset, s.t. $m_\pm$ is mapped to a multiple of the generalized function $(x_4)_\pm^{-5/2}\in C^{-\infty}(U_0)$, and $\Psi_\pm^* x_4=\sigma_2$. Thus the wavefront of $m_\pm$ near $P$ is contained in the conormal bundle to $M_\mu=\{\sigma_2=0\}$. Now for $b\in P$, $T_{b} M_\Phi=\{dz_1+dw_1=0\}$ and $T_{b} M_\mu=\{dz_1-dw_1=0\}$ are different hyperplanes in $T_{b}\Gr_2(\R^{2,2})$. It follows that their annihilators intersect trivially, as claimed.

Thus we may again apply Proposition \ref{prop_Klain_continuous_section_computation} to compute the Klain embeddings of $\phi_{\pm}$.

We now compute for positive values of $\lambda$
\begin{align*}
 \left\langle m_{1,1}^{(\lambda)}, \frac{|z_1+w_1|}{2}\right\rangle & =\frac12 \int_0^1 x^{\lambda} \int_0^{2\pi} \sqrt{1-x}|\cos t+\sin t|dtdx\\
 & =2\sqrt 2\int_0^1 x^{\lambda}\sqrt{1-x}\ dx\\
 & =2\sqrt 2B\left(\lambda+1,\frac32\right),
\end{align*}
where $B$ is Euler's Beta function. By uniqueness of meromorphic extension, we may take $\lambda=-\frac52$ and get
\begin{displaymath}
 \Kl_{\phi_+}(\mathbb R^{1,2})=\left\langle m_+, \frac{|z_1+w_1|}{2}\right\rangle=2\sqrt 2B\left(-\frac32, \frac32\right)=0.
\end{displaymath}

For $\phi_-$, we have
\begin{align*}
\Kl_{\phi_-}(\mathbb R^{2,0})=\frac12\langle m_-, |z_1+w_1|\rangle = \frac18\langle \sign(z_1w_1)(z_1^2+w_1^2-1)^{-5/2}, |z_1+w_1|\rangle. 
\end{align*}

Denote $T:=\{0\leq |w_1|\leq z_1\leq1\}$. Since both sides of the pairing are invariant under the antipodal map and symmetric in $z_1,w_1$
\begin{align*}
\phi_-(S_{12}) & =\frac12\int_{T}|z_1^2+w_1^2-1|^{-5/2}\sign(z_1w_1)(z_1+w_1)dz_1dw_1 \\
& = -\frac16\left( \int_{-1}^{1}x^{-3}\frac{(1+x)dx}{1+x^2}
 + \int_{-1}^{1}x^{-1}\frac{(1+x)dx}{(1+x^2)^2}+\right.\\
&\left. \quad +\int_T \sign(z_1w_1)(z_1^2+w_1^2-1)^{-1/2}\frac{z_1+w_1}{(z_1^2+w_1^2)^2}dz_1dw_1 \right).
\end{align*}

Denote the three integrals by $I_1$, $I_2$, $I_3$. Since integrals over odd function vanish on $[-1,1]$, we get 
\begin{align*}
 I_1 & = \int_{-1}^1 \left(\frac{1}{x^2}-\frac{1}{1+x^2}\right)\ dx=-\frac\pi2-2\\
 I_2 & =\int_{-1}^{1} \frac{dx}{(1+x^2)^2}=\int_{-\pi/4} ^{\pi/4}\cos^2 t dt =\frac{\pi}{4}+\frac12.
\end{align*}

In the integral $I_3$, we substitute $z_1=\cosh s\cos t$, $w_1=\cosh s\sin t$, $dz_1dw_1=\cosh s\sinh sdsdt$ and obtain 
\begin{align*}
I_3 & = \int_{-\pi/4}^{\pi/4}(\cos t+\sin t)\sign(t) \int_0^{\arcosh\frac{1}{\cos t}}\frac{1}{\cosh ^2s}ds dt\\
& = \int_{-\pi/4}^{\pi/4}(\cos t+\sin t)\sign(t) \tanh \left(\arcosh\frac{1}{\cos t}\right)dt\\
& = 2\int_0^{\pi/4}\sin^2 tdt \\
& =\frac{\pi}{4}-\frac12.
\end{align*}

Thus 
\begin{displaymath}
 \phi_-(S_{12})=-\frac16 \left(-2\right)=\frac{1}{3}.
\end{displaymath}
Since $\phi_-$ is $j$-odd, it follows that $\Kl_{\phi_-}$ vanishes on $X^2_{1,1}$, concluding the proof.
\endproof

\begin{Lemma} 
Let $m_\R$ be the $\mathrm{SO}(4)$-invariant probability measure on $\Gr_2(\R^{2,2})=\Gr_2(\R^4)$ and let $m_\C$ be the $\mathrm{U}(2)$-invariant probability measure on the complex Grassmannian $\Gr_1(\mathbb C^2) \subset \Gr_2\R^4$. Then 
\begin{equation} \label{eq_int_kappa_real}
 \int_{\Gr_2 \R^4} \kappa_i(E) dm_\R(E)=\begin{cases} \frac{1-\log(2)}{6} & i=0,2 \\ \frac{\pi}{6}   & i=1, \end{cases}
\end{equation}
and 
\begin{equation} \label{eq_int_kappa_complex}
 \int_{\Gr_2 \R^4} \kappa_i(E) dm_\C(E)=\begin{cases} 0 & i=0,2 \\ \frac{\pi}{4}   & i=1. \end{cases}
\end{equation}
\end{Lemma}

\proof
The orbit $X^2_{1,1}$ lifts to $\{z_1^2+w_1^2<1\}$.  The orbit $X^2_{2,0}$ lifts to $\{z_1^2+w_1^2>1, z_1w_1>0\}$ and similarly $X^2_{0,2}$ lifts to $\{z_1^2+w_1^2>1, z_1w_1<0\}$.

The $\mathrm{SO}(3) \times \mathrm{SO}(3)$-invariant probability measure on $S^2 \times S^2$ projects to the $\mathrm{SO}(4)$-invariant probability measure on $\Gr_2(\R^4)$.

If $E \in X_{1,1}^2$ lifts to $(z,w)$, then, by Lemma \ref{lemma_sigma_on_s2s2}, 
\begin{displaymath}
\kappa_1(E) = | \cos 2\theta(E)|^\frac12= |w_1^2+z_1^2-1|^\frac{1}{2}.
\end{displaymath}

Using Archimedes' theorem, it follows that 

\begin{align*}
 \int_{\Gr_2\R^4} \kappa_1(E) dm_\R(E) & =  \frac{1}{(4\pi)^2}  \int_{z_1^2+w_1^2<1}\sqrt{1-z_1^2-w_1^2}dp_1dp_2\\
& = \frac{1}{4} \int_{z_1^2+w_1^2<1} \sqrt{1-z_1^2-w_1^2} dz_1 dw_1\\
& = \frac{\pi}{6}.
\end{align*}

For the orbits $X_{20}^2 \cup X_{02}^2$, a similar computation yields the integral 
\begin{align*}
 I & =\int_{\Gr_2\R^4} (\kappa_0+\kappa_2)(E) dm_\R(E) \\
 & = \frac {1}{4}\int_{z_1^2+w_1^2>1, -1 \leq z_1,w_1 \leq 1} \sqrt{z_1^2+w_1^2-1} dz_1 dw_1.
\end{align*}

By symmetry considerations, it is enough to consider the first octant. We substitute $x_1=\cos t \cosh s$, $x_2=\sin t \cosh s$. Then $dx_1\ dx_2=\sinh s \cosh s dt\ ds$ and 
\begin{align*}
I & =2 \int_0^{\frac {\pi}{4}}\int_{0}^{\arcosh \frac{1}{\cos t}}\sqrt {\cosh ^2 s-1}\sinh s\cosh s dsdt \\
&=\frac{2}{3} \int_0^{\frac {\pi}{4}}\int_{0}^{\arcosh \frac{1}{\cos t}} d(\sinh^3 s)dt \\
& = \frac{2}{3}\int_0^{\frac {\pi}{4}}\sqrt{\frac {1}{\cos ^2 t}-1}^3 dt \\
& = \frac{2}{3}\int_0^{\frac {\pi}{4}} \tan^3 t\ dt \\
& = \frac{2}{3} \int _0^{\frac\pi 4} \frac{1-\cos^2 t}{\cos^3 t} \sin t dt \\
& = \frac{2}{3}\int_{1/\sqrt 2}^1 \frac{1-r^2}{r^3}dr \\
& =\frac{1-\log(2)}{3}. 
\end{align*}

Since $j$ interchanges the orbits $X_{20}^2$ and $X_{02}^2$ as well as the functions $\kappa_0$ and $\kappa_2$, and preserves the invariant measure, we conclude that 
\begin{displaymath}
 \int_{\Gr_2(\R^4)} \kappa_0 (E) dm_\R(E)=\int_{\Gr_2(\R^4)} \kappa_2 (E) dm_\R(E)=\frac{1-\log(2)}{6},
\end{displaymath}
as claimed. This finishes the proof of \eqref{eq_int_kappa_real}.

Since the supports of $m_\C$ and $\kappa_0, \kappa_2$ intersect at a set of measure zero, the cases $i=0,2$ of \eqref{eq_int_kappa_complex} are trivial. 

The complex Grassmannian is the image of $(0,1,0) \times S^2$, and the $\mathrm{SO}(3)$-invariant probability measure on the second copy of $S^2$ projects to the $\mathrm{U}(2)$-invariant measure. Therefore
\begin{align*}
\int_{\Gr_2(\R^4)} \kappa_1(E) dm_\C(E) & =\frac{1}{4\pi} \int_{S^2} \sqrt{1-w_1^2} dp\\
& = \frac{1}{4\pi} \int_0^{2\pi} \int_0^\pi \sin^2 \theta\ d\theta\ d\phi\\
& = \frac{\pi}{4}.
\end{align*}
\endproof

\begin{Lemma}
Let $\mu_{2,0}^\C, \mu_{2,1}^\C \in \Val_2^{\mathrm{U}(2)}(\C^2)$ be the hermitian intrinsic volumes \cite{bernig_fu_hig}. Then 
\begin{align}
\int_{\Gr_2} \Kl_{\mu_{2,0}^\C} dm_{0,0} & =0 \label{eq_int_kl_mu20_m00}\\
\int_{\Gr_2} \Kl_{\mu_{2,1}^\C} dm_{0,0} & =\frac16 \label{eq_int_kl_mu21_m00}\\
\int_{\Gr_2} \Kl_{\mu_{2,0}^\C} dm_+ & = 0 \label{eq_int_kl_mu20_m+} \\
\int_{\Gr_2} \Kl_{\mu_{2,1}^\C} dm_+ & = -\frac{\pi}{6}. \label{eq_int_kl_mu21_m+}
\end{align}
\end{Lemma}

\proof
Let us compute the pull-back of the Klain functions of $\mu_{2,0}^\C,\mu_{2,1}^\C$ under the map $S^2 \times S^2 \to \Gr_2$. Fix $E \in \Gr_2$ lifting to $(z,w)$. By \eqref{eq_kaehler_angle}, the K\"ahler angle $\theta$ of $E$ is given by $\cos^2 \theta = z_2^2$. It follows that $\Kl_{\mu_{2,0}^\C}=1-\cos^2\theta=1-z_2^2, \Kl_{\mu_{2,1}^\C}=\cos^2\theta=z_2^2$.  

A simple computation shows that the function $\Phi$ from the proof of Proposition \ref{prop_construction_crofton_measures} is given by $\frac{1+z_1^2}{2}$ (in the case of $\mu_{2,0}^\C$) and $\frac{1-z_1^2}{2}$ (in the case of $\mu_{2,1}^\C$). The displayed equations thus follow from \eqref{eq_int_m00} and \eqref{eq_measure_m_+}.
\endproof

\begin{Proposition} \label{prop_klain_mu00}
 The Klain functions of $\phi_+$ and $\phi_{0,0}$ are given by 
\begin{align*}
\Kl_{\phi_+} & =-\frac13 \kappa_1,\\
\Kl_{\phi_{0,0}} & =\frac{1}{3\pi} \kappa_1.
\end{align*}
In particular, $\phi_+=-\pi \phi_{0,0}$.
\end{Proposition}

\proof

As $\phi_+$ and $\phi_{0,0}$ are $j$-even, the Klain functions are of the form $\Kl_{\phi_{0,0}}=a\kappa_2+b \kappa_1+a \kappa_0, \Kl_{\phi_+}=a' \kappa_2+b' \kappa_1+a'\kappa_0.$ 
We could use Proposition \ref{prop_zero_coeffs_of_klain} to conclude that $a=a'=0$. Instead, we give a unified computation of all the coefficients using the hermitian intrinsic volumes. 

The valuation $\frac14 \mu_{2,0}^\C+\frac12 \mu_{2,1}^\C$ has as Crofton measure the $\mathrm{U}(2)$-invariant probability measure $m_\C$ on the complex Grassmannian, as follows from the results in \cite{bernig_fu_hig}. By \eqref{eq_int_kl_mu20_m00}-\eqref{eq_int_kl_mu21_m+} and \eqref{eq_int_kappa_complex}, we obtain 
\begin{align*}
 \frac{1}{12} & =\int_{\Gr_2} \Kl_{\frac14 \mu_{2,0}^\C+\frac12 \mu_{2,1}^\C} dm_{0,0}=\int_{\Gr_2} \Kl_{\phi_{0,0}} dm_\C=b \frac{\pi}{4},\\
 \frac{-\pi}{12} & =\int_{\Gr_2} \Kl_{\frac14 \mu_{2,0}^\C+\frac12 \mu_{2,1}^\C} dm_+=\int_{\Gr_2} \Kl_{\phi_+} dm_\C=b' \frac{\pi}{4},
\end{align*}
i.e. $b=\frac{1}{3\pi}, b'=- \frac13$.

Similarly, the valuation $\frac13 \mu_2^\R=\frac13 \mu_{2,0}^\C+\frac13 \mu_{2,1}^\C$ (which is a multiple of the second intrinsic volume) has as Crofton measure the $\mathrm{SO}(4)$-invariant probability measure $m_\R$ on the real $2$-Grassmannian. Therefore, using \eqref{eq_int_kl_mu20_m00}-\eqref{eq_int_kl_mu21_m+} and \eqref{eq_int_kappa_real},
\begin{align*}
 \frac{1}{18} & =\int_{\Gr_2} \Kl_{\frac13 \mu_{2,0}^\C+\frac13 \mu_{2,1}^\C} dm_{0,0}=\int \Kl_{\phi_{0,0}} dm_\R=a \frac{1-\log 2}{3}+ \frac{1}{18},\\
 \frac{-\pi}{18} & =\int_{\Gr_2} \Kl_{\frac13 \mu_{2,0}^\C+\frac13 \mu_{2,1}^\C} dm_+=\int \Kl_{\phi_+} dm_\R=a' \frac{1-\log 2}{3}- \frac{\pi}{18},
\end{align*}
i.e. $a=0, a'=0$. 
 \endproof

\begin{Corollary}
 Every $\mathrm{O}(2,2)$-invariant generalized translation-invariant valuation admits an invariant generalized Crofton measure.  
\end{Corollary}

The same is true for $\mathrm O(p,q)$ with $\min(p,q)=1$, see \cite{alesker_faifman}. Recently, the statement was generalized to arbitrary $p,q$, based on the results of the present paper \cite[Theorem 2]{faifman_crofton}.  

\begin{appendix}

\section{Invariant generalized sections}
\label{appendix}

 \subsection{Supports of invariant sections}
The following is a technical lemma which allows to bound from above the dimension of the space of invariant sections with a given support under the action of a group $G$. It is well-known, and presented here (with proof) in the form that is adapted to our needs.

Let $X$ be a smooth manifold, and $E$ a smooth vector bundle over $X$.
For any $\alpha \geq 0$ and a locally closed submanifold $Y\subset X$, define the vector bundle $F^\alpha_Y$ over $Y$ with fiber

\begin{displaymath}
F^\alpha_Y|_y=\Sym^\alpha(N_yY) \otimes \Dens^*(N_yY) \otimes E|_y.                                                                                                                                                                                                                                                                                                                                                                                                                                                                                                                                                                                    \end{displaymath}

For a closed submanifold $Y\subset X$, recall the subspaces $\Gamma^{-\infty,\alpha}_{Y}(X,E)\subset\Gamma^{-\infty}_{Y}(X,E)$ of all generalized sections supported on $Y$ with differential order not greater that $\alpha\geq 0$ in directions normal to $Y$. 
One then has a natural isomorphism 
\begin{displaymath}
\Gamma_Y^{-\infty,\alpha}(X,E)/\Gamma_Y^{-\infty,\alpha-1}(X,E)\cong\Gamma^{-\infty}(Y,F_Y^\alpha).
\end{displaymath}

For precise definitions, see \cite[Section 4.4]{alesker_faifman}.

Now let a smooth Lie group $G$ act on $X$ in such a way that there are finitely many orbits, all of which are locally closed submanifolds. We will assume that $E$ is a $G$-vector bundle. If $Y \subset X$ is a $G$-invariant locally closed submanifold, then $F^\alpha_Y$ is naturally a $G$-bundle. If $Y$ is in fact a closed submanifold then $\Gamma^{-\infty,\alpha}_Y(X,E)^G$ form a filtration on $\Gamma^{-\infty}_Y(X,E)^G$.

\begin{Lemma}
\label{lem:LocallyClosedOrbits}
Let $Z\subset X$ be a closed $G$-invariant subset.
Decompose $Z=\bigcup_{j=1}^J Y_j$ where each $Y_j$ is a $G$-orbit. Then
\begin{displaymath}
 \dim \Gamma^{-\infty}_Z(X,E)^G \leq \sum_{\alpha=0}^\infty \sum _{j=1}^J \dim \Gamma^\infty(F^\alpha_{Y_j})^{G}.
\end{displaymath}

More generally, if $Z_1\subset Z_2$ are two $G$-invariant closed subsets of $X$ then
\begin{displaymath}
\dim \Gamma^{-\infty}_{Z_2}(X,E)^G \leq \dim \Gamma^{-\infty}_{Z_1}(X,E)^G+\sum_{\alpha=0}^\infty\sum _{Y_j \subset Z_2\setminus Z_1} \dim \Gamma^\infty(F^\alpha_{Y_j})^{G}. 
\end{displaymath}
\end{Lemma}

\begin{Remark}
Fixing $y_j \in Y_j$, we evidently have 
\begin{displaymath}
 \dim \Gamma^{-\infty}(Y_j, F^\alpha_{Y_j}) ^G= \dim \Gamma^{\infty}(Y_j, F^\alpha_{Y_j}) ^G=\dim \left(F^{\alpha}_{Y_j}|_{y_j}\right)^{\Stab(y_j)}.
\end{displaymath}
\end{Remark}

\proof
Let us start by making two observations. First, let $Y_j \subset Z$ be any $G$-orbit with relative boundary $B:=\mathrm{cl}(Y_j)\setminus Y_j$. Now $Y_j \subset X \setminus B$ is a closed $G$-orbit in $X\setminus B$, and by \cite[Proposition 4.9]{alesker_faifman} we get 
\begin{displaymath}
 \dim \Gamma^{-\infty}_{Y_j}(X \setminus B,E)^G \leq \sum_{\alpha=0}^\infty \dim \Gamma^\infty(Y_j,F^\alpha_{Y_j})^G.
\end{displaymath}

The same upper bound remains valid also for $\dim \Gamma^{-\infty}_{Y_j}(X \setminus C,E)^G$ for any $G$-invariant, closed $C \subset X$ s.t. $B \subset C$, $Y_j \cap C=\emptyset$.

Second, given two $G$-invariant sets $C \subset D$ which are closed in $X$, one has the exact sequence of $G$-modules
\begin{displaymath}
 0 \to \Gamma_C^{-\infty}(X,E) \to \Gamma_D^{-\infty}(X,E) \to \Gamma_{D \setminus C}^{-\infty}(X\setminus C,E)
\end{displaymath}

and therefore also the exact sequence of $G$-invariants
\begin{displaymath}
 0\to \Gamma_C^{-\infty}(X,E)^G\to \Gamma_D^{-\infty}(X,E)^G \to \Gamma_{D\setminus C}^{-\infty}(X\setminus C,E)^G.
\end{displaymath}

Thus one has the inequality
\begin{equation}\label{eqn:second_observation}
 \dim \Gamma_D^{-\infty}(X,E)^G\leq \dim \Gamma_C^{-\infty}(X,E)^G+\dim\Gamma_{D\setminus C}^{-\infty}(X\setminus C,E)^G.
\end{equation}

Now let us prove the statement of the Lemma. Let $Y\subset Z$ be a closed $G$-invariant subset. We will show 
\begin{displaymath}
 \dim \Gamma^{-\infty}_Y(X,E)^G \leq \sum_\alpha \sum_{Y_j \subset Y} \dim \Gamma^\infty(F^\alpha_{Y_j})^G
\end{displaymath}
by induction on $d=\dim Y$.

For the smallest admissible dimension, $Y$ must be a closed $G$-orbit, and the statement follows from \cite[Proposition 4.9]{alesker_faifman}.
Assume we have shown the claim for any closed subset of dimension smaller than $d$, and let $Y \subset Z$ be closed and $G$-invariant of dimension $d$. Let $A \subset Y$ be the maximal $G$-invariant subset of $Y$ of dimension at most $d-1$, which must be closed by our assumptions. 

Denote $B:=Y \setminus A=\bigcup_{i=1}^k B_i$, each $B_i$ being a $G$-orbit which is open in $Y$. Note that $A$ contains the relative boundary of $B_i$ for all $i$. Denoting $A_0:=A, A_r:=A\cup \bigcup_{i=1}^r B_i$, it follows that $A_r \subset X$ is closed. We will show by induction on $r=0,1,\ldots,k$ that 
\begin{displaymath}
 \dim \Gamma^{-\infty}_{A_r}(X,E)^G \leq \sum_\alpha \sum_{Y_j \subset A_r} \dim \Gamma^\infty(F^\alpha_{Y_j})^{G}.
\end{displaymath}

Indeed, for $r=0$, this holds since $\dim A<d$. Assume 
\begin{displaymath}
 \dim \Gamma^{-\infty}_{A_{r-1}}(X,E)^G \leq \sum_\alpha \sum_{Y_j\subset A_{r-1}} \dim \Gamma^\infty(F^\alpha_{Y_j})^{G}.
\end{displaymath}

Since $B_r=\mathrm{cl}(B_r)\setminus A_{r-1}$ is a $G$-orbit with relative boundary $B=\mathrm{cl}(B_r) \setminus B_r \subset A_{r-1}$, it follows by the first observation that for $j$ with  $B_r=Y_j$ one has 
\begin{displaymath}
 \dim\Gamma^{-\infty}_{\mathrm{cl}(B_r) \setminus A_{r-1}}(X\setminus A_{r-1},E)^G \leq \sum_\alpha \dim \Gamma^\infty(F^\alpha_{Y_j})^G.
\end{displaymath}

Now $B_r=\mathrm{cl}(B_r) \setminus A_{r-1}=A_r \setminus A_{r-1}$, so we may apply inequality \eqref{eqn:second_observation} with $C=A_{r-1}$, $D=A_r$ to conclude the inner induction on $r$. Then, taking $r=k$ concludes the outer induction. Taking $Y=Z$ now concludes the proof of the first statement.

The more general statement of the lemma follows from the first by applying inequality \eqref{eqn:second_observation} with $C=Z_1$, $D=Z_2$.
\endproof

\begin{Corollary}
If for all $\alpha \geq0$ and $y \in Z$ we have $(F_{Gy}^\alpha|_y)^{\Stab(y)}=0$, then $E$ admits no non-trivial $G$-invariant generalized sections supported in $Z$.  
\end{Corollary}

\subsection{Wavefront sets of invariant sections}
We make use of the following general description of the wavefront of a generalized section invariant under the action of a group.  It is surely well-known, however we were not able to find appropriate references in the literature. 

We start with a precise description of the wavefront set of a pull-back by submersion.
\begin{Lemma}\label{lem:submersion_wavefront}
	Let $\pi:X\to Y$ be a submersion between smooth manifolds, and $\mu\in \Gamma^{-\infty}(Y, E)$ for a vector bundle $E$ over $Y$. Then $\WF(\pi ^*\mu)=\pi^*\WF(\mu)$.
\end{Lemma}

\proof
The inclusion $\WF(\pi ^*\mu)\subset\pi^*\WF(\mu)$ is a standard statement appearing in the literature, see \cite[Theorem 8.2.4]{hoermander_pde1}. Let us show that $\pi^*\WF(\mu)\subset\WF(\pi ^*\mu)$. 

We may assume $E$ is the trivial line bundle. Fix $x\in X$ and $y=\pi(x)$. Choose coordinates $(z_1,...,z_n)$ in $U_x \ni x$ and $(z_1,\dots, z_k)$ in $U_y\ni y$ such that $\pi(z_1,...,z_n)=(z_1,\dots, z_k)$. For $\xi=\sum_{j=1}^k \xi_j dz_j\in T^*_yY$ we have $\eta=\pi^*\xi=\sum_{j=1}^k \xi_j dz_j\in T_x^*X$. 

If $\eta \not \in\WF(\pi^*\mu)$, we may find $\epsilon>0$ such that for all $\psi_1\in C^\infty_c((-\epsilon, \epsilon)^k)$, $\psi_2\in C^\infty_c((-\epsilon, \epsilon)^{n-k})$, it holds for all natural $N$ that as $t\to\infty$,
\begin{displaymath}
\int_{\R^n}\mu(z_1,\dots,z_k)\prod_{j=1}^n \exp(it \eta'_j z_j)\psi_1(z_1,\dots,z_k)\psi_2(z_{k+1},\dots,z_n)dz_1\dots dz_n=O(t^{-N})
\end{displaymath}
uniformly in a neighborhood $(\eta'_j)_{j=1}^n\in V_\eta$ of $\eta$. Taking $\eta'_{k+1}=\dots=\eta'_n=0$ we conclude that
\begin{displaymath} 
\int_{\R^k}\mu(z_1,\dots,z_k)\prod_{j=1}^k \exp(it \xi'_j z_j)\psi_1(z_1,\dots,z_k)dz_1\dots dz_k=O(t^{-N})
\end{displaymath}
uniformly in a neighborhood $(\xi'_j)_{j=1}^k\in V_\xi$ of $\xi$. This implies $\xi\notin\WF(\mu)$, as required.
\endproof

\begin{Lemma}\label{lem:invariant_wavefront}
	Let $X$ be a smooth manifold, $G$ a compact Lie group acting on $X$, and $Z\subset X$ a $G$-orbit. Let $E$ be a $G$-bundle over $X$, and $\mu\in\Gamma^{-\infty}(X, E)^G$ an invariant generalized section. Then for $z\in Z$ it holds that $\WF(\mu)\cap T^*_zX\subset N_z^*Z$. 
\end{Lemma}
\proof
 Note that $Z$ is an embedded submanifold since $G$ is compact. Let $H$ be the stabilizer of $z$ in $G$, and $V=N_zZ=T_zX/T_zZ$. The normal bundle $NZ$ is $G$-equivariantly diffeomorphic to the quotient $G\times_H V$ of $G\times V$ by the action of $H$ given for $h\in H$, $g\in G$, $v\in V$ by $h(g,v)=(gh^{-1}, hv)$. By the equivariant slice theorem \cite[Theorem 2.4.1]{duistermaat_kolk}, there is a $G$-equivariant map $A: G\times_H V \to X$ which induces a diffeomorphism $\alpha:U\to W$ between an open $G$-invariant neighborhood $U$ of the zero section and an open $G$-invariant neighborhood $W$ of $Z$ and restricts to $(g, 0)\mapsto gz$ on the zero section. We may then restrict $\mu$ to $W$ and consider $\alpha^*\mu\in\Gamma^{-\infty}(U, A^*E)^G$.
 
 Denote by $\pi_H: G\times V\to G\times_H V$, $\pi_V:G\times V\to V$ the natural projections, and set $\tilde E:=\pi_H^* \alpha^*E$. Since $G$ acts freely on $G\times V$, we may choose a bundle $E_V$ over $V$ such that $\tilde E=\pi_V^* E_V$ (e.g. by taking $E_V=i_0^*\tilde E$ where $i_0(v)=(e,v)\in G\times V$). Then $\pi_V^*:\Gamma^{-\infty}(V, E_V)\to \Gamma^{-\infty}(G\times V, \tilde E)^G$ is onto. To see this, fix a Lebesgue measure $dv$ on $V$ and a Haar measure $dg$ on $G$. Fix $\phi\in \Gamma^\infty_c(V, E_V^*)$ and let $w(g)\in C^\infty_c(G)$ be arbitrary. For $\tilde \mu\in \Gamma^{-\infty}(G\times V, \tilde E)^G$ we find $\langle \tilde \mu, w(g)\phi(v)dgdv\rangle$ is a continuous $G$-invariant functional of $w$, and hence $\langle \tilde \mu, w(g)\phi(v)dgdv\rangle=\nu(\phi dv)\int_Gwdg$ for some $\nu\in \Gamma^{-\infty}(V, E_V)$, which readily implies $\tilde \mu =\pi_V^*\nu$. 
 
Hence, the $G$-invariance of $\alpha^*\mu$ implies that $\tilde \mu:=\pi_H^*\alpha^*\mu=\pi_V^*\nu$ for some $\nu$.

 Writing $T^*(G\times V)=G\times V\times\mathfrak g^*\times V^*$, it then follows that $\WF(\tilde \mu)\subset G\times V\times \{0\}\times V^*$. Moreover, it follows from Lemma \ref{lem:submersion_wavefront} that $\WF(\tilde \mu)=\pi_H^*\alpha^*\WF(\mu)$. Now if $\xi\in \WF(\mu)\cap T_z^*Z$ and $v\in T_zZ$, choose $u\in\mathfrak g$ with $v=\left.\frac{d}{dt}\right|_{t=0}\exp(tu)z$. Then $\langle \xi, v\rangle=\langle \pi_H^*\alpha^*\xi, u\rangle=0$ since $\pi_H^*\xi \in \WF(\tilde \mu)$. Thus $\xi\in N_z^*Z$ as claimed.
 
\endproof

\begin{Corollary}\label{cor:invariant_wavefront}
	Let $G$ be a real semisimple Lie group, $E\to X$ a $G$-equivariant vector bundle, $\mu\in\Gamma^{-\infty}(X, E)^G$, and $Z\subset X$ a $G$-orbit. Then for $z\in Z$ it holds that $\WF(\mu)\cap T^*_zX\subset N_z^*Z$.
\end{Corollary}
\begin{Remark}
	If $Z$ fails to be a submanifold at $z$, $N_z^*Z$ has to be interpreted as the annihilator of the subspace of $T_zX$ spanned by the infinitesimal action of $G$.
\end{Remark}
\proof
	 Let $B$ denote the Killing form on $\mathfrak g$, the Lie algebra of $G$. By assumption, $B$ is non-degenerate. Let $\mathfrak g=\mathfrak k+\mathfrak p$ be a Cartan decomposition, which is then $B$-orthogonal: $\mathfrak p=\mathfrak k^B$. 
	 
We claim that $\Span\left(\cup_{g\in G}\Ad_g(\mathfrak k)\right)=\mathfrak g$. Indeed, assume that the subspaces $\Ad_p(\mathfrak k)$, $p\in \exp(\mathfrak p)$ are $B$-orthogonal to some element $b\in \mathfrak g$. Then $b\in (\Ad_p(\mathfrak k))^B= \Ad_p (\mathfrak p) \Rightarrow \Ad_p b\in  \mathfrak p$ for all $p\in \exp(\mathfrak p)$, in particular $b\in\mathfrak p$. Deriving with respect to $p$ yields $[h, b]\in \mathfrak p$ for all $h\in\mathfrak p$. But since $[\mathfrak p, \mathfrak p]\subset \mathfrak k$, we conclude that $[h,b]=0$  $\forall h\in\mathfrak p$. Then for all $k\in\mathfrak k$ and $h\in\mathfrak p$ we have $B([k, b], h)=B(k, [b,h])=0$, so that $[k, b]\in \mathfrak p^B=\mathfrak k$. However, $[\mathfrak k, \mathfrak p]\subset \mathfrak p$, so we conclude that $[k, b]=0$ $\forall k\in\mathfrak k$. Thus $[b,\mathfrak g]=0$. As $\mathfrak g$ is semisimple, $b=0$, proving the claim.
	 
It follows that the tangent spaces at $z$ to the $K$-orbits through $z$, as $K$ ranges over all possible maximal compact subgroups of $G$, span $\Im(\mathfrak g\to T_zX)$ (which is $T_zZ$ if $Z$ is a manifold at $z$). By Lemma \ref{lem:invariant_wavefront}, $\WF(\mu)\cap T^*_zX\subset \cap_{K\subset G} N_z^*(K z)=N_z^*Z$. 
\endproof

\end{appendix}
\def\cprime{$'$}


\end{document}